\documentclass{elsart}
\usepackage{graphicx}
\usepackage{latexsym}
\usepackage{amssymb}
\usepackage{subfigure}
\newcommand{\PM}{Petviashvili }

\usepackage{graphicx,epsfig,amssymb,amsmath,amsfonts}
\journal{Communications in Nonlinear Science and Numerical Simulation}

\begin{document}
\begin{frontmatter}
\title{An extended Petviashvili method for the numerical generation of traveling and localized waves}
\author{J. \'Alvarez}
\address{Department of Applied Mathematics,
University of Valladolid, Paseo del Cauce 59, 47011,
Valladolid, Spain.}
\address{
IMUVA, Institute of Mathematics of University of Valladolid; Spain.
Email: joralv@eii.uva.es}
\author{A. Dur\'an \thanksref{au}}
\address{Department of Applied Mathematics, University of
Valladolid, Paseo de Belen 15, 47011-Valladolid, Spain.}
\address{
IMUVA, Institute of Mathematics of University of Valladolid; Spain.
Email:
angel@mac.uva.es }
\thanks[au]{Corresponding author}

\begin{abstract}
A family of fixed-point iterations is proposed for the numerical computation of traveling waves and localized ground states. The methods are extended versions of Petviashvili type, and they are applicable when the nonlinear term of the system contains homogeneous functions of different degree. The methods are described and applied to several examples of interest, that calibrate their efficiency.
\end{abstract}
\end{frontmatter}
\section{Introduction}
Introduced here is an extended Petviashvili family of methods, suitable for
the numerical approximation of solutions of systems of the form
\begin{eqnarray}
Lx=N(x),\quad x\in \mathbb{R}^{m}, \quad m>1,\label{lab11}
\end{eqnarray}
where $L$ is a nonsingular real $m\times m$
matrix and $N:\mathbb{R}^{m}\rightarrow \mathbb{R}^{m}$ is a nonlinear function that consists of homogeneous functions with several degrees. The paper will formulate the methods, analyze conditions for the convergence and explore their application in some examples of generation of solitary waves.

The iterative techniques presented here are somehow related to the so-called Petviashvili method, \cite{petviashvili}. This method is usually applied for the numerical resolution of systems of the form (\ref{lab11}), but when $N$ is homogeneous with degree $p$ such that $|p|>1$. It is formulated in the following form: The so-called {\em stabilizing factor}
\begin{eqnarray}
m(x)=\frac{\langle Lx_{n},x_{n}\rangle}{\langle N(x_{n}),x_{n}\rangle},\label{Stb}
\end{eqnarray}
is defined. Then, given an initial iteration
$x_{0}$, the step $n\mapsto n+1$ is implemented as
\begin{eqnarray}
  Lx_{n+1}=m(x_{n})^{\gamma}N(x_{n}), \quad n=0,1,\ldots\label{lab13}
\end{eqnarray}
where $\langle\cdot,\cdot\rangle$ denotes the Euclidean inner
product and $\gamma$ is a free parameter.

The Petviashvili method (\ref{lab13}) is a fixed-point algorithm that originally
appeared to generate numerically lump solitary wave profiles of
the KP-I equation, \cite{petviashvili}. It is usually included as
one of the methods for the numerical generation of solitary waves, a
family that many other techniques belong to, such as shooting
methods, some variants of Newton's method, \cite{lakoba,yang},
variational procedures, \cite{baod,caliario}, squared operator
methods,  \cite{yangl2} or imaginary-time evolution methods,
\cite{yangl1}. From the original paper \cite{petviashvili}, several
convergence studies and generalizations of the method have been
done, \cite{ablowitzm,lakobay,lakobay2,pelinovskys,yang2}.

As analyzed in \cite{alvarezd1,alvarezd2}, from the point of view of the convergence, the Petviashvili method is a modification of the classical fixed-point algorithm
\begin{eqnarray}
  Lx_{n+1}&=&N(x_{n}), \quad n=0,1,\ldots,\label{lab14}
\end{eqnarray}
which overcomes the harmful directions for which (\ref{lab14}) is not convergent. For systems (\ref{lab11}) with $N$ homogeneous of degree $p$ ($|p|>1$), it turns out that $\lambda=p$ is an eigenvalue of the iteration matrix 
\begin{eqnarray}
S=L^{-1}N^{\prime}(x^{\ast}),\label{lab14b} 
\end{eqnarray}
at a fixed point $x^{\ast}$ with an eigenvector given by $x^{\ast}$. (The prime denotes the Jacobian of $N$.) Then the iteration matrix of (\ref{lab13}) consists of a deflaction that moves this eigenvalue to some below one in magnitude (which is zero for the optimal choice of the parameter $\gamma$) and preserves the rest of the spectrum of $S$ in (\ref{lab14b}). Thus, if $\lambda=p$ is the only eigenvalue of $S$ with modulus greater than one, the Petviashvili method leads to convergence.

Introduced here is a family of iterative techniques, based on the philosophy of the Petviashvili method, that can be applied to systems (\ref{lab11}), but where $N$ is a combination of homogeneous functions with different degree and the iteration matrix (\ref{lab14b}) contains one eigenvalue with modulus greater than one (and, consequently the classical fixed-point algorithm does not converge). This type of systems appears in many contexts of interest, with particular emphasis on the numerical generation of solitary waves. For this kind of systems, and contrary to the homogeneous case, the fixed point is not an eigenvector of (\ref{lab14b}) anymore. However, if (\ref{lab14b}) still contains an eigenvalue with modulus above one, the strategy of the Petviashvili method, in order to reduce the magnitude of the eigenvalues, is still applicable, leading to modified versions of the algorithm. This paper concerns the case of nonlinear terms containing different homogeneities, while the possibility of adapting the idea to general nonlinearities will be a subject of future work.

The structure of the paper is as follows. Section \ref{sec2} is devoted to the description of the methods and the structure of the corresponding iteration matrix at the fixed point. For simplicity, the study will be done for a nonlinearity $N$ in (\ref{lab11}), consisting of two different homogeneous terms. The generalization of these algorithms to more than two homogeneities will be done in the expected way. The form of the iteration matrix allows to design the methods with the goal of transforming (\ref{lab14b}) to get convergence. Several examples to illustrate this are shown in Section \ref{sec3}. They include the generation of ground state solutions in different nonlinear Schr\"{o}dinger (NLS) models, with and without potentials and of solitary wave solutions of extended versions of classical nonlinear dispersive equations in water waves.
\section{An extended version of \PM type methods}
\label{sec2}
\subsection{Formulation}
We assume that in (\ref{lab11}) the nonlinear term can be written as
\begin{eqnarray}
\label{lab21} N(x)=N_{1}(x)+N_{2}(x),
\end{eqnarray}
where for $j=1,2$, $N_{j}(x)$ is an homogeneous function with
degree $p_{j}$ such that $|p_{j}|>1, j=1,2$ and $p_{1}\neq p_{2}$.
The following methods for the numerical approximation of (\ref{lab11}) are proposed. We consider $C^{1}$ functions $s_{j}:\mathbb{R}^{m}\rightarrow\mathbb{R}, j=1,2$, homogeneous of degree $q_{j}, j=1,2$ and generate a sequence of the form
\begin{eqnarray}
  Lx_{n+1}&=&s_{1}(x_{n})N_{1}(x_{n})+s_{2}(x_{n})N_{2}(x_{n}), \quad n=0,1,\ldots\label{lab22}
\end{eqnarray}
In the case that $N$ contains more than two homogeneities
\begin{eqnarray*}
N(x)=\sum_{j=1}^{L} N_{j}(x),
\end{eqnarray*}
with $N_{j}$ and homogeneous function of degree $p_{j}, j=1,\ldots
L$ and $|p_{1}|>|p_{2}|>\cdots >|p_{L}|>1$, then the corresponding formulation
must substitute (\ref{lab22}) by
\begin{eqnarray}
Lx_{n+1}=\sum_{j=1}^{L}s_{j}(x_{n}) N_{j}(x_{n}),\label{lab23}
\end{eqnarray}
for some $C^{1}$ homogeneous functions $s_{j}:\mathbb{R}^{m}\rightarrow\mathbb{R}, j=1,\ldots L$.

\subsection{Choice of the stabilizing factors}
The original idea of the \PM method is somehow present in (\ref{lab22}) and (\ref{lab23}): the functions $s_{j}$ could play the role of stabilizing factors and this would guide their choice. First, it is indeed required that the fixed points of the system
\begin{eqnarray*}
  Lx&=&s_{1}(x)N_{1}(x)+s_{2}(x)N_{2}(x),
\end{eqnarray*}
contain fixed points of (\ref{lab11}). Therefore, if $x^{\ast}$ solves (\ref{lab11}), then $s_{1}(x^{\ast})=s_{2}(x^{\ast})=1$. Inversely, if $x_{n}$, defined by (\ref{lab22}) (or (\ref{lab23})) converges to some $x$, then $x$ must be a solution of (\ref{lab11}).

Note that (\ref{lab22}) can be written as a fixed-point algorithm for the iteration function
\begin{eqnarray}
\label{lab28}
F(x)=s_{1}(x)L^{-1}N_{1}(x)+s_{2}(x)L^{-1}N_{2}(x),
\end{eqnarray}
and the associated iteration matrix at $x^{\ast}$ has the form
\begin{eqnarray}
F^{\prime}(x^{\ast})&=&S+\frac{1}{p_{2}-p_{1}}(p_{2}I-S)x^{\ast}\nabla s_{1}(x^{\ast})^{T}\nonumber\\&&-
\frac{1}{p_{2}-p_{1}}(p_{1}I-S)x^{\ast}\nabla s_{2}(x^{\ast})^{T},\label{lab29}
\end{eqnarray}
where $S$ is the iteration matrix (\ref{lab14b}) of the classical fixed-point algorithm (\ref{lab14}).

Some choices of the homogeneous stabilizing factors $s_{j}$ look natural. Two examples in this sense would be as follows: 
then:
\begin{itemize}
\item[(i)] If we take:
\begin{eqnarray}
s_{1}(x)=s_{2}(x)=m(x)^{\gamma},\label{lab24}
\end{eqnarray}
for some $\gamma$ and where $m$ is the stabilizing factor (\ref{Stb}), then the resulting method (\ref{lab22}) consists of implementing the \PM scheme (\ref{lab13}) with  more general nonlinearity (\ref{lab21}).
\item[(ii)] The choice
\begin{eqnarray}
s_{j}(x)=m(x)^{\gamma_{j}},\label{lab25}
\end{eqnarray}
for some $\gamma_{j}, j=1,2$, can be seen as a generalization of (\ref{lab24}). It also looks natural to think that the choice of the parameters $\gamma_{1}, \gamma_{2}$ should have to do with the degrees of homogeneity $p_{1}$ and $p_{2}$. In this sense, we remind that in the case of the \PM method, the optimal choice of $\gamma$ is $\gamma=p/(p-1)$, where $p$ is the degree of homogeneity of the nonlinear term, \cite{lakobay,lakobay2,pelinovskys}.
\end{itemize}
As mentioned in the Introduction, the \PM method is mainly used in systems (\ref{lab11}) where $N$ is an homogeneous term with degree of magnitude above one. This kind of systems is very frequent in the numerical generation of solitary wave profiles in nonlinear dispersive equations. This is the reason, in our opinion, for the relative popularity of the method in that research area. (The origin of the method is also there, \cite{petviashvili}.) These special systems have the key property that the degree of homogeneity is an eigenvalue of the iteration matrix (\ref{lab14b}) of the classical fixed-point algorithm at the fixed point $x^{\ast}$ and with $x^{\ast}$ itself as eigenvector. Thus, the effect of the \PM method is filtering the eigenspace given by $x^{\ast}$. 

Now, in the case of systems (\ref{lab11}) satisfying (\ref{lab21}), the
presence of the fixed point $x^{*}$ as an eigenvector of the
iteration matrix (\ref{lab14b}) is not guaranteed.
This means that any $p_{1}, p_{2}$ is not necessarily an eigenvalue. The
following example illustrates the typical case. We consider the
generation of localized ground states for nonlinear
Schr\"{o}dinger models of the form
\begin{eqnarray}
  iu_{t}+u_{xx}+F(|u|^{2})u=0,\quad F(|u|^{2})=\alpha
  |u|^{m_{1}}+\beta |u|^{m_{2}},\label{lab26}
  \end{eqnarray}
  where $\alpha, \beta, m_{1}, m_{2}$ are real constants. The physical context where (\ref{lab26}) appears and several
  results of existence of localized ground state solutions
  $u(x,t)=U(x)e^{i\mu t}, U(x)\rightarrow 0, |x|\rightarrow\infty,
  \mu>0$, can be seen in \cite{yang2} and references therein. The
  equation for $U$ is of the form
  \begin{eqnarray}
  -\mu U+U^{\prime\prime}+F(U^{2})U=0.\label{lab27}
  \end{eqnarray}
  Explicit formulas are known in some cases. For example, when
  $m_{1}=\sigma, m_{2}=2\sigma, \alpha, \beta> 0$, we have, \cite{yang2}
  \begin{equation}
  \label{exact}
  \begin{array}{l}
U_{\rm exact}(x)=\left(\frac{A}{B+\cosh(Dx)}\right)^{1/\sigma}\\
A=\frac{(2+\sigma)\beta \mu}{\alpha},\quad
D=\sigma\sqrt{\mu},\quad
B={\rm
sgn}(\alpha)\left(1+\frac{(2+\sigma)^{2}\beta}{(1+\sigma)\alpha^{2}}\mu\right)^{-1/2}.
\end{array}
  \end{equation}
Equation (\ref{lab27}) is now discretized. We consider the corresponding periodic problem on a sufficiently long interval $(-l,l)$ and discretize (\ref{lab27}) with Fourier collocation techniques, \cite{Boyd,Canutohqz}. The discrete system will have the form   
  (\ref{lab11}) with
\begin{equation}
\label{lab27b}
\begin{array}{l}
L=\mu I-D_{h}^{2},\quad
 N(U_{h})=N_{1}(U_{h})+N_{2}(U_{h}),\\
N_{1}(U_{h})=\alpha\left(|U_{h}|.^{m_{1}}\right). U_{h}, \quad N_{2}(U_{h})=\beta \left(|U_{h}|.^{m_{2}}\right). U_{h}.
\end{array}
\end{equation}
In (\ref{lab27b}), $I$ is the $m\times m$ identity matrix, $D_{h}$ is the
pseudospectral differentiation matrix and $U_{h}\in \mathbb{R}^{m}$
stands for an approximation to the values of the exact solution
$U(x_{j})$ at the grid points $x_{j}=-l+jh, h=2l/m$ on $(-l,l)$. The dots in $N$ stand for the Hadamard product of the vectors. This Fourier collocation procedure will be taken as the discretization method for all the experiments in the present paper.

\begin{table}
\begin{center}
\begin{tabular}{|c|c|c|c|}\hline
$S=L^{-1}N^{\prime}(U_{\rm exact})$&$F^{\prime}(U_{\rm exact})$&$F^{\prime}(U_{\rm exact})$&$F^{\prime}(U_{\rm exact})$\\\hline
3.479415E+00&9.999999E-01&9.999999E-01&9.999999E-01\\
9.999999E-01&4.836366E-01&4.808735E-01&4.833482E-01\\
4.841875E-01&2.871905E-01&3.840844E-01&2.871905E-01\\
2.871905E-01&-2.383647E-01&2.871905E-01&1.904167E-01\\
1.904415E-01&1.904284E-01&1.904659E-01&1.356336E-01\\
1.356336E-01&1.356336E-01&1.356336E-01&1.015442E-01\\\hline
\end{tabular}
\end{center}
\caption{Six largest magnitude eigenvalues of the iteration matrices for (\ref{lab27}) in the cubic-quintic case ($p_{1}=3, p_{2}=5$), at the exact solution (\ref{exact}). They correspond to: (\ref{lab14b}) (first column); 
(\ref{lab28}), (\ref{lab24}) with $\gamma=p_{1}/(p_{1}-1)$ (second column);
(\ref{lab28}), (\ref{lab24}) with $\gamma=p_{2}/(p_{2}-1)$ (third column);
(\ref{lab28}), (\ref{lab25}) with $\gamma_{1}=p_{1}/(p_{1}-1), \gamma_{2}=p_{2}/(p_{2}-1)$ (fourth column). \label{tav_1}}
\end{table}
Table \ref{tav_1} shows the six largest
magnitude eigenvalues of the iteration matrix
(\ref{lab14b}), evaluated at the exact values
$\tilde{U}=(U(x_{0}),\ldots,U(x_{m-1}))$ for the cubic-quintic
case, that is, with $\sigma=2$ ($m_{1}=2,
m_{2}=4$) (first column).
We observe that the degrees of homogeneities in this case,
$p_{j}=m_{j}+1, j=1,2$, do not appear as eigenvalues. Instead,
there exists a dominant, simple eigenvalue $\lambda^{\ast}=3.479415E+00$,
greater than one. The eigenvalue $\lambda=1$ also appears and it
is simple as well. (This is due to the symmetry of (\ref{lab27}), consisting of spatial translations and was explained in \cite{alvarezd3}. Its effect is an {\em orbital convergence}, that is, a convergence to a possible translated profile.) The rest of the spectrum is within the interval $(0,1)$. Thus, the divergence of the classical fixed-point algorithm in this case is only due to the greater than one eigenvalue.

The rest of the columns in 
Table \ref{tav_1} illustrates the effect of the use of methods of the form (\ref{lab22}). The six largest magnitude eigenvalues of the corresponding iteration matrix (\ref{lab29}) at the exact solution are computed for several choices of (\ref{lab28}), namely: (\ref{lab24}) with $\gamma=p_{1}/(p_{1}-1)$ (second column); (\ref{lab24}) with $\gamma=p_{2}/(p_{2}-1)$ (third column) and (\ref{lab25}) with $\gamma_{1}=p_{1}/(p_{1}-1), \gamma_{2}=p_{2}/(p_{2}-1)$ (fourth column). Note that for the three cases, the effect of the new iteration functions is a translation of the spectrum of the corresponding Jacobian that enables the associated fixed-point algorithm to converge (at least in the previously mentioned orbital sense). The eigenvalue $\lambda^{\ast}$ of $S$ is transformed to a new one with magnitude below one and the rest of the spectrum continues to be in the same range. In the experiments of Section \ref{sec3}, the last method will be implemented in all the examples, as a representative of the family (\ref{lab22}) (or its generalization (\ref{lab23})). This does not rule out, however, other several choices. 

\section{Some applications of the methods in solitary wave generation}
\label{sec3}
In this section, some examples of application of the methods (\ref{lab22}) will be shown. They will illustrate two different situations: the case of isolated fixed points and the case of equations with symmetries (where the fixed points are not isolated and they are gathered in orbits). The examples concern problems of wave generation.
\subsection{Equations with symmetries. Example 1}
\label{sec31}
Considered here are two examples of wave generation for equations with symmetries. The first example considers again the equation (\ref{lab26}). We have taken
$m_{1}=\sigma, m_{2}=2\sigma, \sigma=2$ (cubic-quintic case), to compare with the exact solution (\ref{exact}). The numerical experiments in this example have been performed by applying (\ref{lab28}), (\ref{lab25}) with $\gamma_{j}=p_{j}/(p_{j}-1), p_{j}=m_{j}+1, j=1,2,$ on a Fourier collocation discretization of (\ref{lab26}) of the form (\ref{lab27b}).

\begin{figure}[htbp]
\centering
\subfigure[]{
\includegraphics[width=6.5cm]{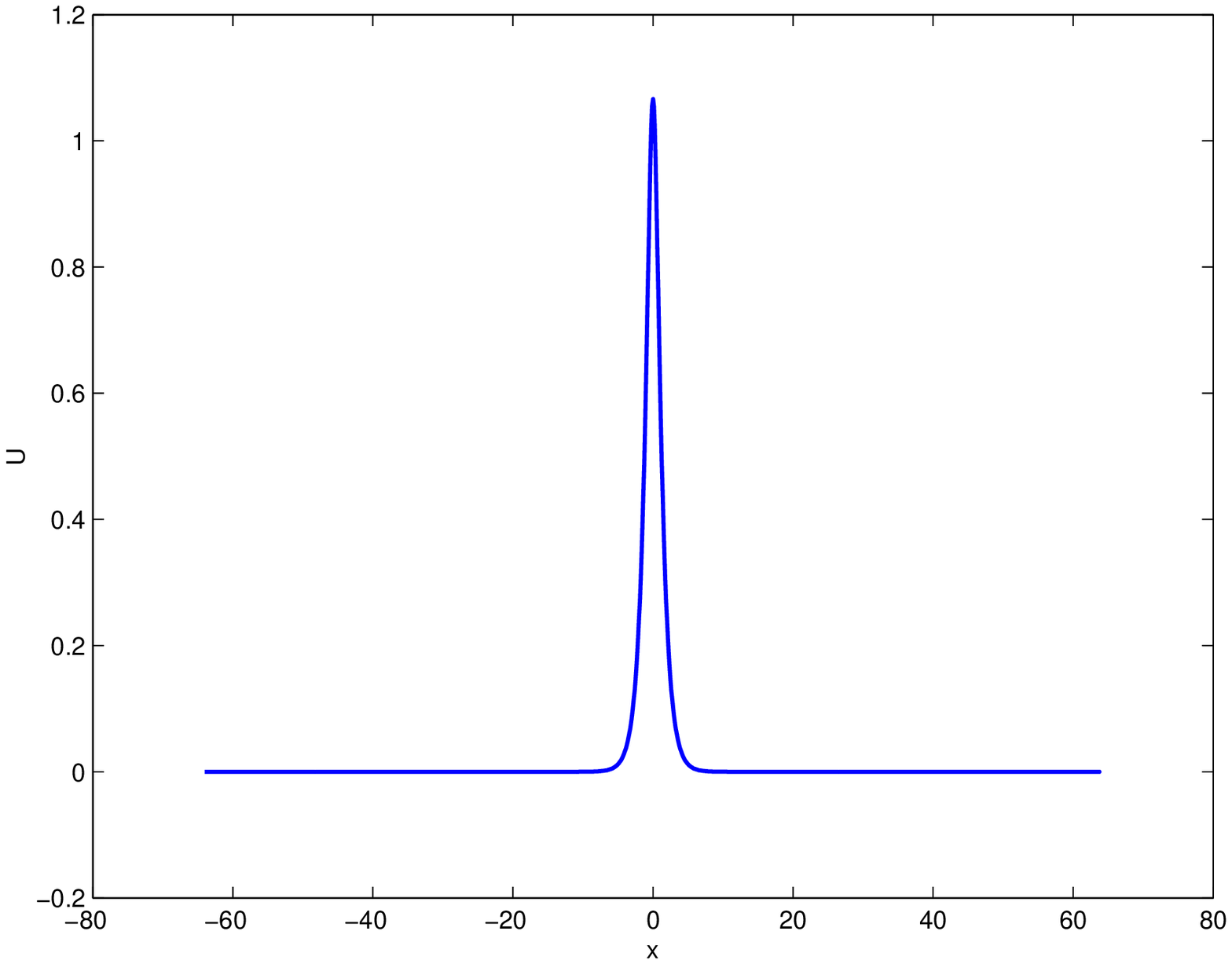} }
\subfigure[]{
\includegraphics[width=6.5cm]{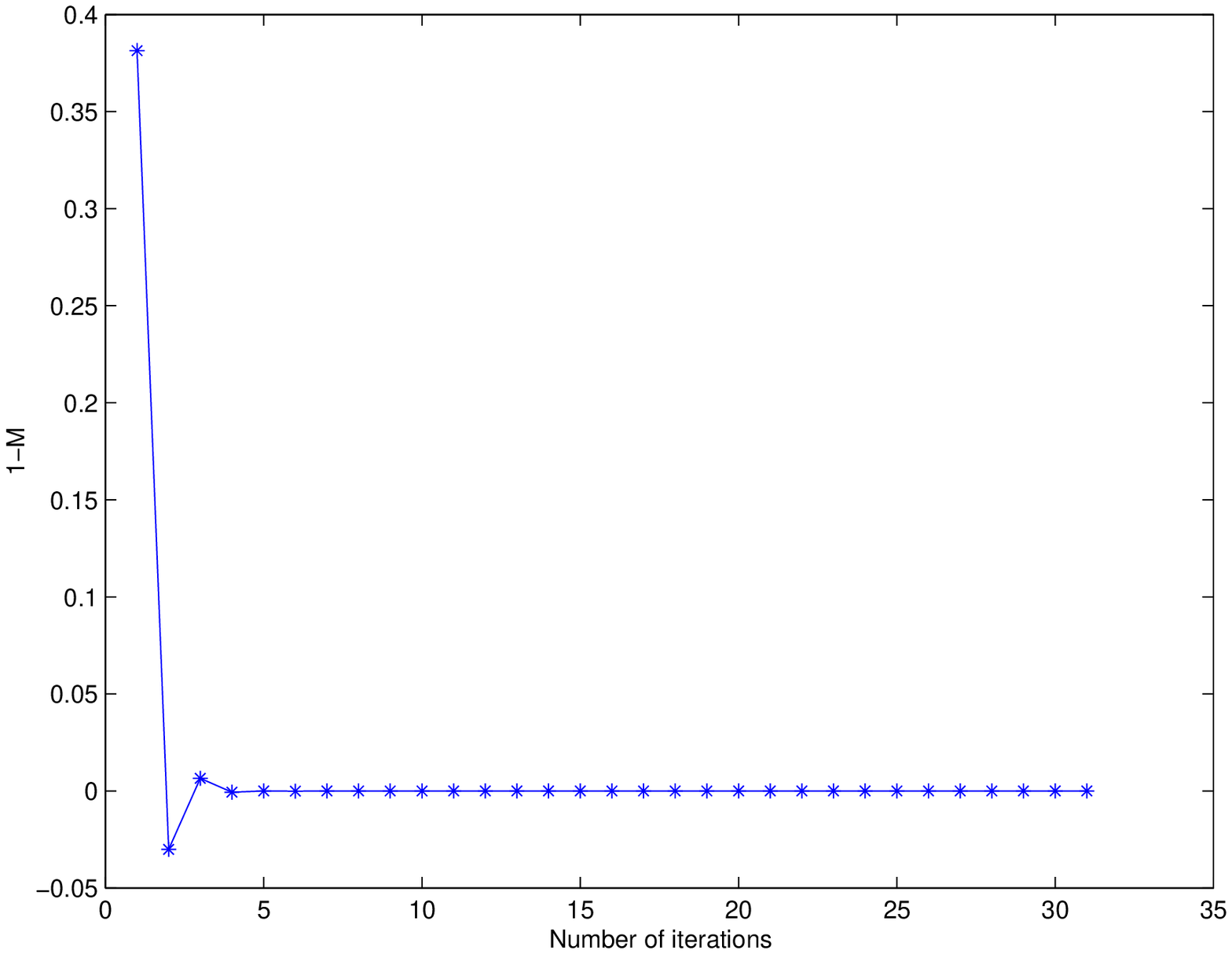} }
\caption{Localized wave solution of (\ref{lab26}) in the
cubic-quintic case ($m_{1}=2, m_{2}=4$), with $\alpha=\beta=1, \mu=1$: (a) Numerical profile obtained with (\ref{lab28}), (\ref{lab25}) and $\gamma_{j}=p_{j}/(p_{j}-1), p_{j}=m_{j}+1, j=1,2$. (b) Discrepancy in the stabilizing
factor (\ref{Stb}) vs number of iterations.} \label{fexample_1a}
\end{figure}

\begin{figure}[htbp]
\centering
\subfigure[]{
\includegraphics[width=6.5cm]{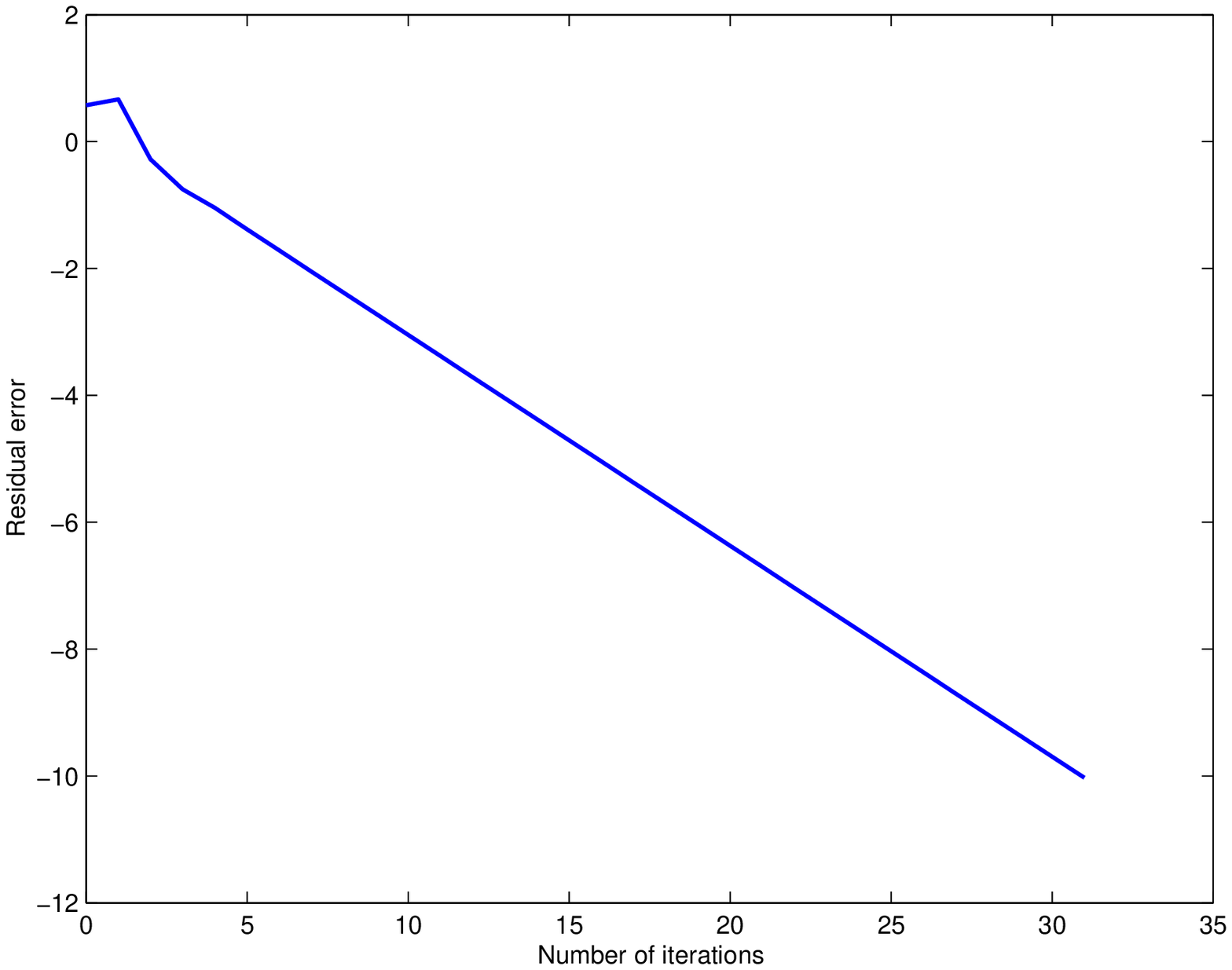} }
\subfigure[]{
\includegraphics[width=6.5cm]{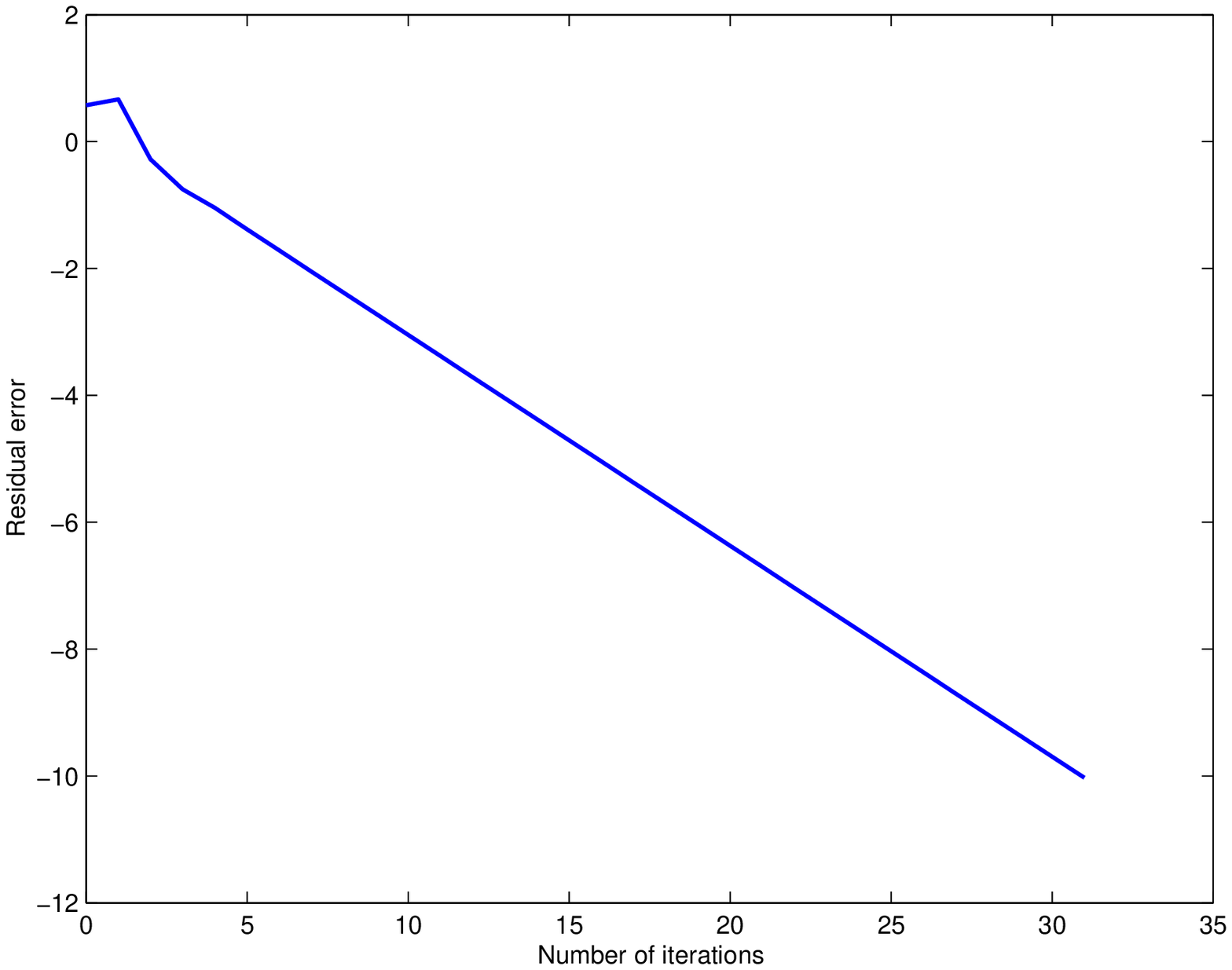} }
\caption{Localized wave solution of (\ref{lab26}) in the
cubic-quintic case ($m_{1}=2, m_{2}=4$) with $\alpha=\beta=1,\mu=1$: (a) Logarithm of the
residual errors (\ref{res}) vs number of iterations. (b) Logarithm of the
errors (\ref{er}) with respect to the exact solution (\ref{exact}) vs number
of iterations.} \label{fexample_1b}
\end{figure}

Figure \ref{fexample_1a}(a) shows the form of the approximated
localized wave. The accuracy of this computed profile is measured by the following results. Observe first that in the case of convergence, the stabilizing factor (\ref{Stb}) evaluated at the iterates generates a sequence that must tend to one. This is observed, for this case, in Figure \ref{fexample_1a}(b). On the other hand,
in Figure \ref{fexample_1b}, two errors, in
semilog scale and with Euclidean norm, are displayed: 
\begin{itemize}
\item[(i)] The (relative) residual error:
\begin{eqnarray}
\label{res}
RE_{n}=||LU_{n}-N(U_{n})||/||\widetilde{U}||,\quad n=0,1,\ldots
\end{eqnarray}
In (\ref{res}), $U_{n}$ is the $n-th$ iterate of the procedure, $L$ and $N$ are given by (\ref{lab27b}).
\item[(ii)] The (relative) error with respect to the exact profile $\widetilde{U}$ at the grid points:
\begin{eqnarray}
\label{er}
RE_{n}=||U_{n}-\widetilde{U}||/||\widetilde{U}||,\quad n=0,1,\ldots.
\end{eqnarray}
\end{itemize}
In
both cases, we obtain an error of order $10^{-10}$ in $30$
iterations, approximately. On the other hand, the order of the
method is linear. This is observed from Table \ref{tav_2}, which displays several ratios between two consecutive values of (\ref{er}).
\begin{table}
\begin{center}
\begin{tabular}{c|c|c|c|c}
$n$&16&18&20&25\\\hline
$||e_{n}||/||e_{n-1}||$&4.650469E-01&4.650471E-01&4.650469E-01&4.650473E-01\\
\end{tabular}
\end{center}
\caption{Localized wave solution of (\ref{lab26}) in the
cubic-quintic case ($m_{1}=2, m_{2}=4$) with $\alpha=\beta=1,\mu=1$. Quotients of several consecutive errors (\ref{er}) with respect to
the exact solution (\ref{exact}). \label{tav_2}}
\end{table}

\begin{figure}[htbp]
\centering
\subfigure[]{
\includegraphics[width=6.5cm]{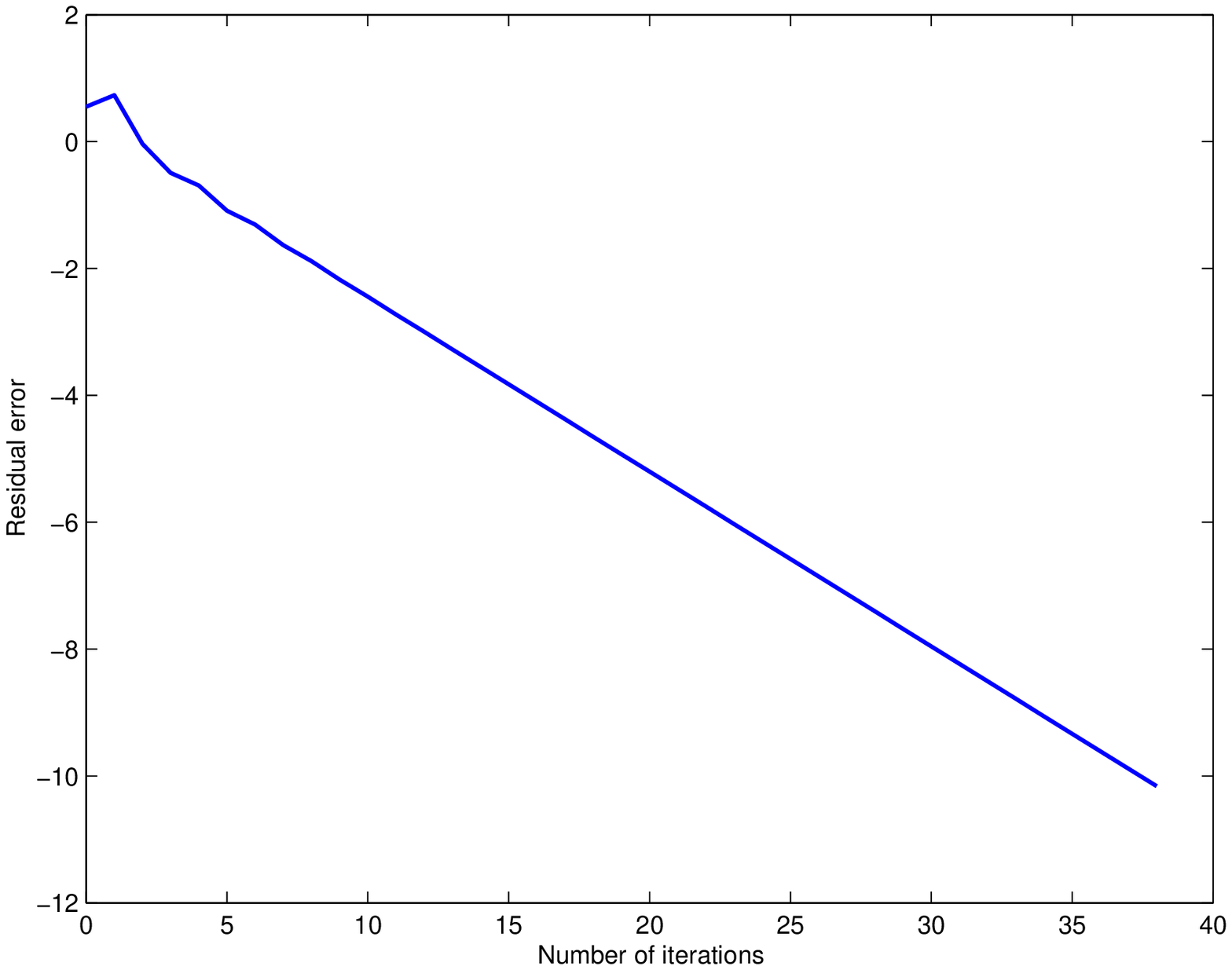} }
\subfigure[]{
\includegraphics[width=6.5cm]{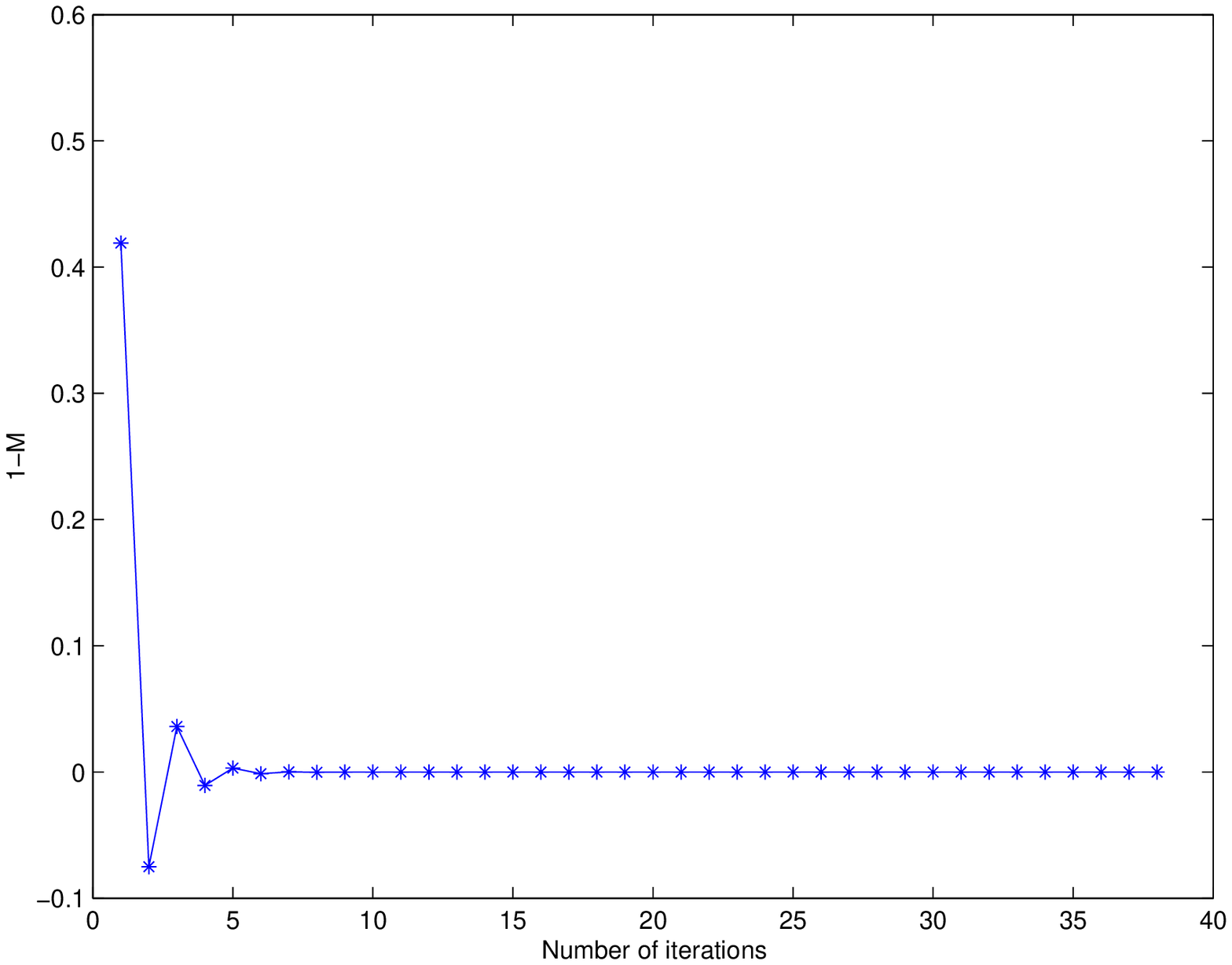} }
\caption{Localized wave solution of (\ref{lab26}) with $m_{1}=1,
m_{2}=3$ and $\alpha=\beta=1,\mu=2\pi$: (a) Logarithm of the residual
errors (\ref{res}) vs number of iterations. (b) Discrepancy in the stabilizing
factor (\ref{Stb}) vs number of iterations.} \label{fexample_2a}
\end{figure}

\begin{figure}[htbp]
\centering
\subfigure[]{
\includegraphics[width=6.5cm]{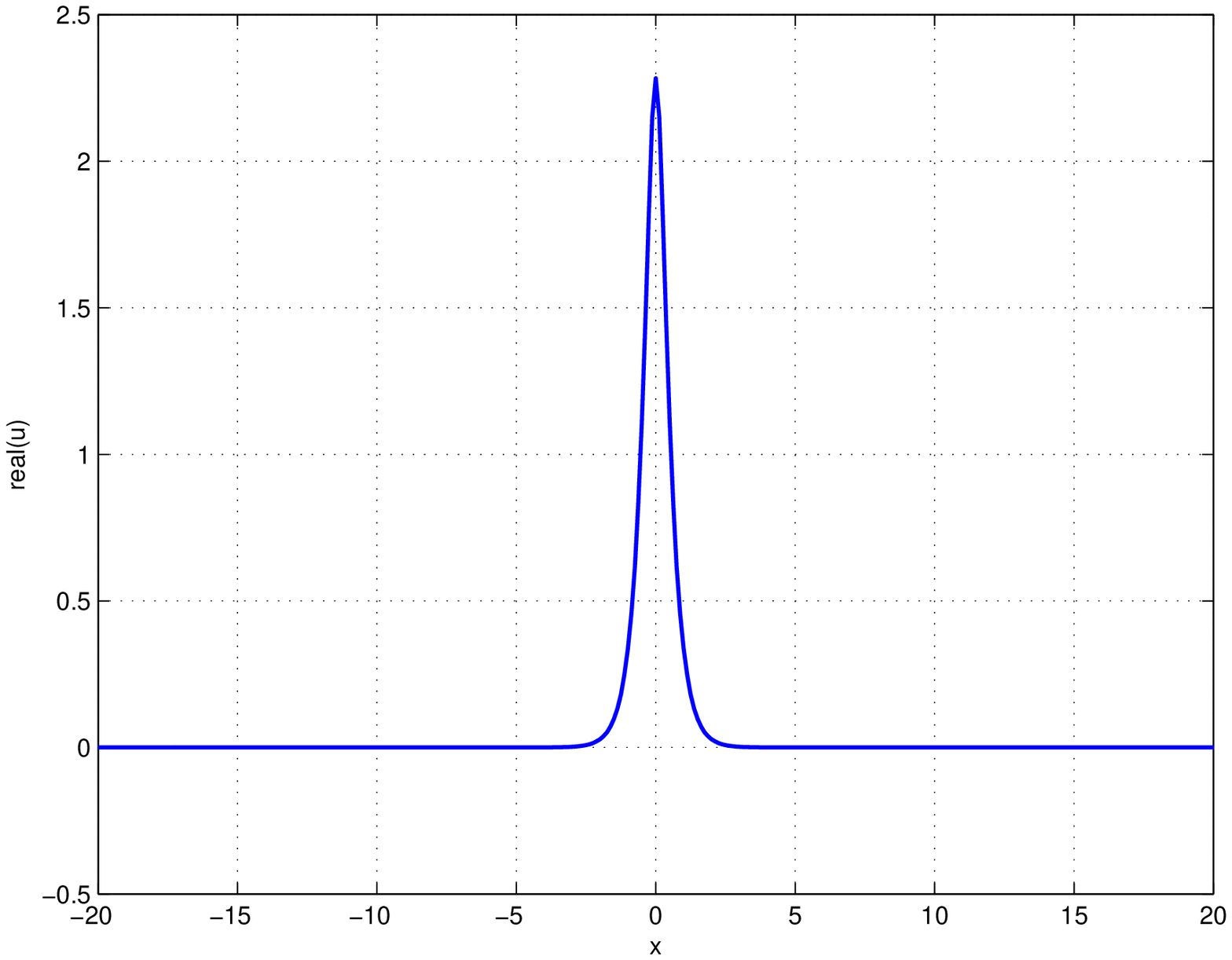} }
\subfigure[]{
\includegraphics[width=6.5cm]{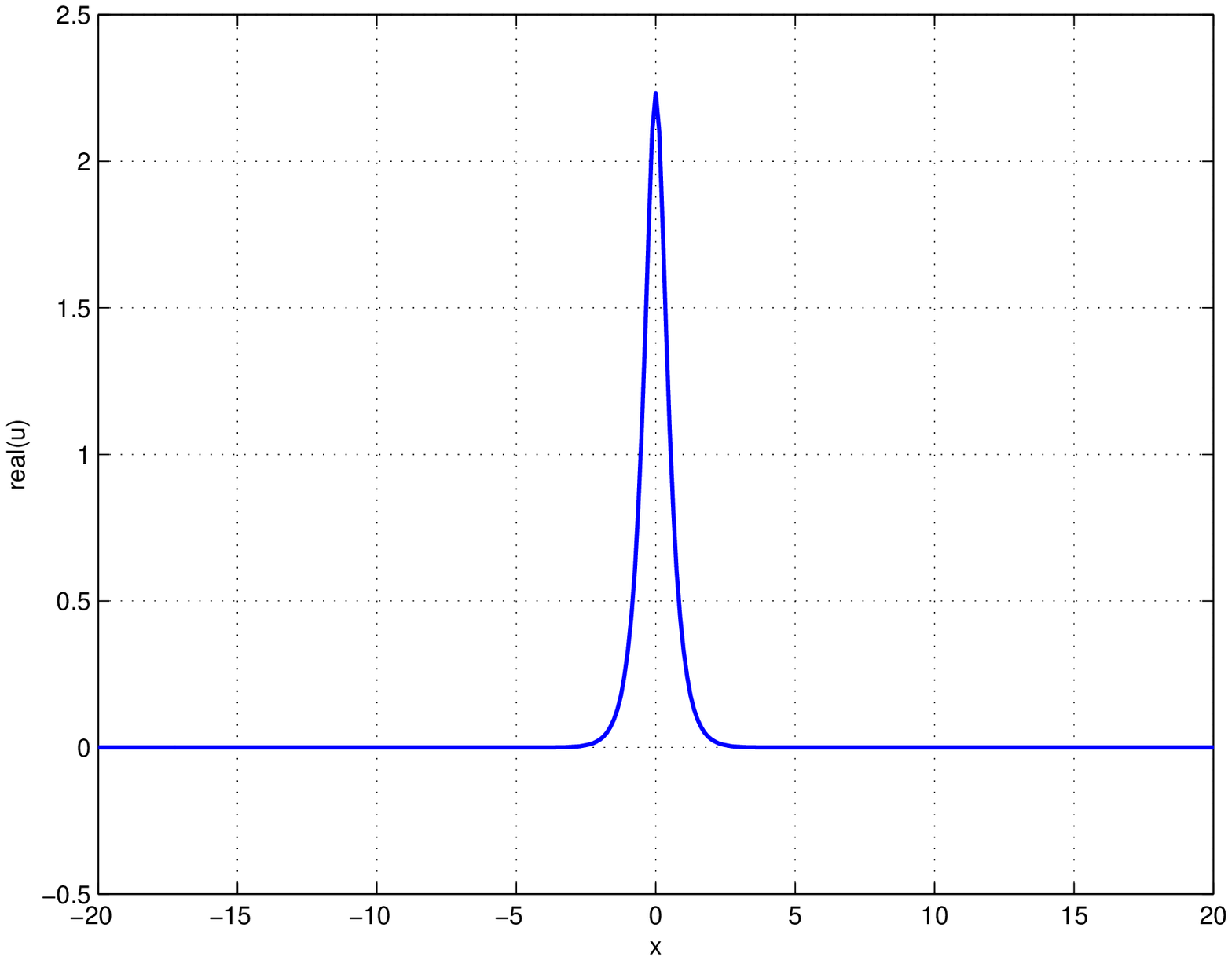} }
\subfigure[]{
\includegraphics[width=6.5cm]{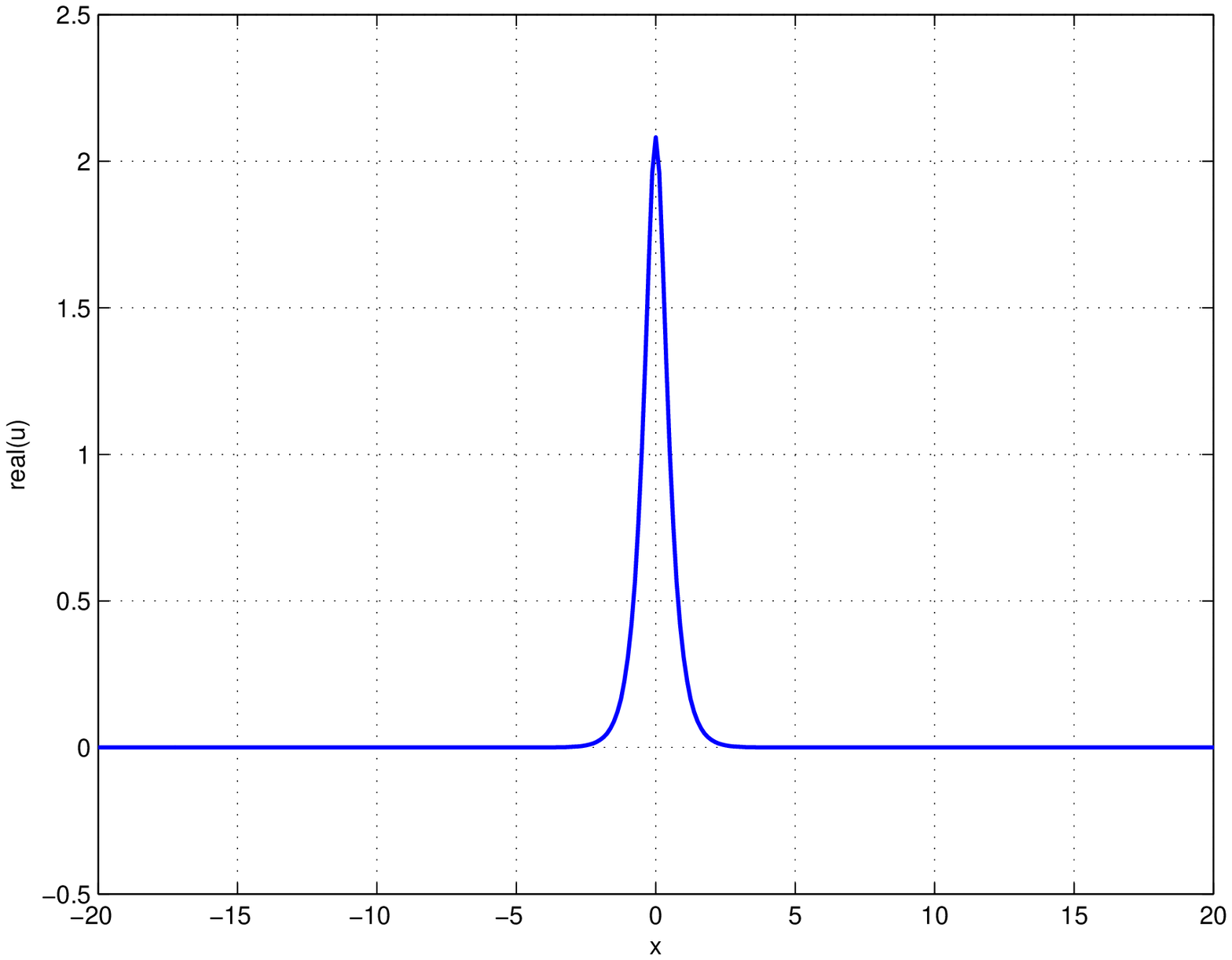} }
\caption{Evolution of the real part of the approximate localized wave 
of (\ref{lab26}) with $m_{1}=1, m_{2}=3$ and
$\alpha=\beta=1,\mu=2\pi$: (a) $t=0$; (b) $t=20$; (c) $t=40$.}
\label{fexample_2b}
\end{figure}

Within this example, we still consider the equations (\ref{lab26}), but now the
parameters are $m_{1}=1, m_{2}=3$, in such a way that the
exact solution is not analytically known. As for the eigenvalues
of the iteration matrix (\ref{lab14b}), the situation is
very similar to that of the example of section \ref{sec2}, compare Tables \ref{tav_1} and \ref{tav_3} (first
column); we have a dominant, simple eigenvalue, greater than one
(in this case, between $p_{2}=2$ and $p_{1}=4$). The next one is the
eigenvalue $\lambda=1$, simple, and the rest is below one. The
convergence for this case is shown in Figure \ref{fexample_2a} and
the evolution of the last computed iterate, as initial condition
of a time stepping code for (\ref{lab26}), is illustrated in
Figure \ref{fexample_2b}, where the real part has been taken. In this case $\mu=2\pi$ and the numerical solution has been displayed at values $t=0, 20, 40$, where $e^{i\mu t}=1$. That is why the profile is approximately the same.
\begin{table}
\begin{center}
\begin{tabular}{|c|c|}
\hline $S=L^{-1}N^{\prime}(U_{f})$ &$F^{\prime}(U_{f})$
\\\hline\hline
3.590177E+00&1.000000E+00\\
9.999999E-01&4.765339E-01\\
4.780777E-01&2.831579E-01\\
2.831579E-01&-2.135475E-01\\
1.886059E-01&1.881397E-01\\
1.353368E-01&1.353368E-01\\
\hline\hline
\end{tabular}
\end{center}
\caption{Six largest magnitude eigenvalues of the 
iteration matrix (\ref{lab14b}) (first column) and  (\ref{lab28}), (\ref{lab25}) with $\gamma_{j}=p_{j}/(p_{j}-1), p_{j}=m_{j}+1, j=1,2$ (second column) for (\ref{lab26}) with $m_{1}=1,
m_{2}=3$ and $\alpha=\beta=\mu=2\pi$. $U_{f}$ stands for the last computed iterate. \label{tav_3}}
\end{table}

\subsection{Equations with symmetries. Example 2}
\label{sec32}
The purpose of the second example is illustrating the performance of the methods (\ref{lab22}) when generating solitary wave profiles of the so-called  e-Boussinesq system
\begin{eqnarray}
\eta_{t}&=&-d_{1}W_{x}-d_{2}W_{xxx}-d_{4}(W\eta)_{x}+d_{5}(W\eta^{2})_{x},\label{ebou1}\\
W_{t}&=&-\frac{1}{d_{1}}\eta_{x}-d_{3}W_{xxt}-\frac{d_{4}}{2}(W^{2})_{x}+d_{5}(W^{2}\eta)_{x},\label{ebou2}
\end{eqnarray}
with
\begin{eqnarray*}
&&d_{1}=\frac{H}{r+H},\quad
d_{2}=\frac{H^{2}}{2(r+H)^{2}}(s+\frac{2}{3}(1+r H)),\nonumber\\
&&d_{3}=\frac{sd_{1}}{2},\quad
d_{4}=\frac{H^{2}-r}{(r+H)^{2}},\quad
d_{5}=\frac{r(1+H)^{2}}{(r+H)^{3}}, \label{ebou3}
\end{eqnarray*}
and some parameters $r, H, s$, with physical meaning. Equations
(\ref{ebou1}), (\ref{ebou2}) are derived in \cite{Nguyend} as a
Boussinesq system for two-way propagation of interfacial waves
under certain physical conditions of the model. Smooth solitary wave
solutions $\eta=\eta(x-c_{s}t), W=W(x-c_{s}t)$, vanishing at
infinity, must satisfy the system
\begin{eqnarray}
c_{s}\eta-(d_{1}+d_{2}\partial_{xx})W&=&-W\eta(-d_{4}+d_{5}\eta),\label{ebou4}\\
-\frac{1}{d_{1}}\eta+c_{s}(1+d_{3}\partial_{xx})W&=&-W^{2}(-\frac{d_{4}}{2}+d_{5}\eta).\label{ebou5}
\end{eqnarray}
\begin{figure}[htbp]
\centering
\subfigure[]{
\includegraphics[width=6.5cm]{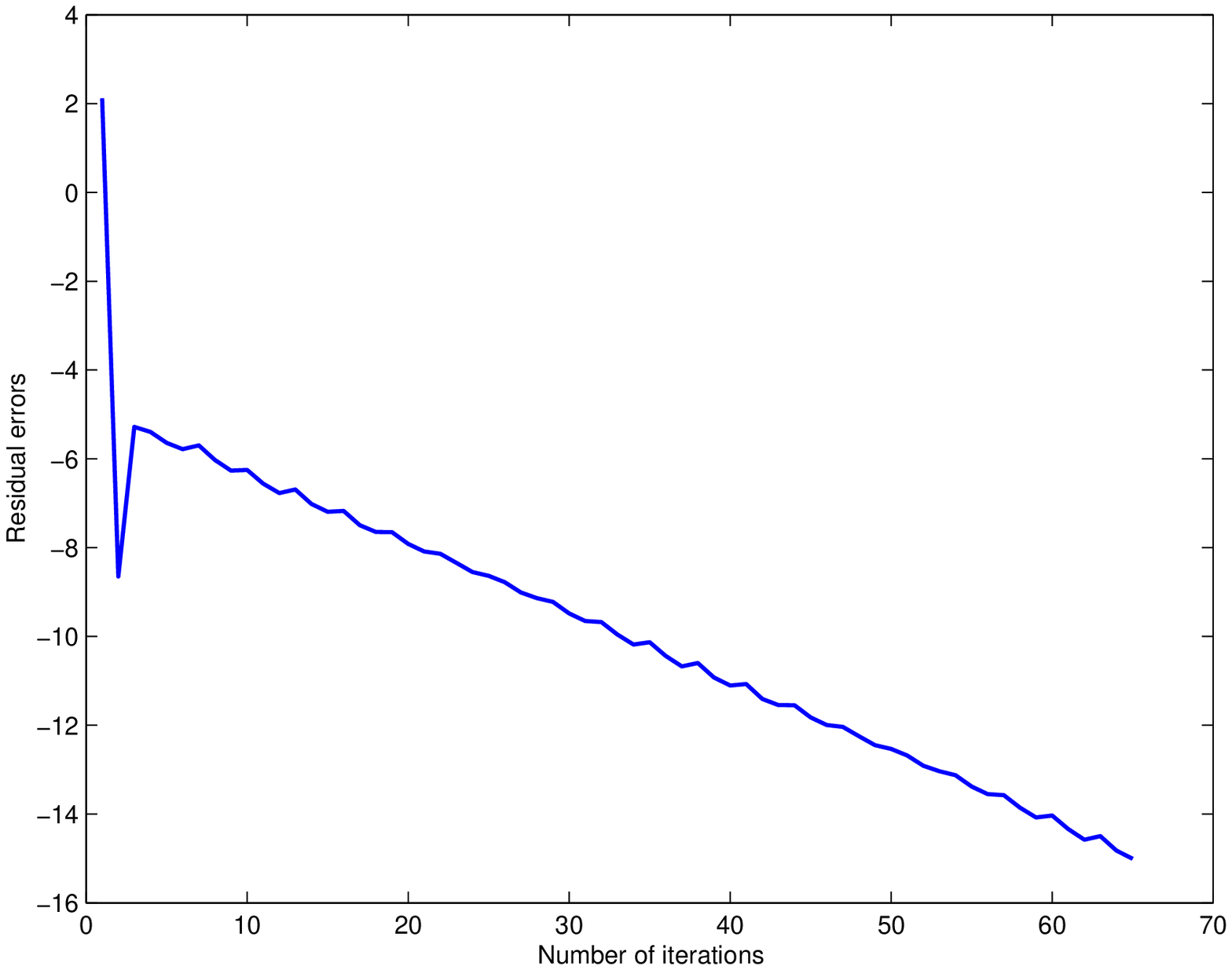} }
\subfigure[]{
\includegraphics[width=6.5cm]{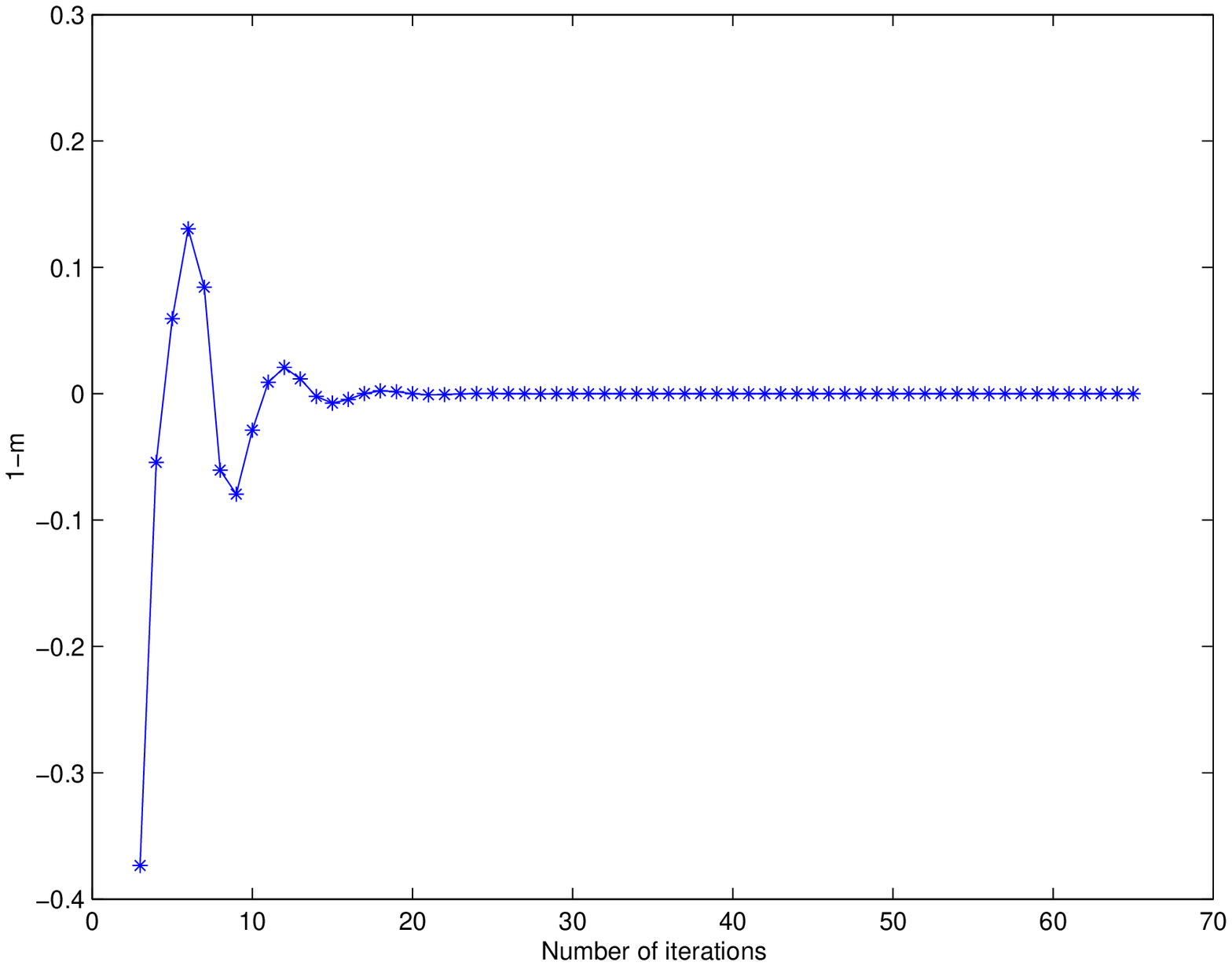} }
\caption{Solitary wave solution of (\ref{ebou1}), (\ref{ebou2}) with  $r=0.8,
H=0.95, c_{s}=1.02$: (a) Logarithm of the residual errors (\ref{res}) vs number of
iterations. (b) Discrepancy in the stabilizing factor (\ref{Stb}) vs number of
iterations.} \label{fexample_3a}
\end{figure}

\begin{figure}[htbp]
\centering
\subfigure[]{
\includegraphics[width=6.5cm]{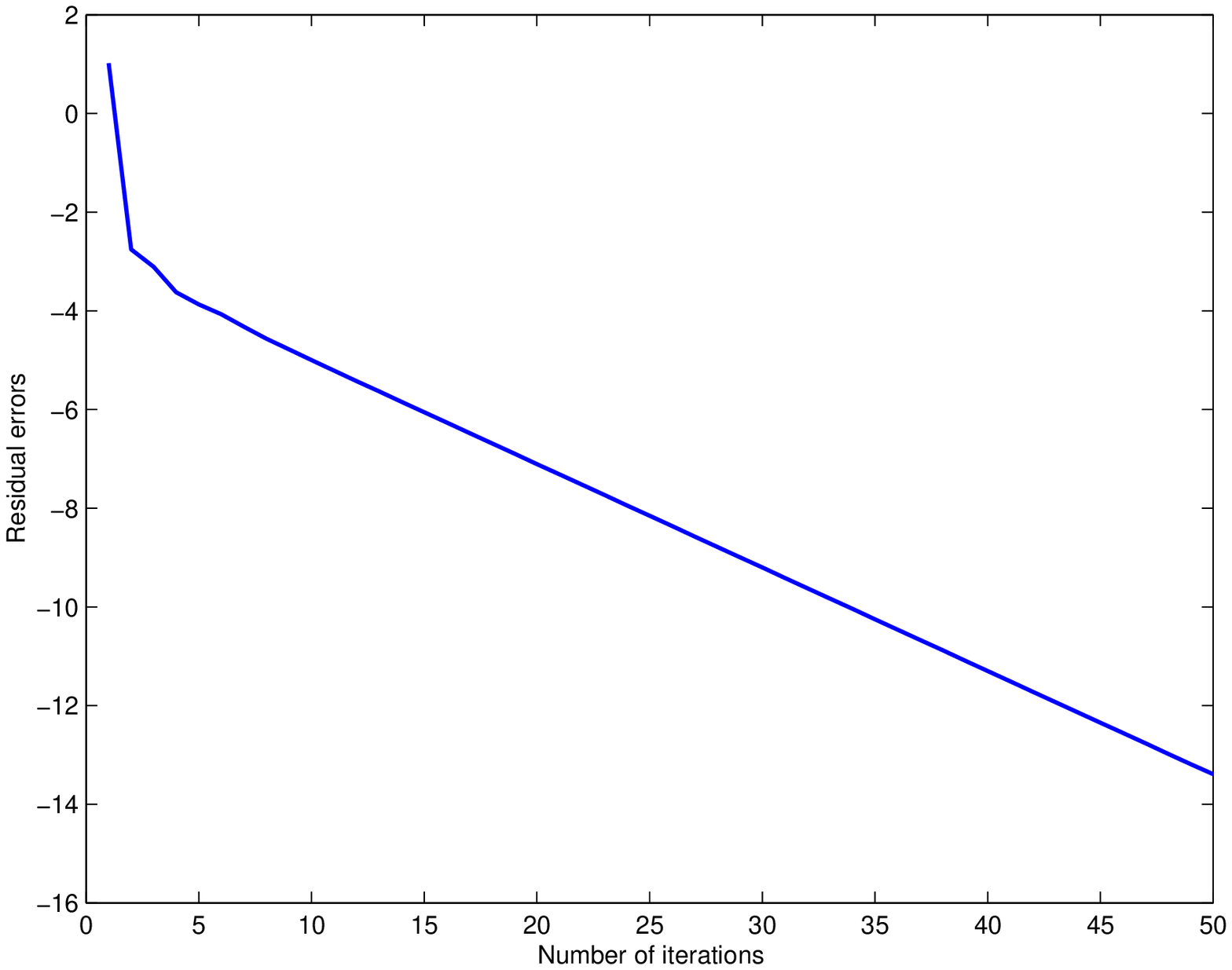} }
\subfigure[]{
\includegraphics[width=6.5cm]{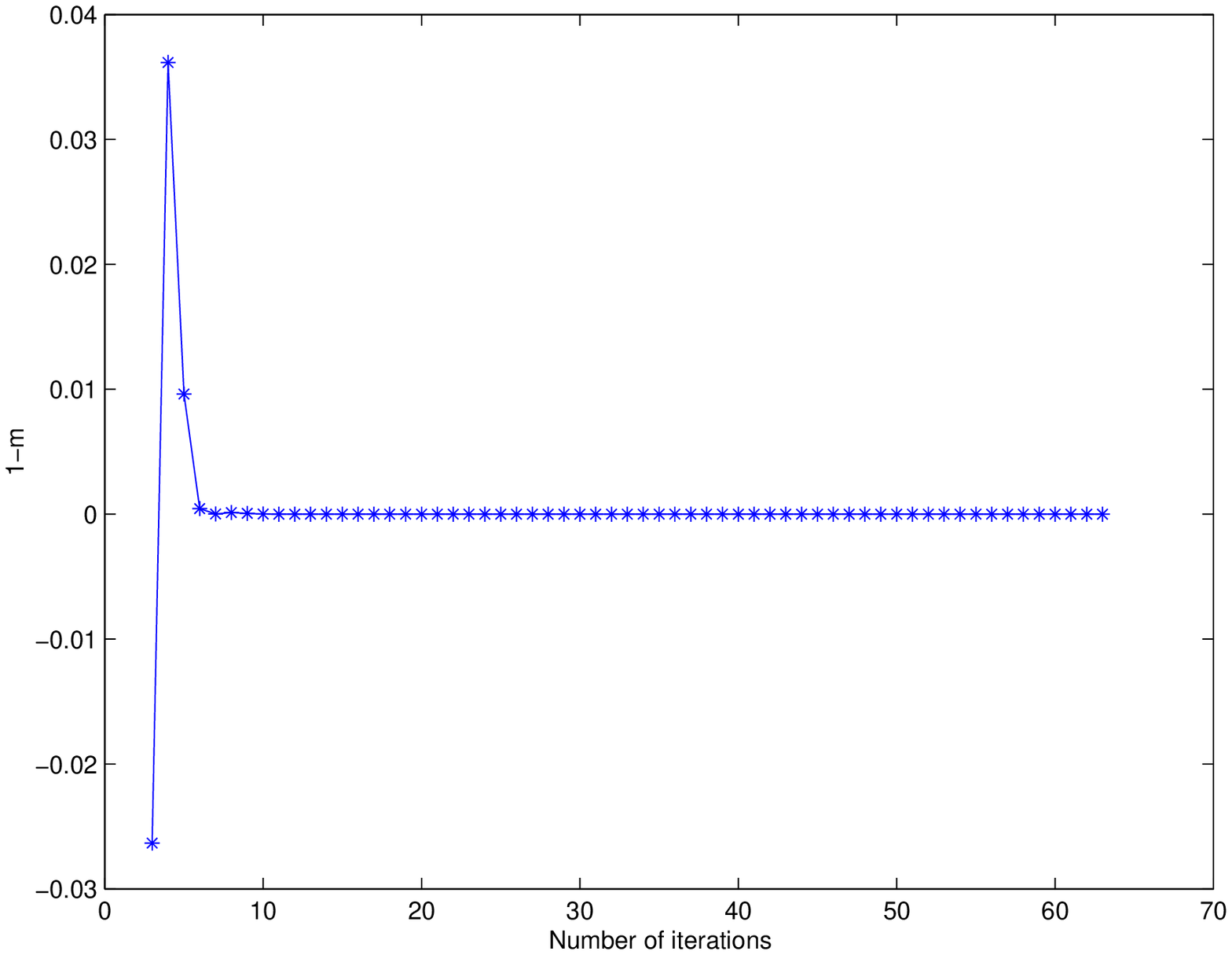} }
\caption{Solitary wave solution of (\ref{ebou1}), (\ref{ebou2}) with  $r=0.8,
H=1.8, c_{s}=1.05$: (a) Logarithm of the residual errors (\ref{res}) vs number of
iterations. (b) Discrepancy in the stabilizing factor (\ref{Stb}) vs number of
iterations.} \label{fexample_3aa}
\end{figure}
\begin{figure}[htbp]
\centering
\subfigure[]{
\includegraphics[width=6.5cm]{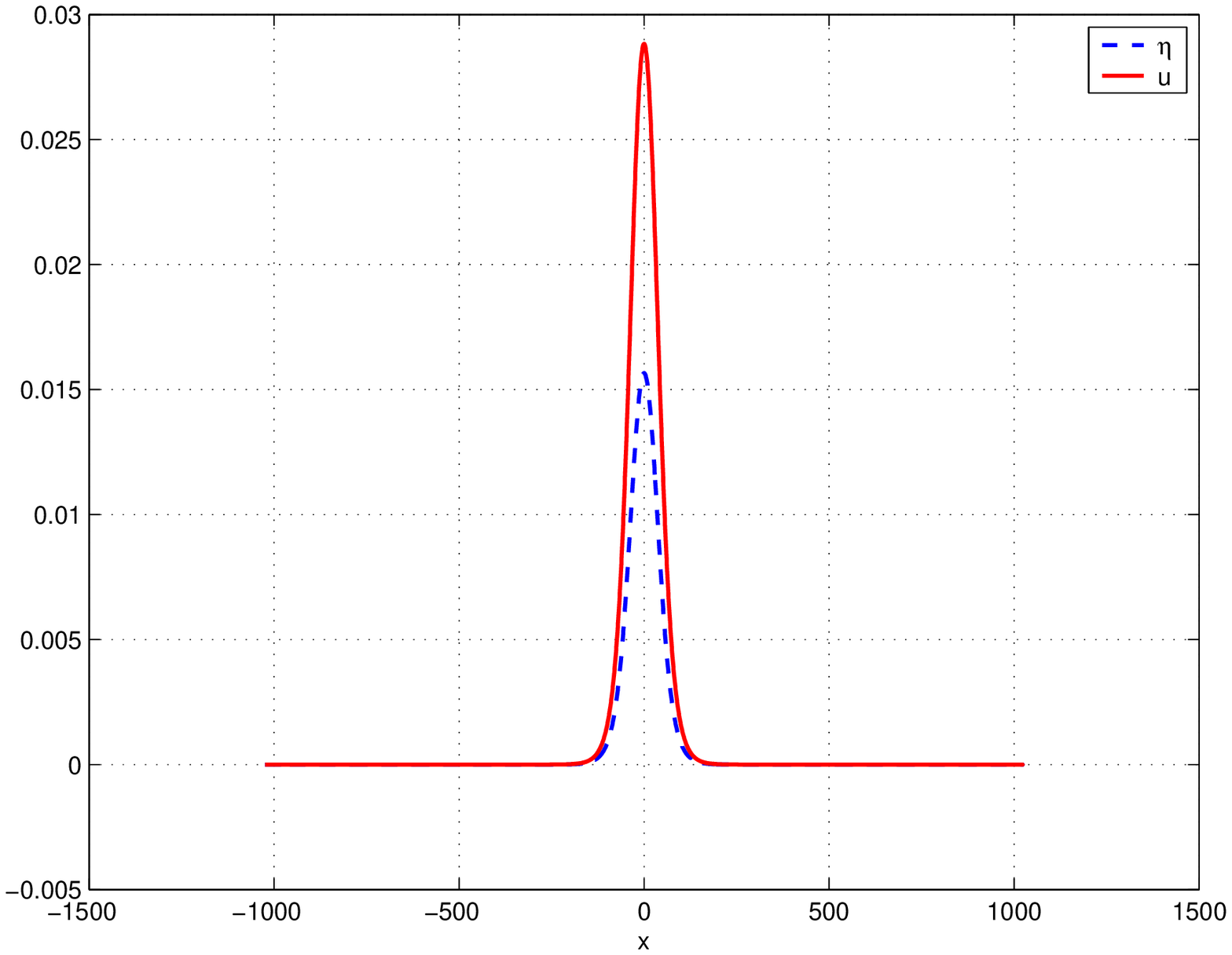}}
\subfigure[]{
\includegraphics[width=6.5cm]{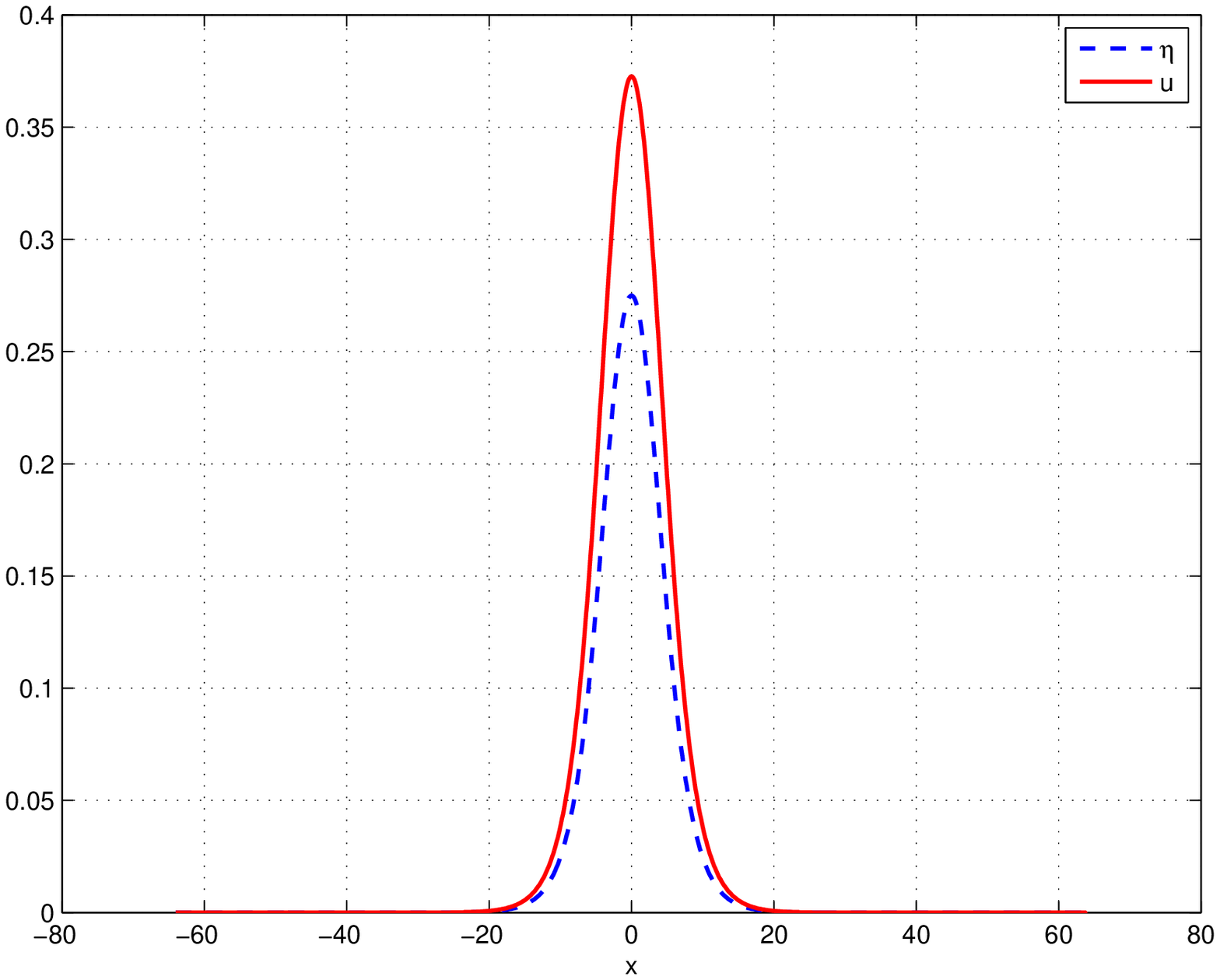}}
\caption{Approximate solitary wave solution of (\ref{ebou1}), (\ref{ebou2}). $u$ profile with solid line and $\eta$ profile with dashed line: (a) $r=0.8,
H=0.95, c_{s}=1.02$; (b) $r=0.8, H=1.8, c_{s}=1.05$.} \label{fexample_3b}
\end{figure}

Again, a Fourier collocation discretization of (\ref{ebou4}),
(\ref{ebou5}) has been considered, leading to the system
\begin{eqnarray}
\label{ebsystemd} \underbrace{ \begin{pmatrix}
    c_{s}I&-(d_{1}I+d_{2}D_{h}^{2} \\
    -\frac{1}{d_{1}}I&c_{s}(I+d_{3}D_{h}^{2})\\
    \end{pmatrix}}_{L}\begin{pmatrix}
    \eta_{h} \\
    W_{h}\\
    \end{pmatrix}
    =-\underbrace{\begin{pmatrix}
    W_{h}.\eta_{h}.(-d_{4}+d_{5}\eta_{h}) \\
    W_{h}.^{2}(-\frac{d_{4}}{2}+d_{5}\eta_{h})\\
    \end{pmatrix}
    }_{N(\eta_{h},W_{h})}.
\end{eqnarray}
In this case, we have a quadratic$+$cubic nonlinearity ($p_{1}=3,
p_{2}=2$).  System (\ref{ebsystemd}) is iteratively solved by using the iteration (\ref{lab22}) with $s_{1}, s_{2}$ given by (\ref{lab25}) and $\gamma_{j}=p_{j}/(p_{j}-1), p_{j}=m_{j}+1, j=1,2$. Two experiments are considered, corresponding to different
values of $r$ and $H$, but with $s=-(1+rH)$,
\cite{Nguyend}. Figures \ref{fexample_3a} and \ref{fexample_3aa}
display the convergence in both cases, as for the residual error (\ref{res}) (with $L$ and $N$ given by (\ref{ebsystemd})), 
and the stabilizing factor (\ref{Stb}). In computational terms, the case
$r=0.8, H=0.95$ is harder, since a larger interval of integration
is considered, see the approximate profiles $\eta, W$ in Figure
\ref{fexample_3b}. As far as the eigenvalues are concerned, for the
case $r=0.8, H=1.8$, the results shown in Table \ref{tav_4} provide similar information (one larger than one
eigenvalue, the eigenvalue one, both simple, and the rest below
one) and, in this case, the dominant eigenvalue is below the
lowest homogeneity $p_{2}=2$. 
\begin{table}
\begin{center}
\begin{tabular}{|c|c|}\hline\hline
$S=L^{-1}N^{\prime}(x_{f})$&$F^{\prime}(x_{f})$\\\hline
1.558592E+00&9.999999E-01\\
9.999999E-01&6.252658E-01\\
6.456383E-01&4.556396E-01\\
4.556396E-01&3.506633E-01-1.562707E-02i\\
3.406177E-01&3.506633E-01+1.562707E-02i\\
2.671308E-01&12.671308E-01\\\hline
\end{tabular}
\end{center}
\caption{Six largest magnitude eigenvalues of the 
iteration matrix (\ref{lab14b}) (first column) and (\ref{lab29}) (second column) for $\gamma_{1}=2, \gamma_{2}=3/2$,  where $x_{f}=(\eta_{f}, u_{f})$ is the last computed iterate, in (\ref{ebou1}), (\ref{ebou2})  with $r=0.8,
H=1.8, c_{s}=1.05$. \label{tav_4}}
\end{table}
\begin{figure}[htbp]
\centering \subfigure[]{
\includegraphics[width=6.5cm]{eb1.eps}}
\subfigure[]{
\includegraphics[width=6.5cm]{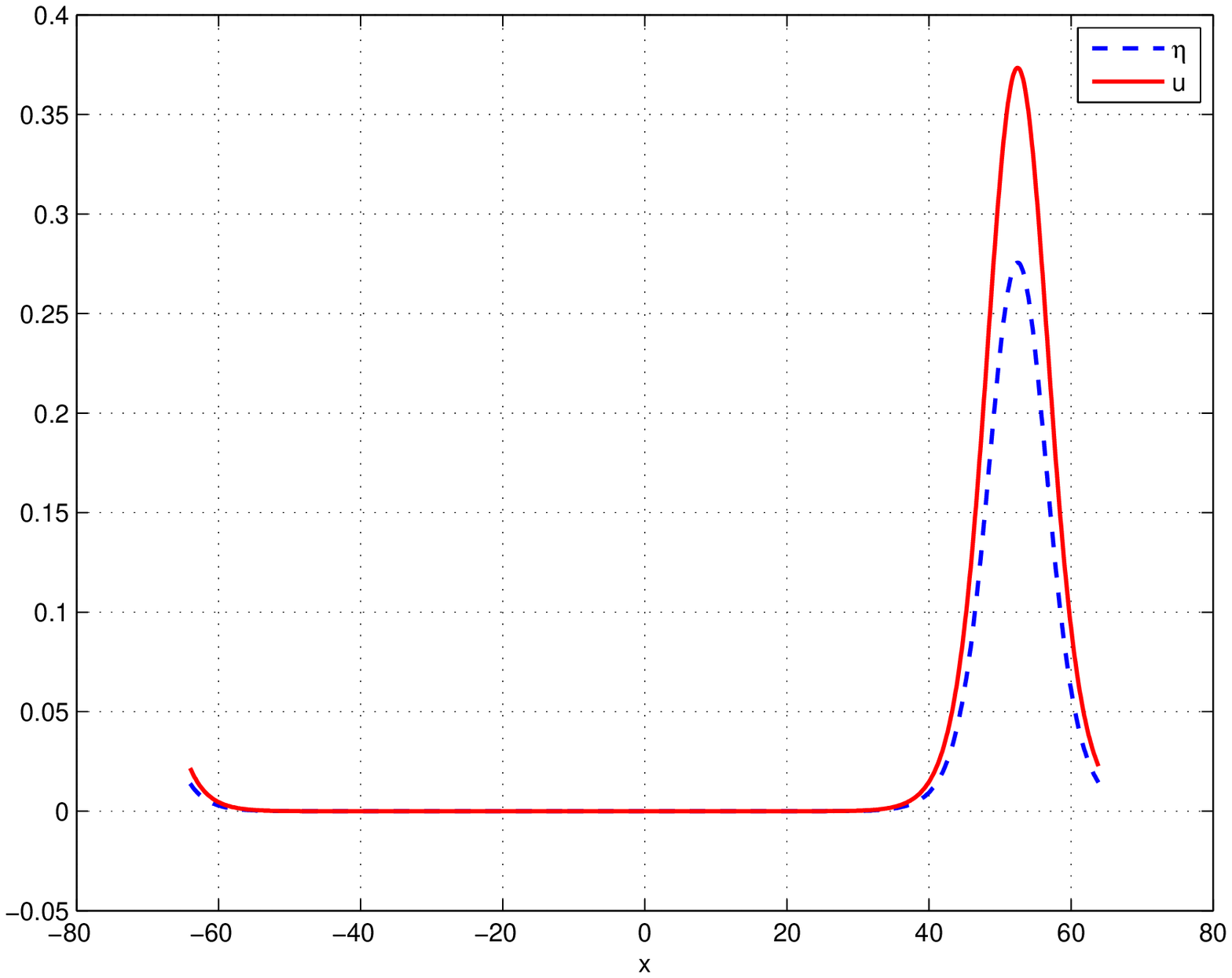}}
\subfigure[]{
\includegraphics[width=6.5cm]{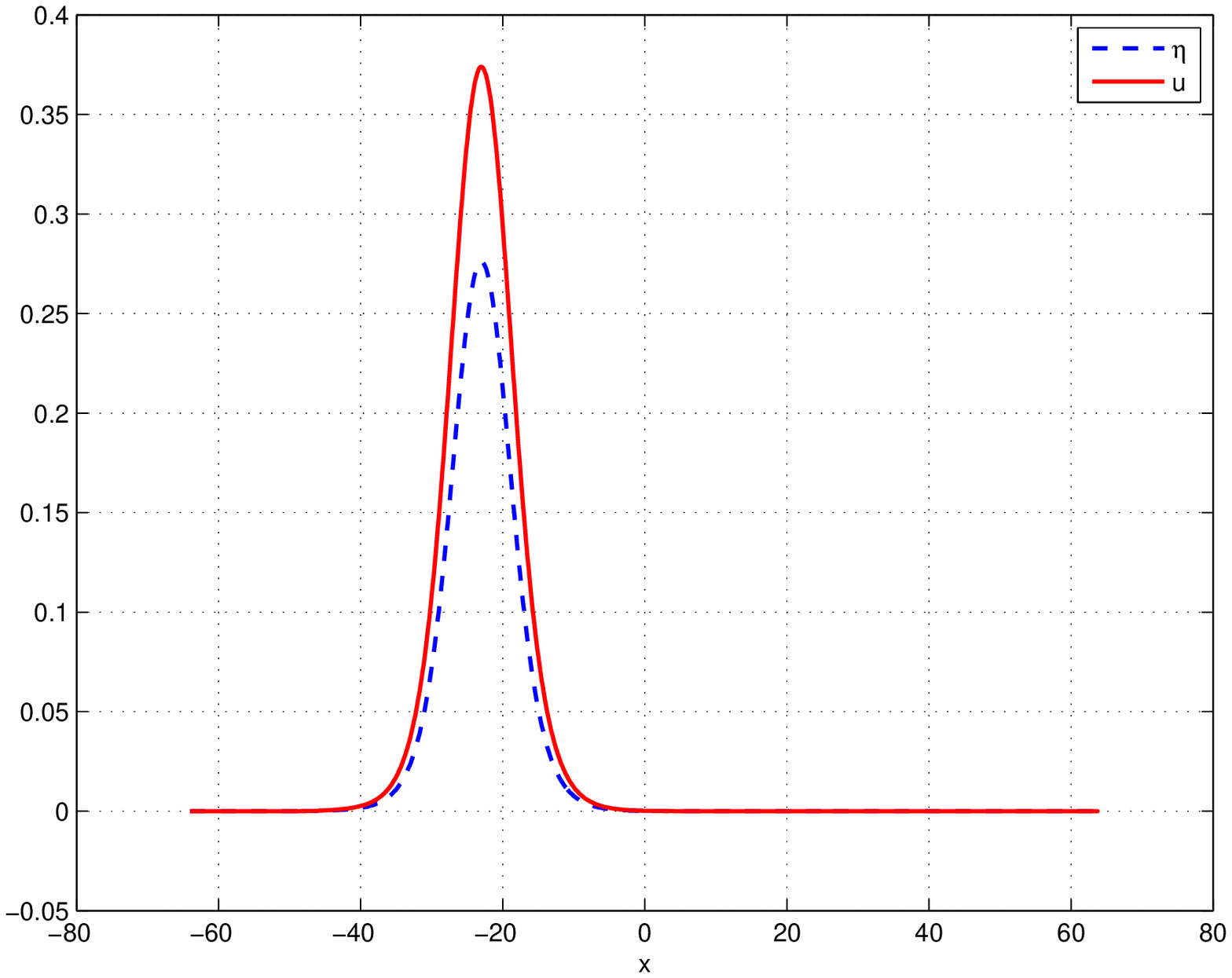}}
\subfigure[]{
\includegraphics[width=6.5cm]{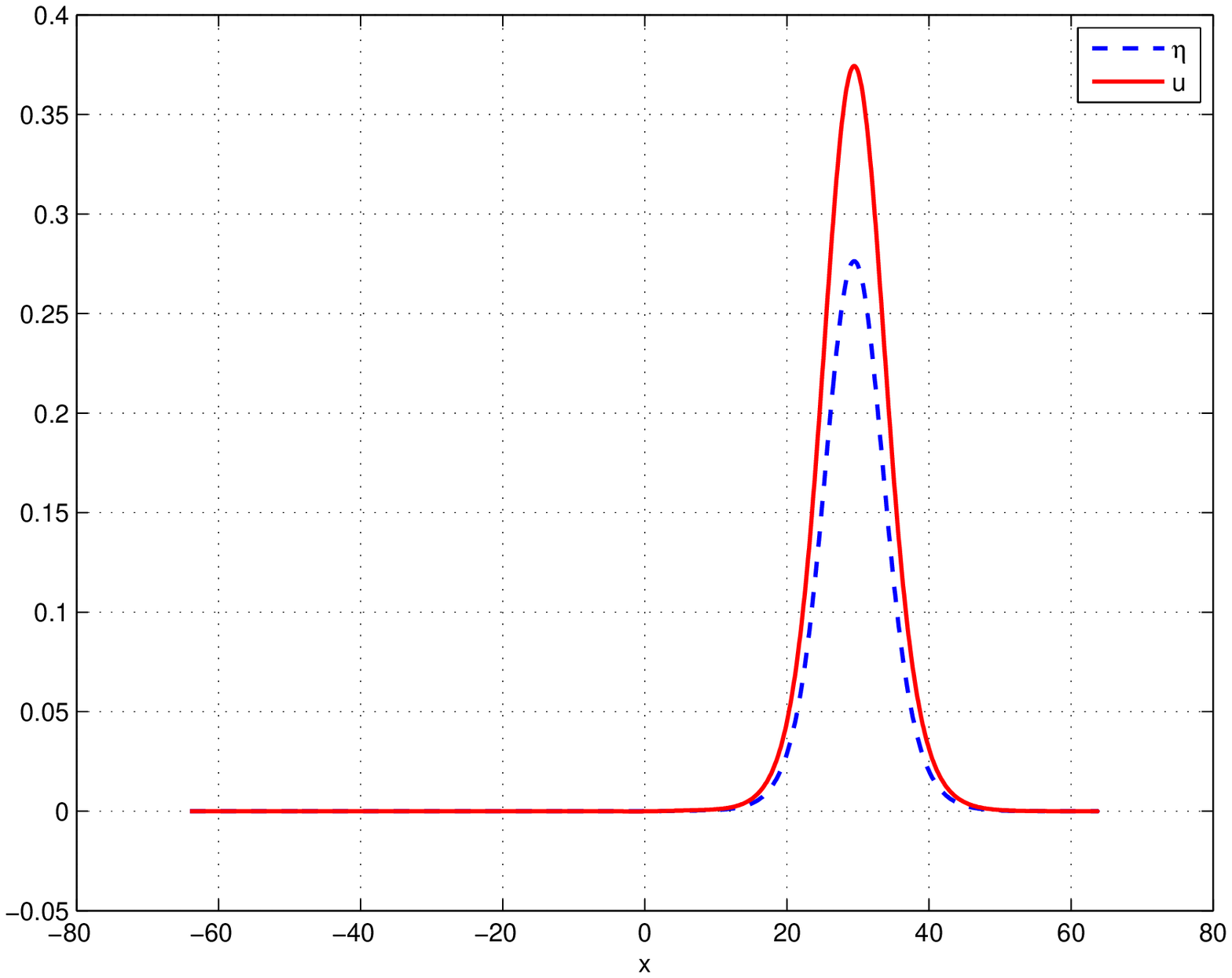}}
\subfigure[]{
\includegraphics[width=6.5cm]{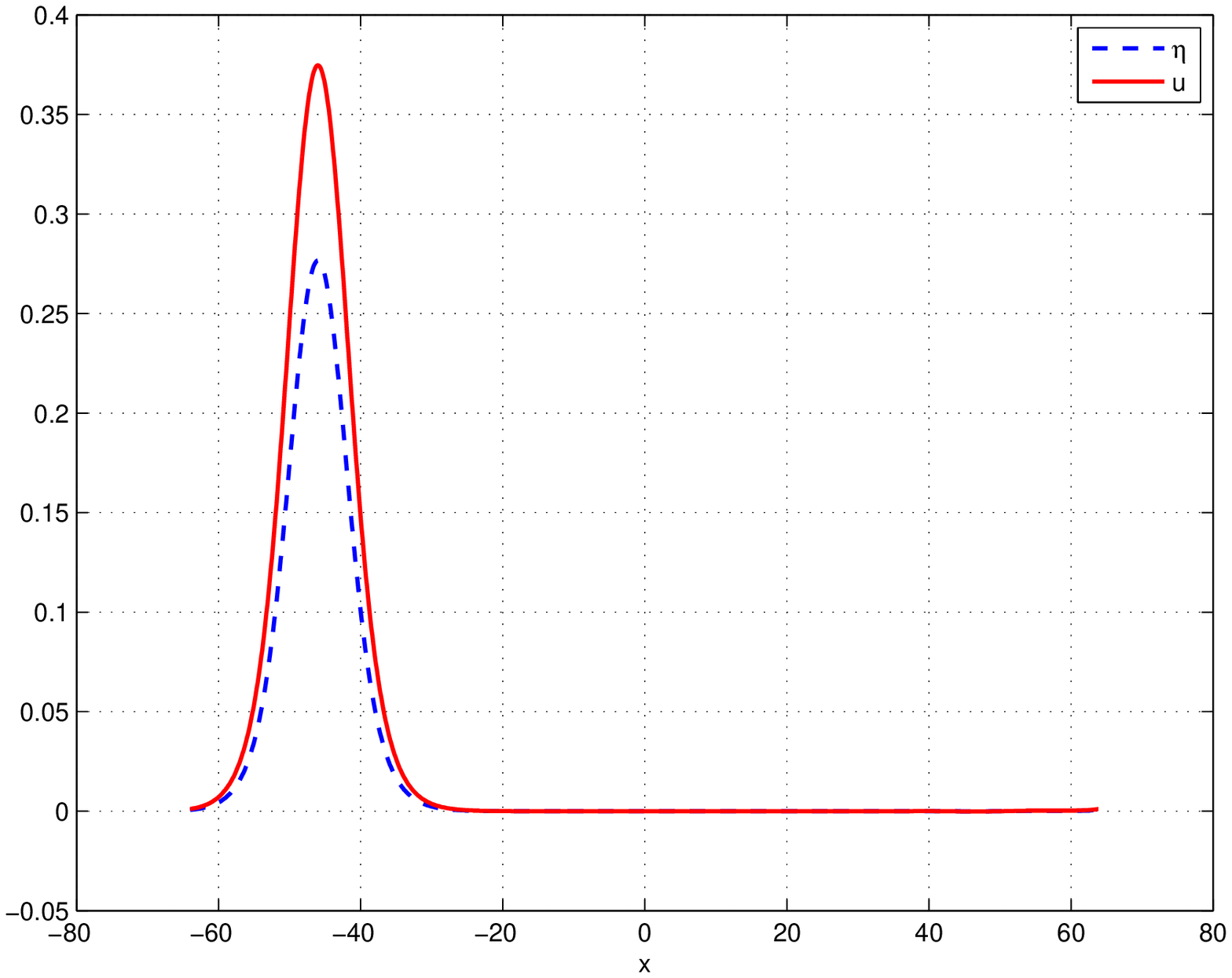}}
\caption{ Evolution of the profiles of Figure \ref{fexample_3b}(b). Approximate $u$ (solid lines) and $\eta$ (dashed lines) at times (a) $t=0$, (b) $t=50$, (c) $t=100$, (d) $t=150$, (e) $t=200$. The speed is $c_{s}=1.05$.} \label{fexample_5}
\end{figure}
Finally, in order to check the accuracy of the profiles, these have been taken as initial conditions of a time-stepping code to integrate (\ref{ebou1}), (\ref{ebou2}) numerically. The evolution of the numerical approximation is illustrated in Figure \ref{fexample_5}. (The periodic boundary conditions forces the numerical solution, traveling to the right, to go out of the computational window and reappear on the left.) We observe that the profiles propagate without any disturbance, behind or in front of them.

\subsection{Equations without symmetries. Example 1}
\label{sec:33}
A second group of experiments concerns the generation of ground state solutions in nonlinear Schr\"{o}dinger (NLS) equations with potentials, in such a way that the presence of that function breaks the symmetry and the localized ground state solutions can be obtained as isolated fixed points of a differential system.
From the point of view of the iteration, this means that the eigenvalue $\lambda=1$, that appeared in the experiments performed in sections \ref{sec31} and \ref{sec32}, will not be present here, \cite{alvarezd3}.

The first example of this group involves the generation of ground state solutions of a generalized NLS equation of the form
\begin{eqnarray}
iu_{t}+u_{xx}-V(x)u+|u|^{2}-\gamma|u|^{4}u=0,\label{lab31}
\end{eqnarray}
where $V(x)$ is a symmetric double-well potential
\begin{eqnarray}
V(x)=-V_{0}\left({\rm sech}^{2}(x+x_{0})+{\rm sech}^{2}(x-x_{0})\right),\label{dwp}
\end{eqnarray}
with $V_{0}>0, x_{0}\in\mathbb{R}$ and $\gamma>0$. Equation (\ref{lab31}) is studied in \cite{yang3}, where bifurcations of solitary waves are analyzed. Localized wave solutions
$u(x,t)=U(x)e^{i\mu t}$ satisfy
\begin{eqnarray}
-\mu U+u^{\prime\prime}-V(x)U+|U|^{2}U-\gamma|U|^{4}U=0.\label{lab31b}
\end{eqnarray}
The system (\ref{lab11}) for the corresponding
Fourier collocation approximation consists of
\begin{equation}
\label{lab32}
\begin{array}{l}
L=\mu I-D_{h}^{2}+diag(V(x_{0}),\ldots,V(x_{m-1})),\\
N(U_{h})=N_{1}(U_{h})+N_{2}(U_{h})=\left(|U_{h}|.^{2}\right).U_{h}-\gamma \left(|U_{h}|.^{4}\right).U_{h},
\end{array}
\end{equation}
where $diag(V(x_{1}),\ldots,V(x_{m}))$ stands for the $m\times m$ diagonal matrix with diagonal entries the values of the potential (\ref{dwp}) at the grid points, $V(x_{j}), j=0,\ldots,m-1$.
In \cite{yang3}, two types of bifurcations are predicted. They can be identified from the power curve of a family of positive, symmetric solitary wave solutions. It is the curve $(\mu, P(\mu))$ where $P$ is the power
\begin{eqnarray}
P(\mu)=\int_{-\infty}^{\infty} U^{2}(x,\mu)dx,\label{power}
\end{eqnarray}
Considered here is the numerical resolution of (\ref{lab11}), (\ref{lab32}) by using (\ref{lab22}) with (\ref{lab25}) for three values of $\mu=1.9, 2.69$. They correspond to values of $\mu$ close to the two types of bifurcation.

The nonlinear term in (\ref{lab32}) contains two homogeneities with degrees $p_{1}=3, p_{2}=5$. For the experiments below, the parameters in (\ref{lab31}), (\ref{dwp}) take the values $V_{0}=2.8, x_{0}=1.5, \gamma=0.25$, \cite{yang3}.

\begin{table}
\begin{center}
\begin{tabular}{|c|c|c|c|}
\hline\hline
\multicolumn{2}{|c|}
{$\mu=1.9$} & \multicolumn{2}{|c|}{$\mu=2.69$}\\
\hline
$S=L^{-1}N^{\prime}(U_{f})$&$F^{\prime}(U_{f}),$&$L^{-1}N^{\prime}(U_{f})$&$F^{\prime}(U_{f})$\\\hline
2.935028E+00&6.686554E-01&1.305101E+00&8.633083E-01\\
6.686554E-01&1.1592453E-01&8.633083E-01&7.559223E-01\\
1.159253E-01&6.877290E-02&4.854323E-01&4.855206E-01\\
6.877290E-02&4.180900E-02&2.575429E-01&2.575429E-01\\
4.187650E-02&2.903208E-02&1.500531E-01&2.314120E-01\\
2.903208E-02&2.118866E-02&1.179335E-01&1.179335E-01\\
\hline
\end{tabular}
\end{center}
\caption{Six largest magnitude eigenvalues of the iteration matrices
(\ref{lab14b}) and (\ref{lab28}), (\ref{lab25}) with $\gamma_{1}=3/2, \gamma_{2}=5/4$, for $\mu=1.9$ (first and second columns) and $\mu=2.69$ (third and fourth columns). $U_{f}$ stands for the last computed iterate. \label{tav_5}}
\end{table}

\begin{figure}[htbp]
\centering \subfigure[]{
\includegraphics[width=6.5cm]{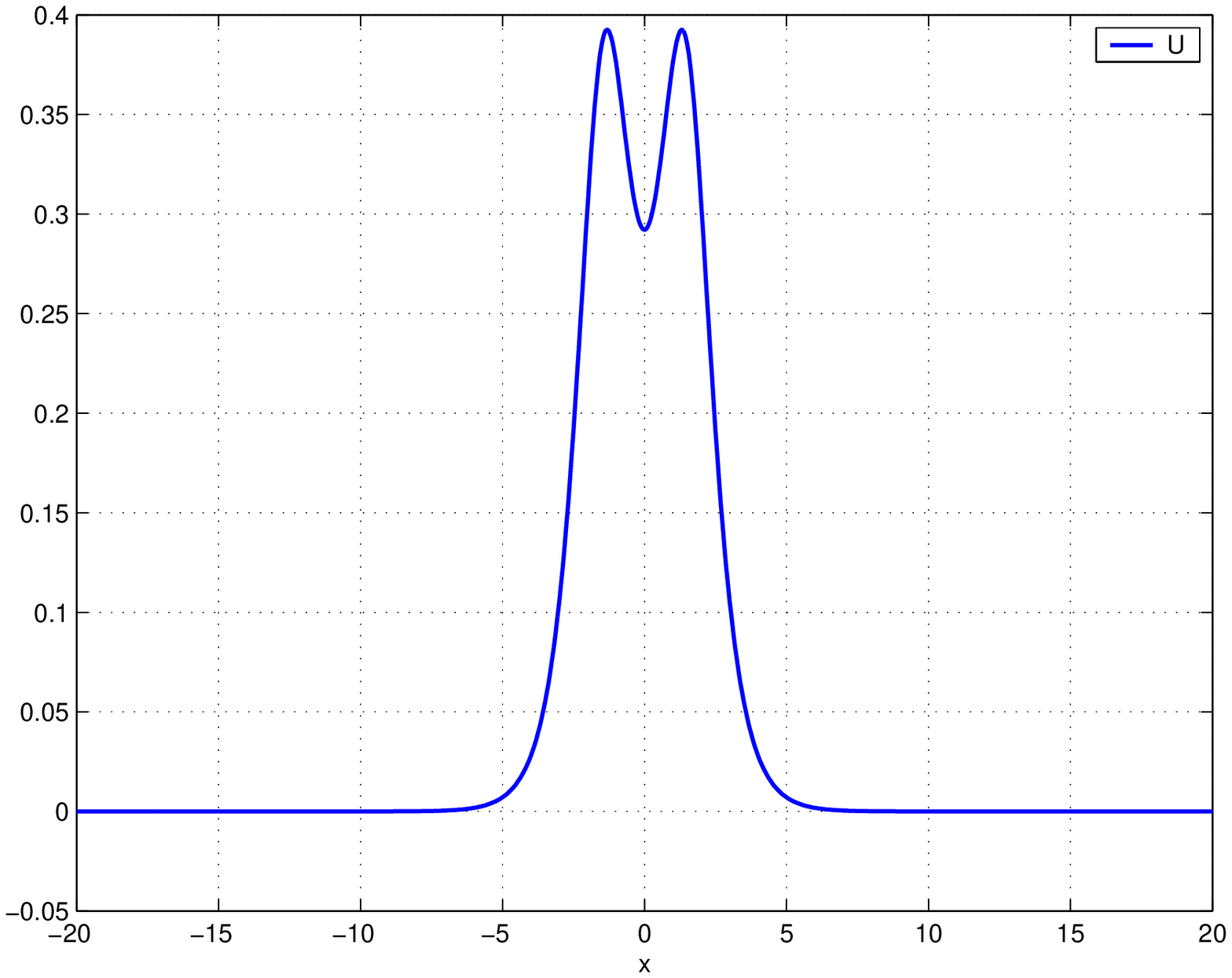}}
%\subfigure[]{
%\includegraphics[width=6.5cm]{nlsatp1d.eps}}
\subfigure[]{
\includegraphics[width=6.5cm]{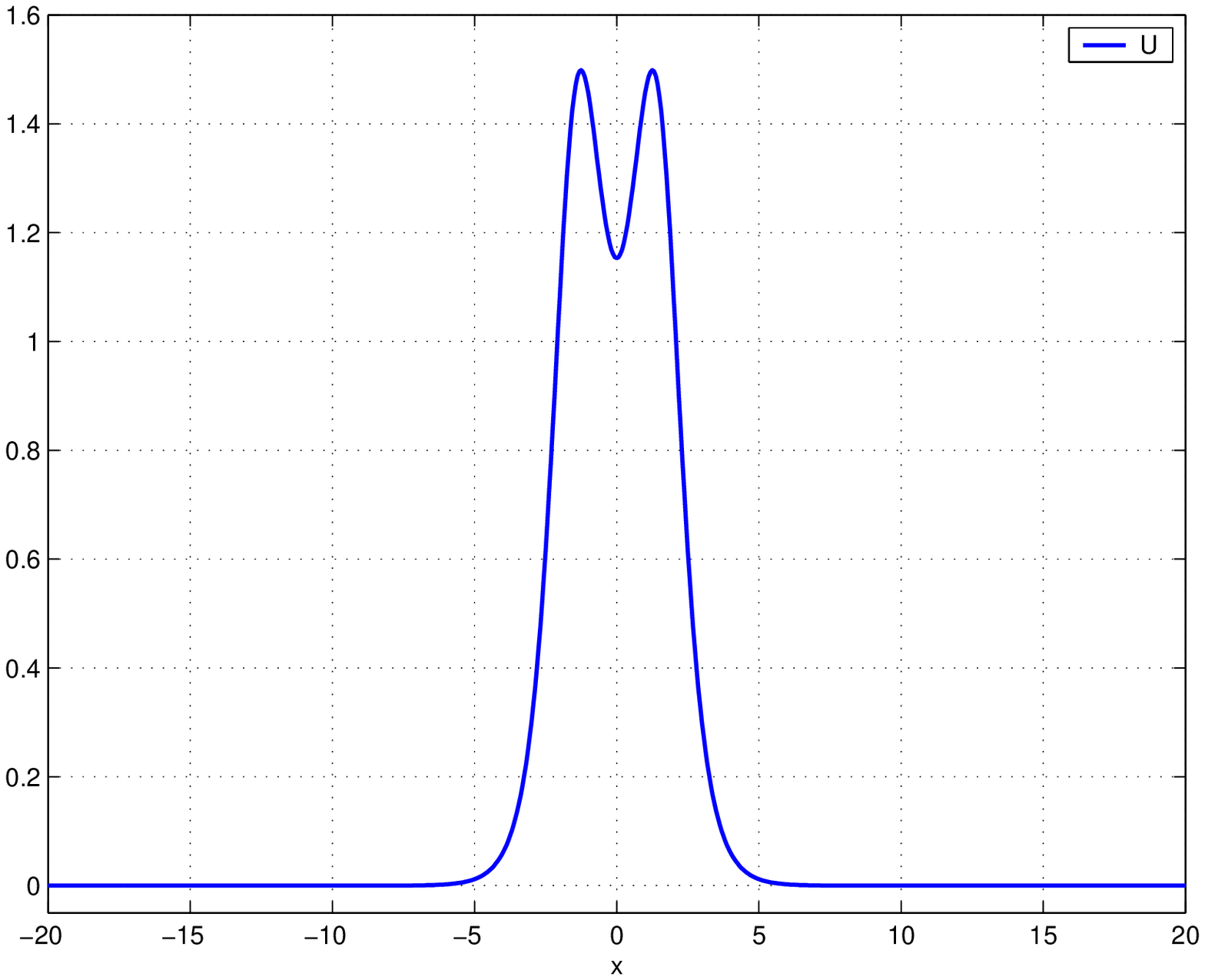}}
%\subfigure[]{
%\includegraphics[width=6.5cm]{nlsatp1d3.eps}}
\caption{Approximate profiles for (\ref{lab31b}) with $V_{0}=2.8, x_{0}=1.5, \gamma=0.25$. (a) $\mu=1.9$; (b) $\mu=2.69$.} \label{fexample_6}
\end{figure}
The convergence of the method is illustrated by the following results. Table \ref{tav_5} shows the six largest magnitude eigenvalues of the iteration matrix (\ref{lab14b}) of the classical fixed point algorithm (\ref{lab14}) and of the iteration matrix (\ref{lab28}), (\ref{lab25}) with $\gamma_{1}=3/2, \gamma_{2}=5/4$, both at the last computed iterate $U_{f}$. The results correspond to the two values of $\mu$ considered, near to two types of bifurcation: a symmetry breaking pitchfork bifurcation ($\mu=1.9$) and a saddle-node bifurcation ($\mu=2.69$), \cite{yang3}. In both cases, the presence of a {\em unique} eigenvalue of magnitude above one in the spectrum of $S$ (first and third columns) explains the nonconvergence of the classical fixed-point algorithm (\ref{lab14}). The extended method (\ref{lab28}), (\ref{lab25}) modifies the spectrum, in such a way that the harmful eigenvalue is ruled out and the rest of the spectrum is retained to be below one in magnitude. The resulting profiles where the iteration matrices are evaluated at are displayed in Figures \ref{fexample_6}(a) (for $\mu=1.9$) and (c) (for $\mu=2.69$). They are positive and symmetric, \cite{yang3}. 

\begin{figure}[htbp]
\centering
\subfigure[]{
\includegraphics[width=6.5cm]{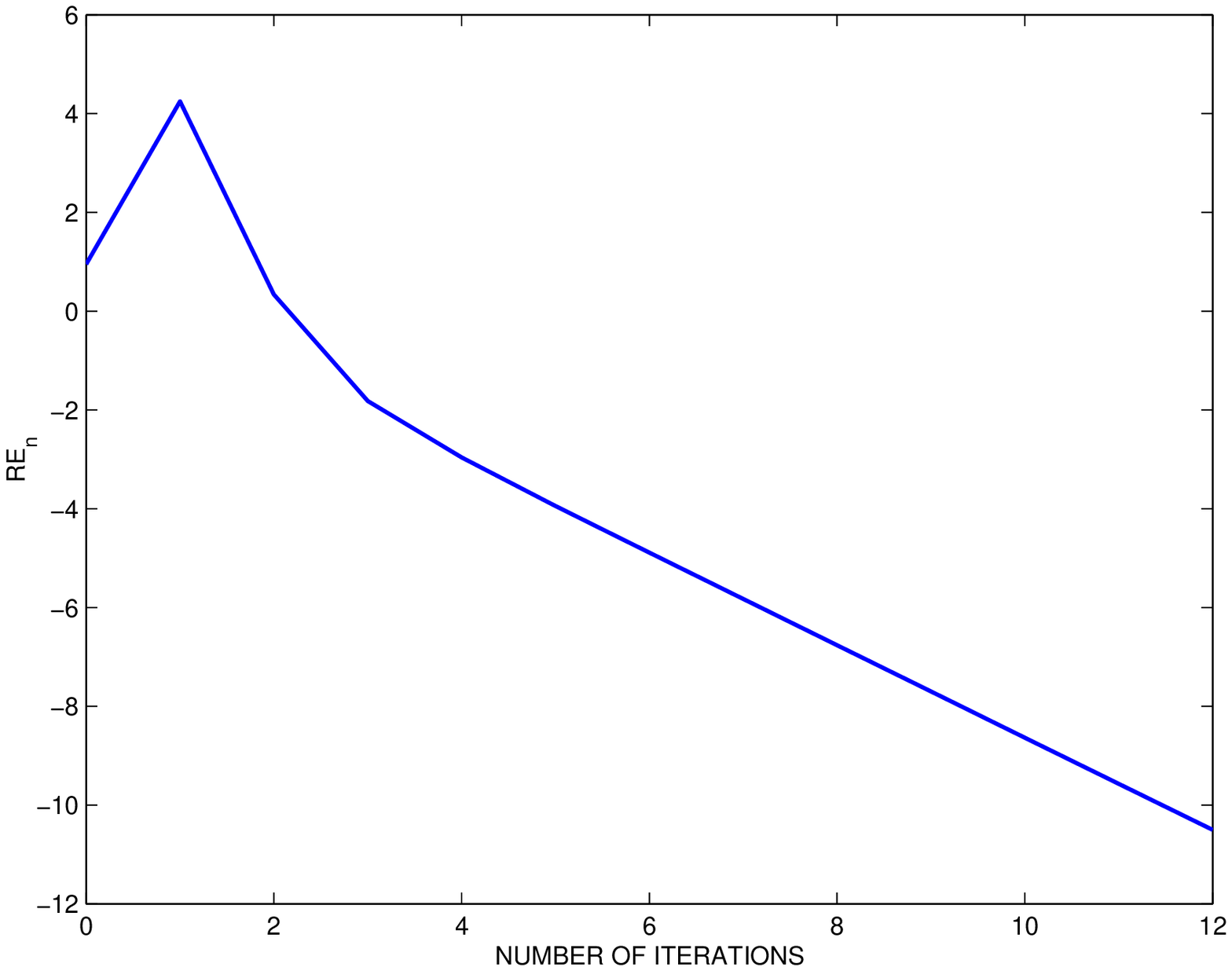}}
\subfigure[]{
\includegraphics[width=6.5cm]{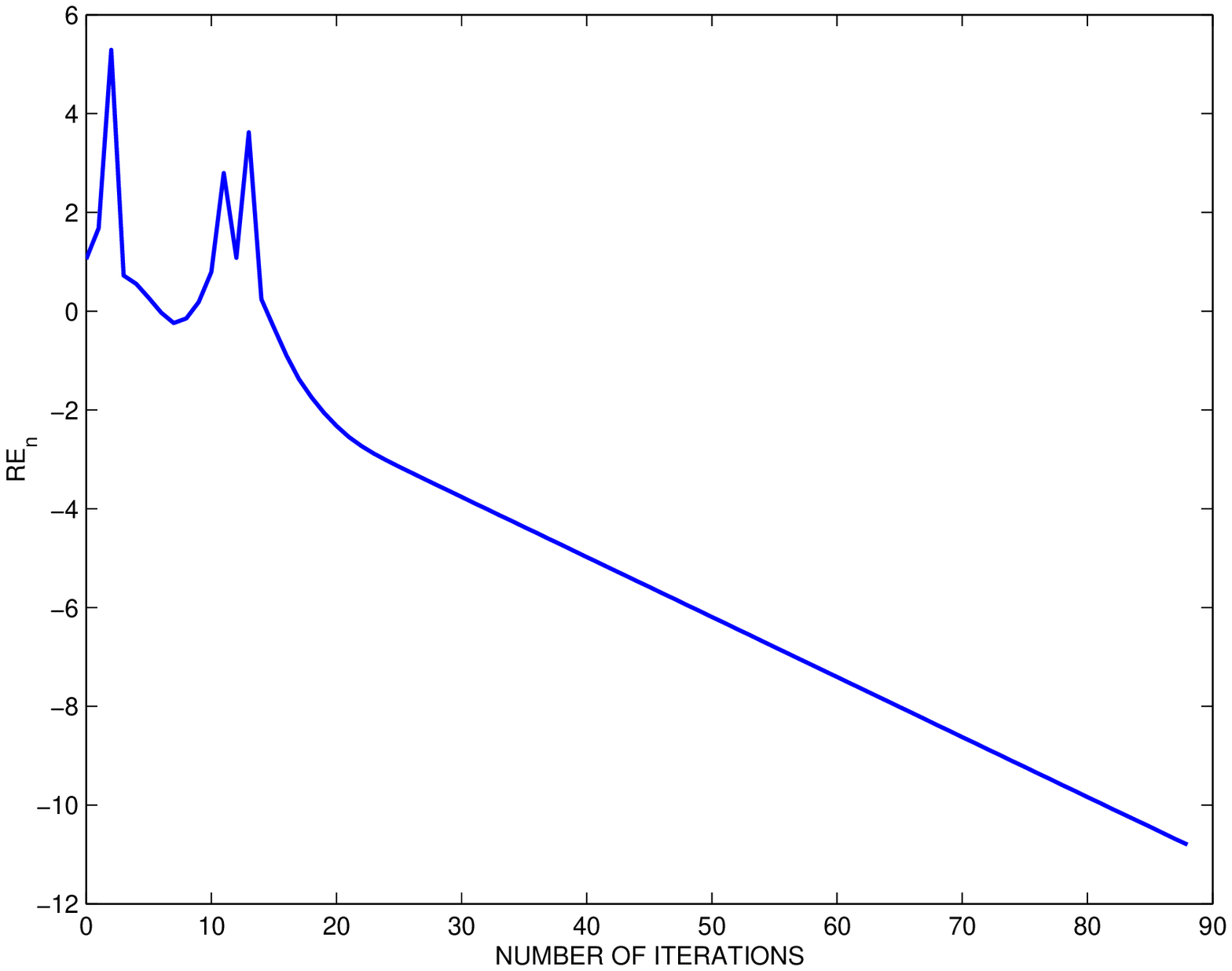}}
\caption{Residual error (\ref{res}) vs number of iterations for the generation of the profiles of Figure \ref{fexample_6}. (a) $\mu=1.9$; (b) $\mu=2.69$.} \label{fexample_6a}
\end{figure}

\begin{figure}[htbp]
\centering \subfigure[]{
\includegraphics[width=6.5cm]{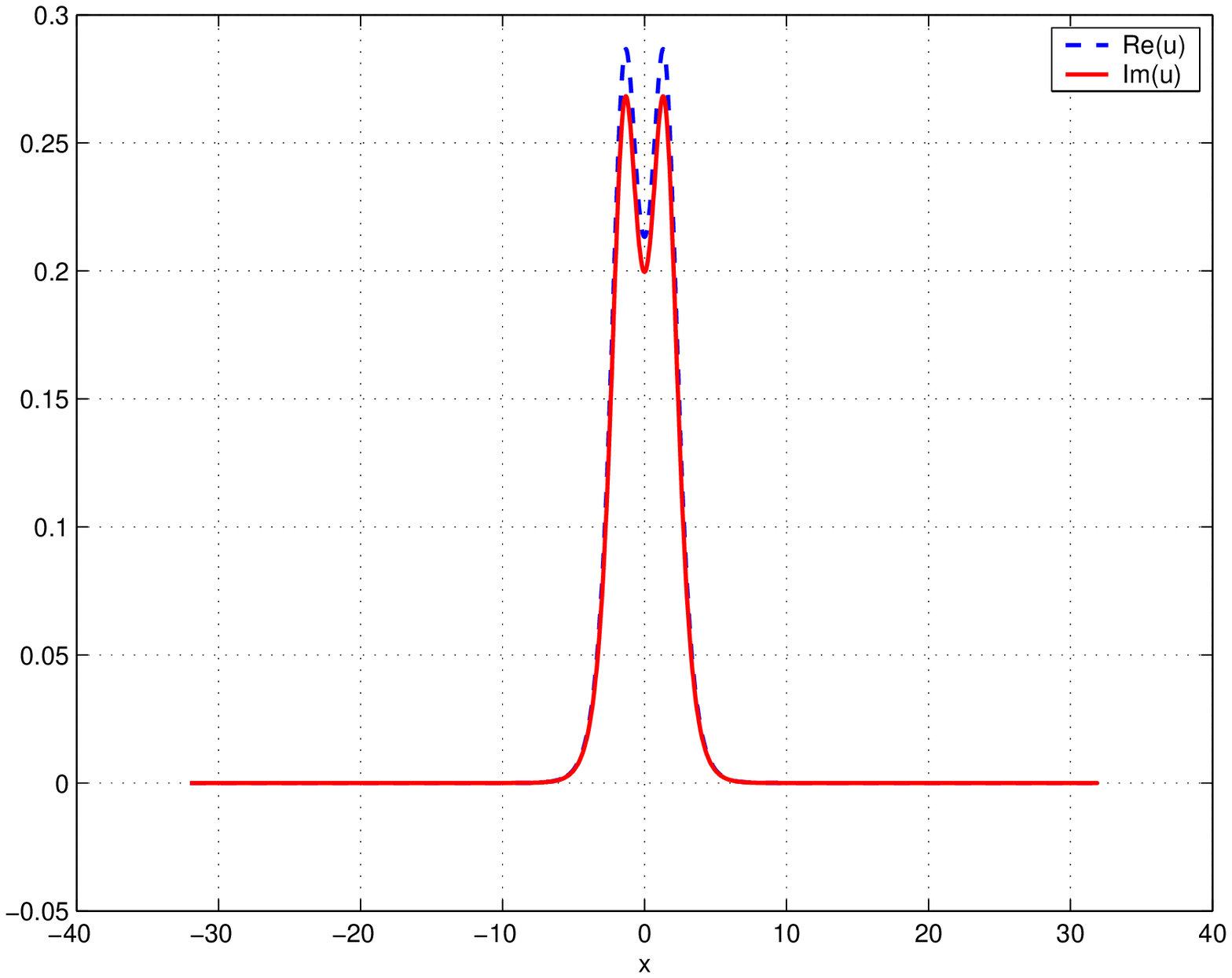}}
\subfigure[]{
\includegraphics[width=6.5cm]{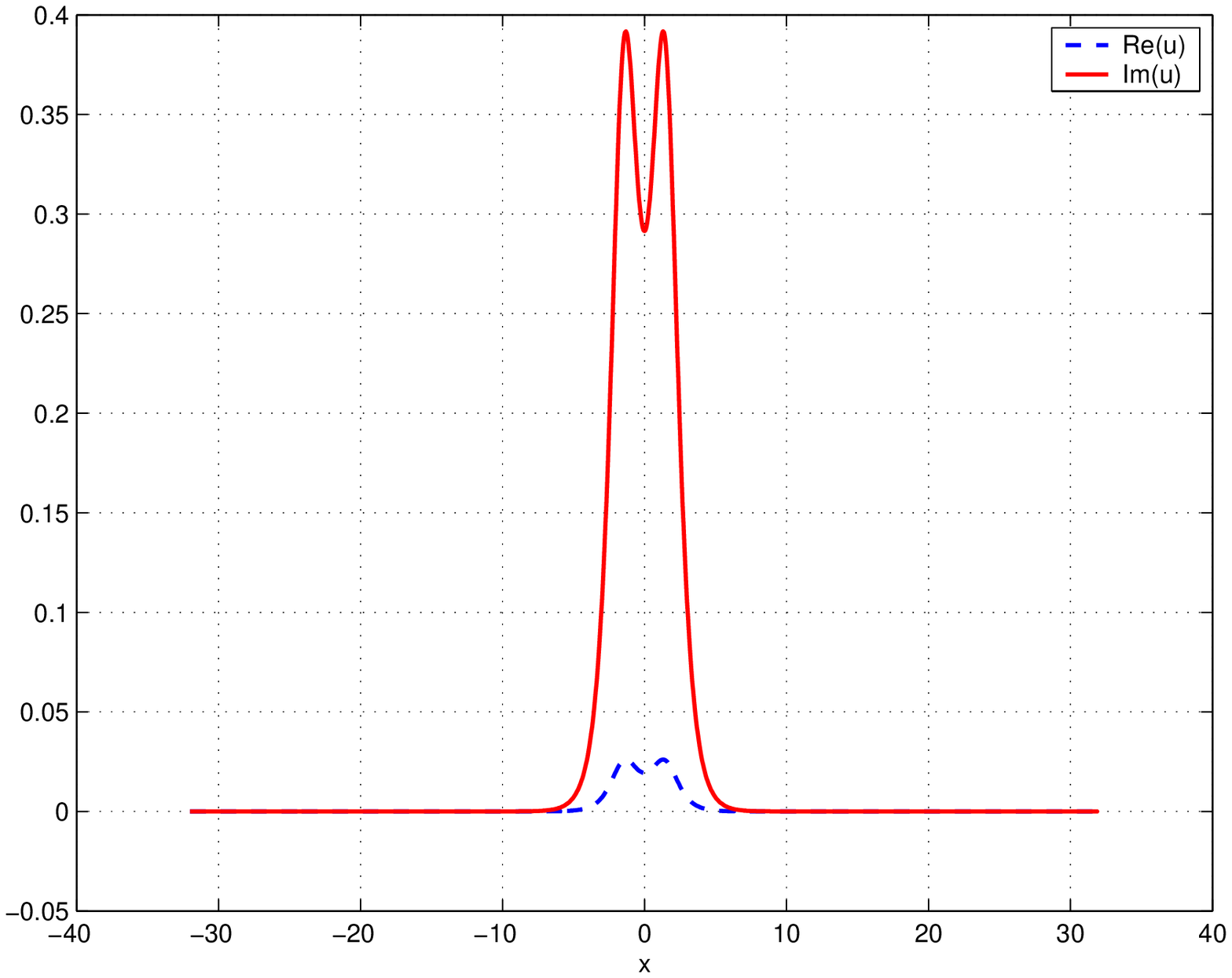}}
\subfigure[]{
\includegraphics[width=6.5cm]{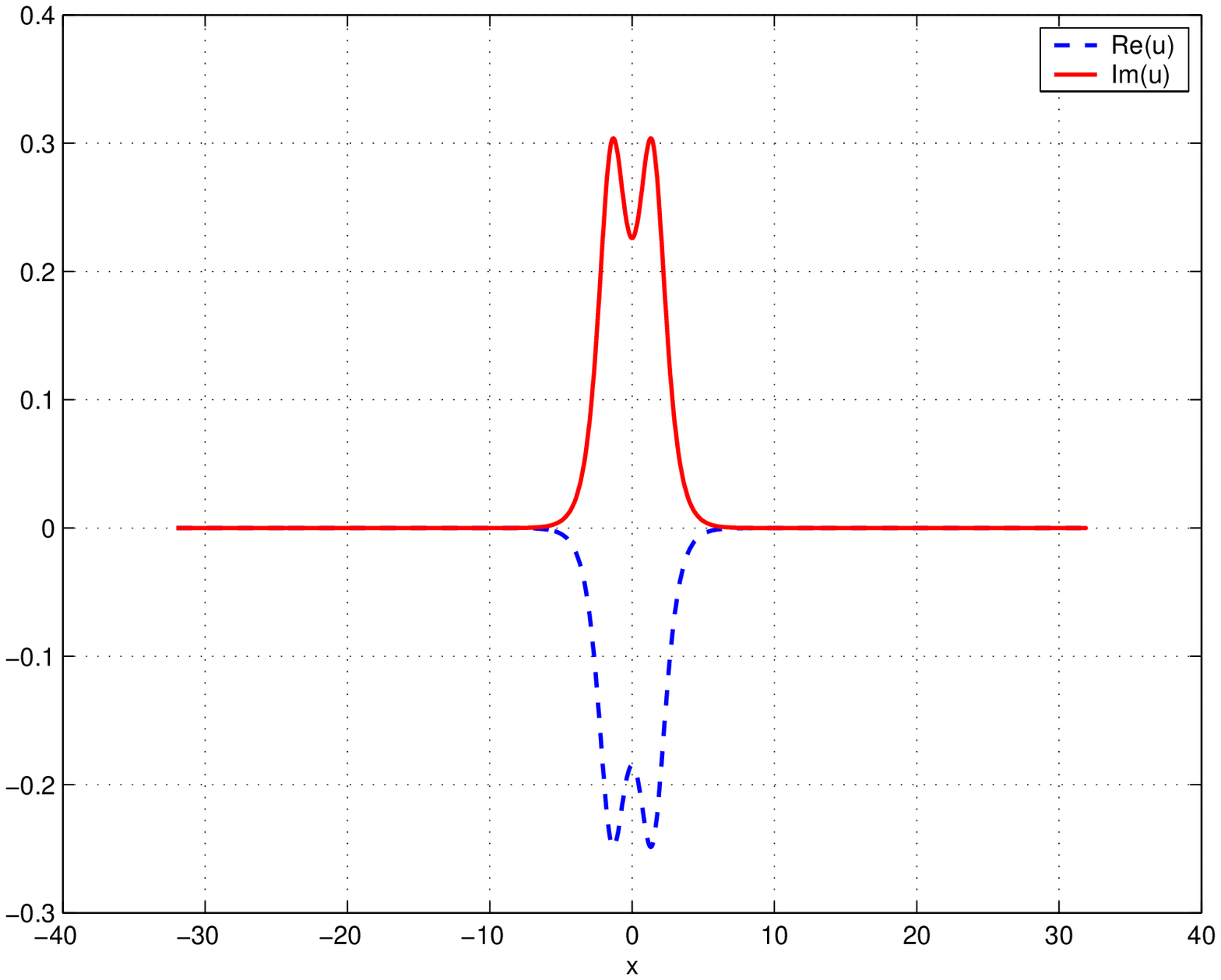}}
\subfigure[]{
\includegraphics[width=6.5cm]{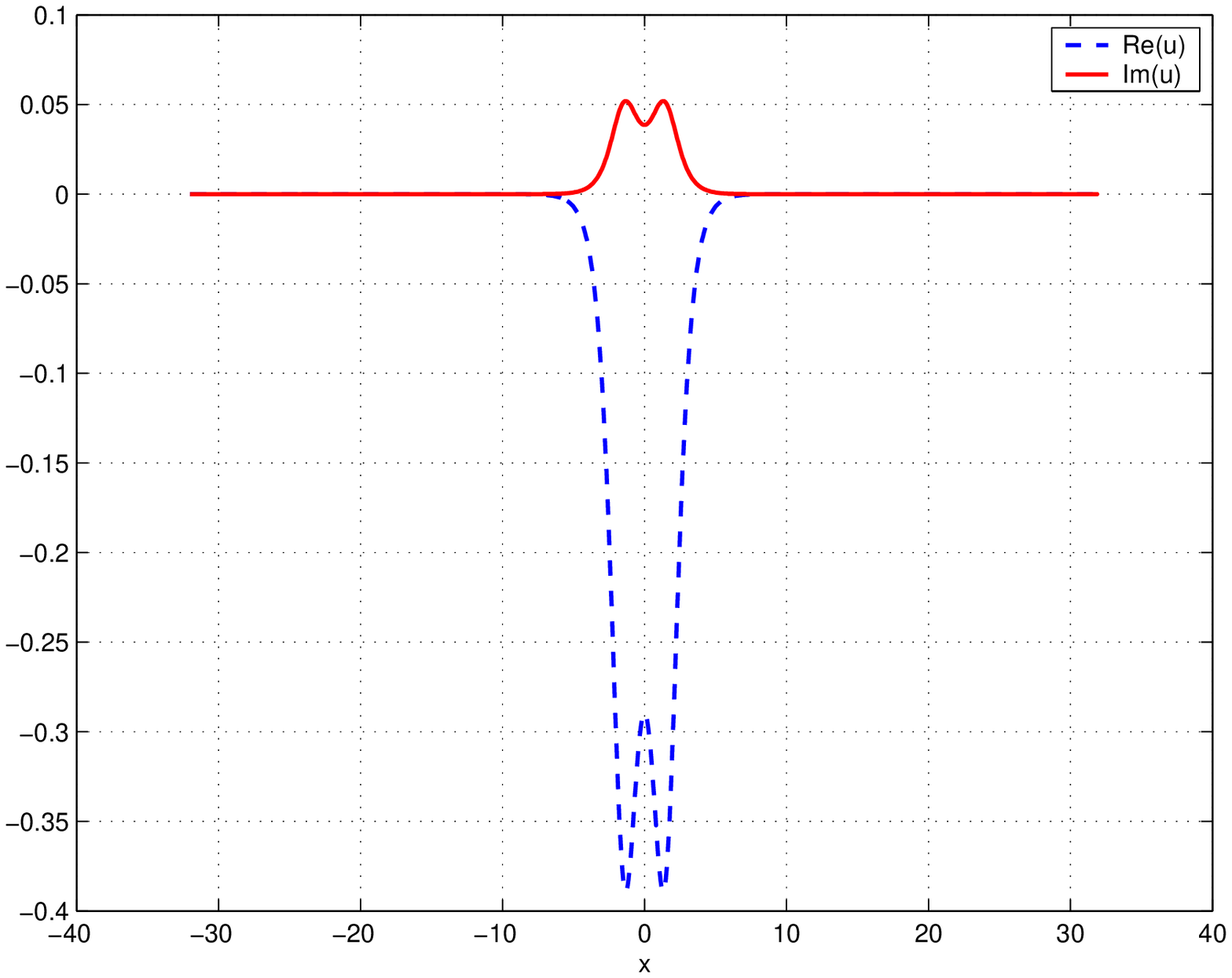}}
\caption{Evolution of the profile of Figure \ref{fexample_6}(a). Numerical solution (real part with solid line and imaginary part with dashed line) at times $t=50, 100, 150, 200$.} \label{fexample_7}
\end{figure}

\begin{figure}[htbp]
\centering \subfigure[]{
\includegraphics[width=6.5cm]{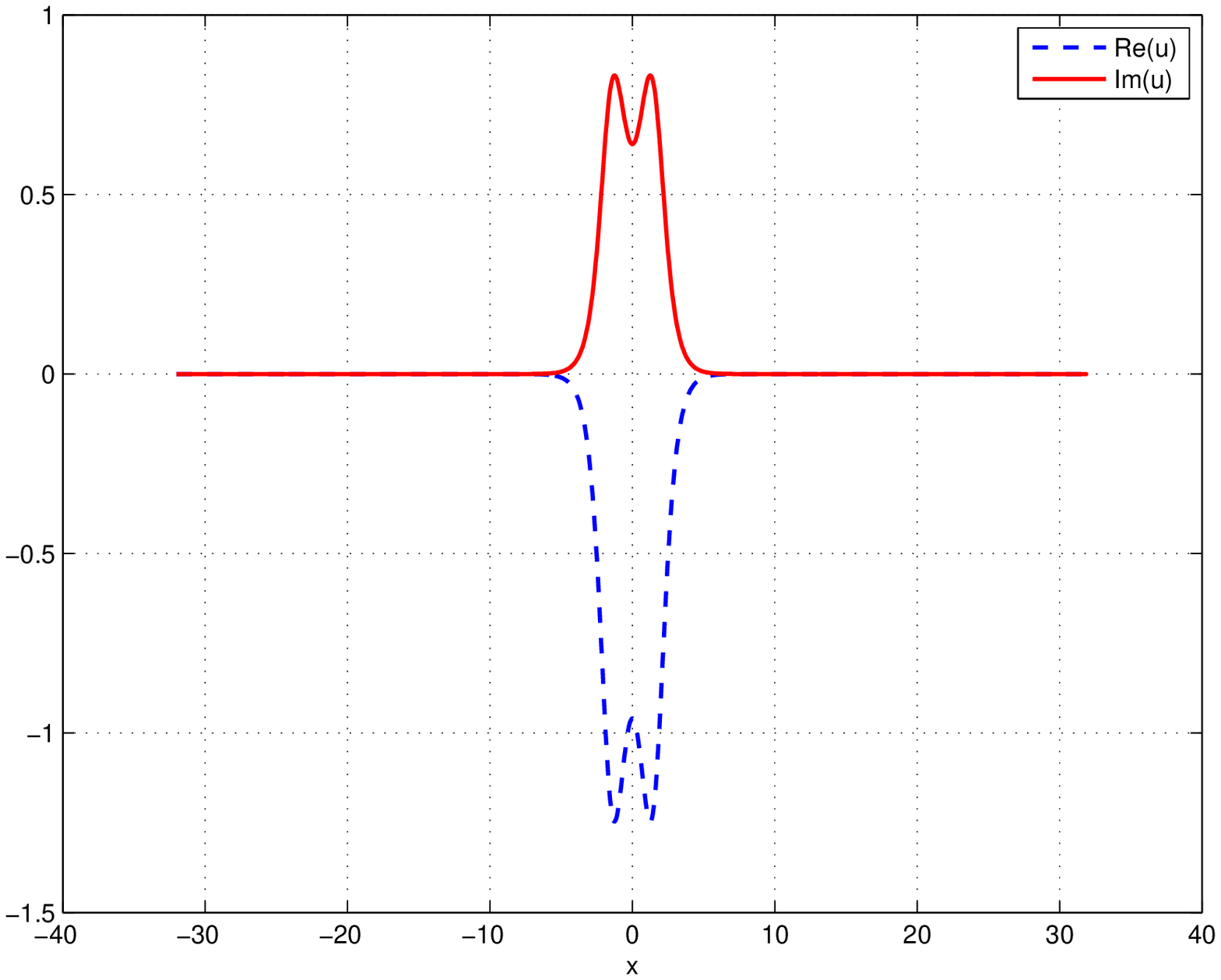}}
\subfigure[]{
\includegraphics[width=6.5cm]{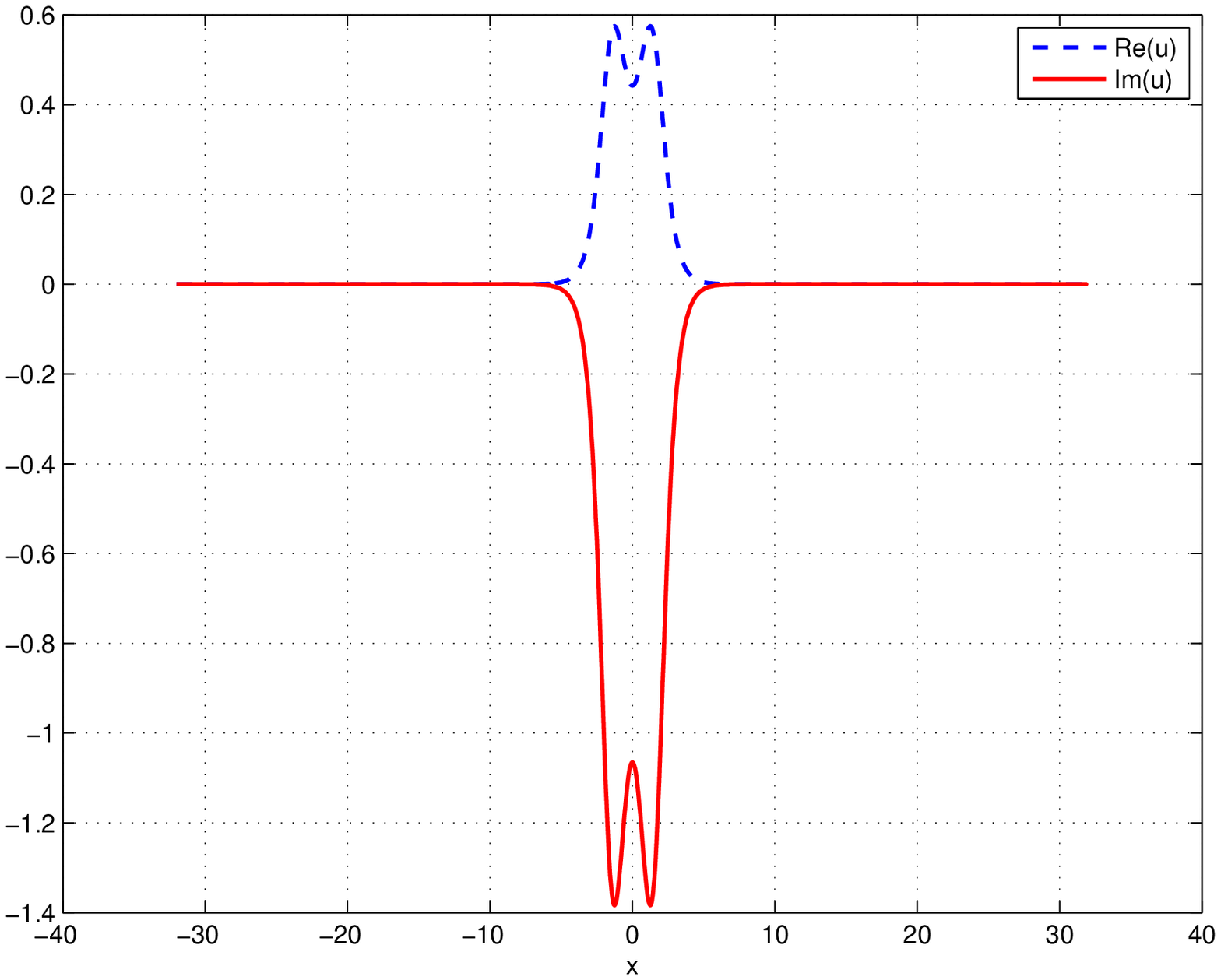}}
\subfigure[]{
\includegraphics[width=6.5cm]{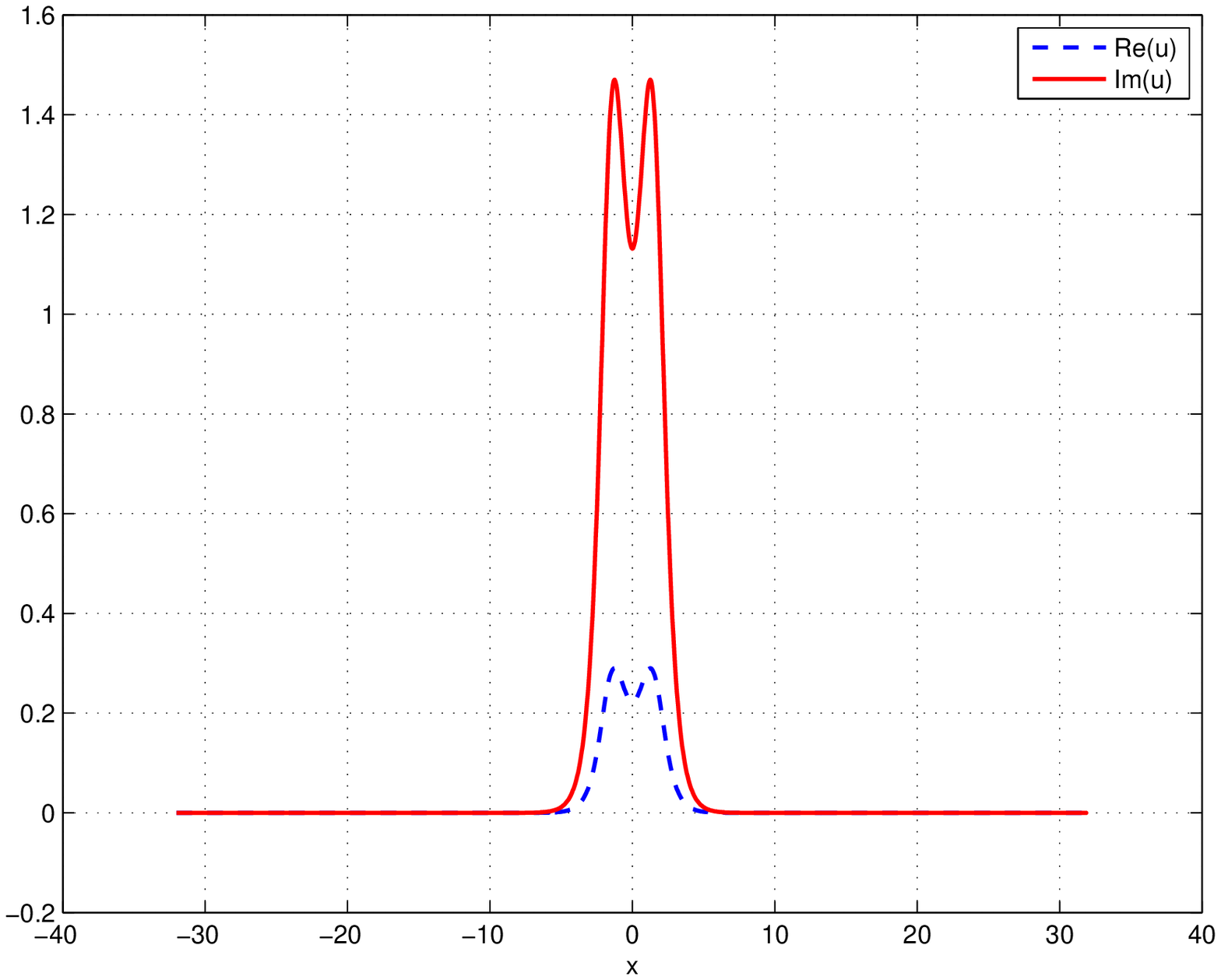}}
\subfigure[]{
\includegraphics[width=6.5cm]{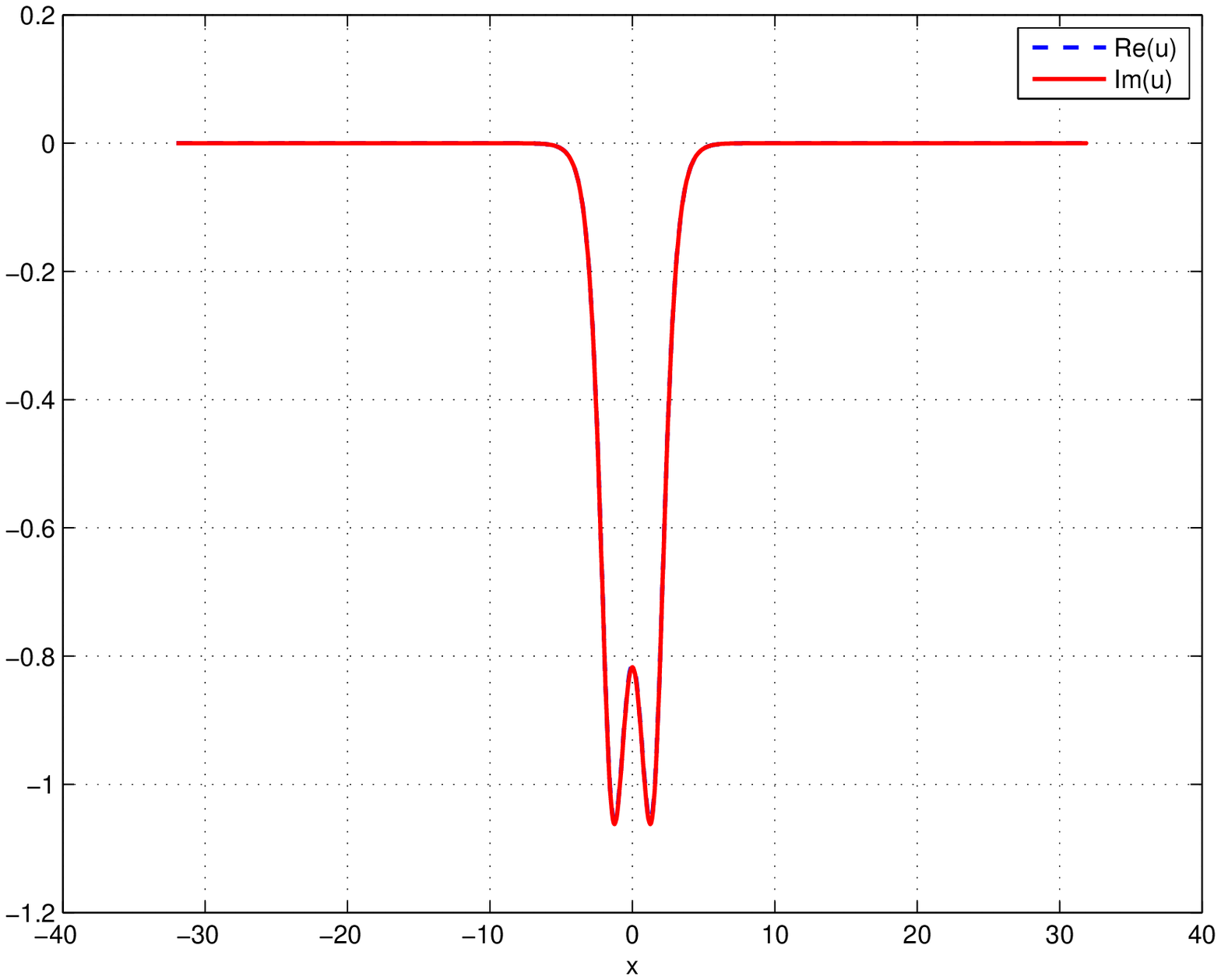}}
\caption{Evolution of the profile of Figure \ref{fexample_6}(b). Numerical solution (real part with solid line and imaginary part with dashed line) at times $t=50, 100, 150, 200$.} \label{fexample_10}
\end{figure}

The convergence is also confirmed by the next two experiments. Figure \ref{fexample_6a} shows the behaviour of the residual error (\ref{res}), where $L$ and $N$ are now given by (\ref{lab32}).
In both cases ($\mu=1.9$ for Figure \ref{fexample_6a}(a) and $\mu=2.69$ for Figure \ref{fexample_6a}(b)) the decrease of the residual is observed, with a higher computational cost in the second case.
A final test for the accuracy of the computed waves is shown in Figures \ref{fexample_7} and \ref{fexample_10}. The final iterates have been taken as initial conditions of a time-stepping code to integrate (\ref{lab31}) numerically. Figure \ref{fexample_7} illustrates the evolution of the resulting numerical solution, which is displayed (with the real and imaginary parts in a separate way) at different times for the case $\mu=1.9$. Figure \ref{fexample_10} corresponds to $\mu=2.69$. In both cases, the profile evolves as a localized ground state with a high accuracy.

\begin{figure}[htbp]
\centering
\subfigure[]{
\includegraphics[width=6.5cm]{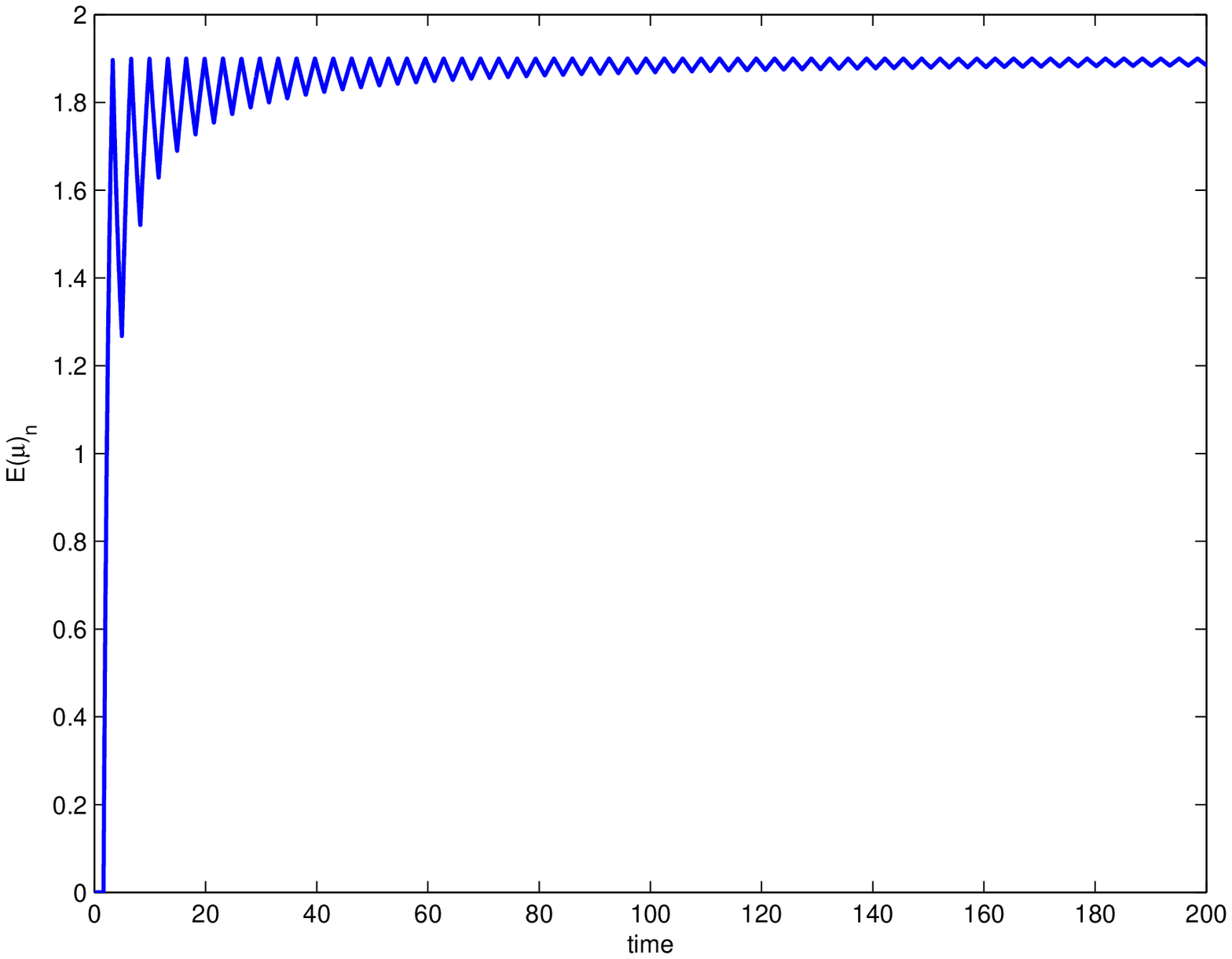}}
\subfigure[]{
\includegraphics[width=6.5cm]{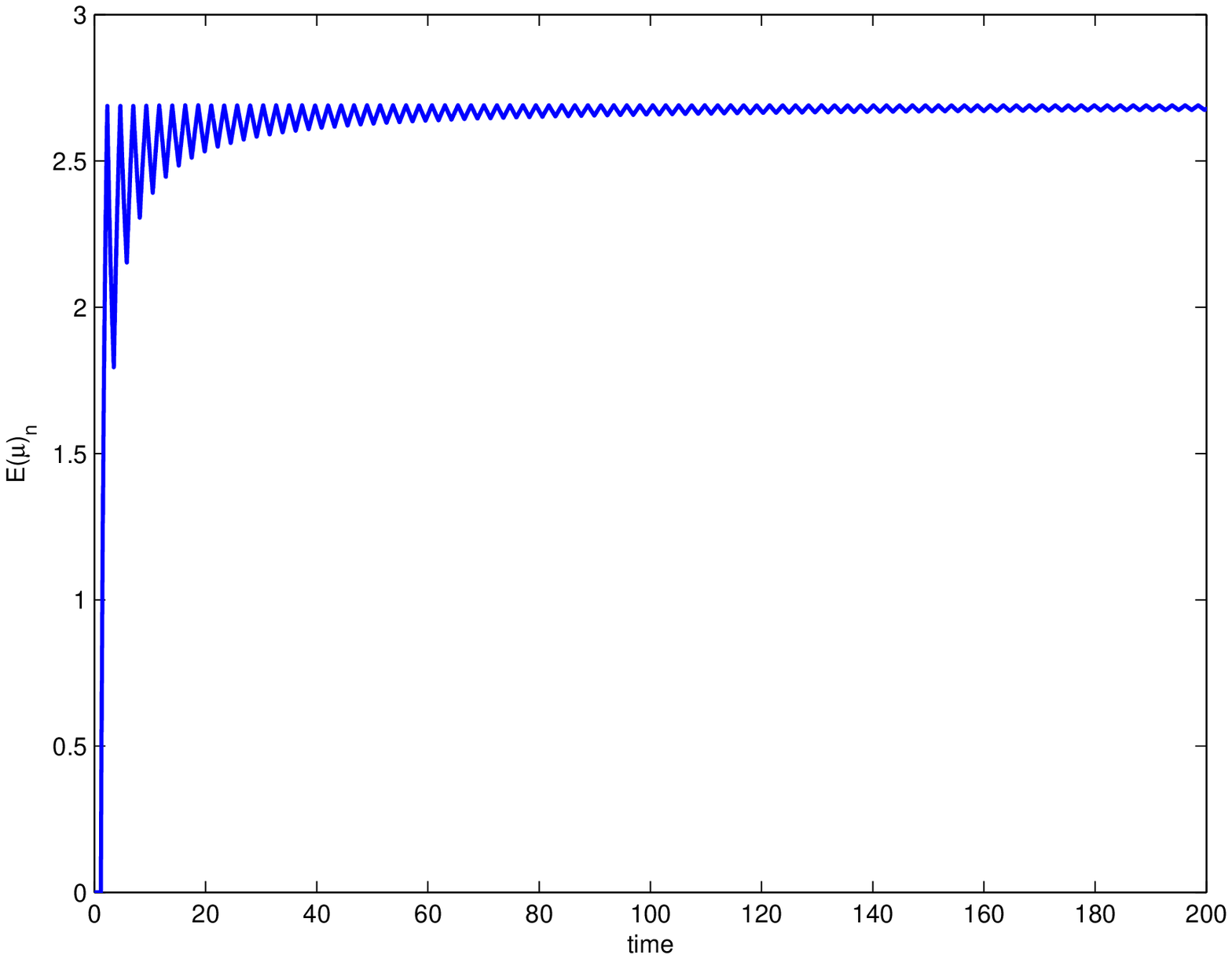}}
\caption{Evolution of the computed phase speed of the numerical solution of Figure \ref{fexample_6}. (a) $\mu=1.9$; (b) $\mu=2.69$.\label{fexample_10b}}
\end{figure}
This is also confirmed by Figure \ref{fexample_10b}, which shows the evolution of the \lq angular velocity\rq\ of the numerical solution. It is clearly observed that the computed quantity evolves to the corresponding theoretical phase speed $\mu$ of the ground state ($\mu=1.9$ in Figure \ref{fexample_10b}(a) and $\mu=1.9$ in Figure \ref{fexample_10b}(b)).

\subsection{Equations without symmetries. Example 2}
\label{sec:34}
In this last example, a generalized NLS equation (GNLS) with cubic,
quintic and seventh power nonlinearities
\begin{eqnarray}
iu_{t}+u_{xx}-V(x)u+|u|^{2}-0.2|u|^{4}u+\kappa
|u|^{6}u=0,\label{gnls2}
\end{eqnarray}
is considered. It also contains a potential $V(x)$ and a real
constant $\kappa$. The equation (\ref{gnls2}), with its physical motivation, is also
studied in \cite{yang3} (see also references therein) and the same elements have been taken here. In particular, $V$ is the 
asymmetric double-well potential
\begin{eqnarray*}
V(x)=-3.5{\rm sech}^{2}(x+1.5)-3{\rm sech}^{2}(x-1.5).\label{dwp2}
\end{eqnarray*}
and $\kappa=\kappa_{c}\approx0.01247946$.
%\item $\kappa=\kappa_{c}\approx0.01247946$.
%\end{enumerate}
Localized ground state solutions $u(x,t)=U(x)e^{i\mu t}$ now satisfy
\begin{eqnarray}
-\mu U+u^{\prime\prime}-V(x)U+|U|^{2}U-0.2|U|^{4}U+\kappa
|U|^{6}U=0,\label{gnls2b}
\end{eqnarray}
and the form of the system (\ref{lab11}) for the corresponding
Fourier collocation approximation is
\begin{eqnarray}
L&=&\mu I-D_{h}^{2}+diag(V(x_{0}),\ldots,V(x_{m-1})),\nonumber\\
N(U_{h})&=&N_{1}(U_{h})+N_{2}(U_{h})+N_{3}(U_{h})\nonumber\\
&=&\left(|U_{h}|.^{2}\right).U_{h}-0.2\left(|U_{h}|.^{4}\right).U_{h}+\kappa
\left(|U_{h}|.^{6}\right).U_{h}.\label{lab41}
\end{eqnarray}
Observe that now the nonlinearity in (\ref{lab41}) contains three homogeneous terms with degrees $p_{1}=3, p_{2}=5, p_{3}=7$. Our aim here is analyzing the performance of the extended fixed-point algorithm (\ref{lab23}) with $L=3$
\begin{eqnarray}
  Lx_{n+1}&=&s_{1}(x_{n})N_{1}(x_{n})+s_{2}(x_{n})N_{2}(x_{n})+s_{3}(x_{n})N_{3}(x_{n}), \quad n=0,1,\ldots,\label{lab22e}
\end{eqnarray}
for some factors $s_{j}, j=1,2,3$. In particular, the experiments are focused on the extension of (\ref{lab25}) by taking
\begin{eqnarray}
s_{j}(x)=m(x)^{\gamma_{j}},\quad \gamma_{j}=\frac{p_{j}}{p_{j}-1},\quad j=1,2,3,\label{lab25e}
\end{eqnarray}
where $m$ is the stabilizing factor (\ref{Stb}). Other alternatives are indeed possible.

\begin{figure}[htbp]
\centering \subfigure[]{
\includegraphics[width=6.5cm]{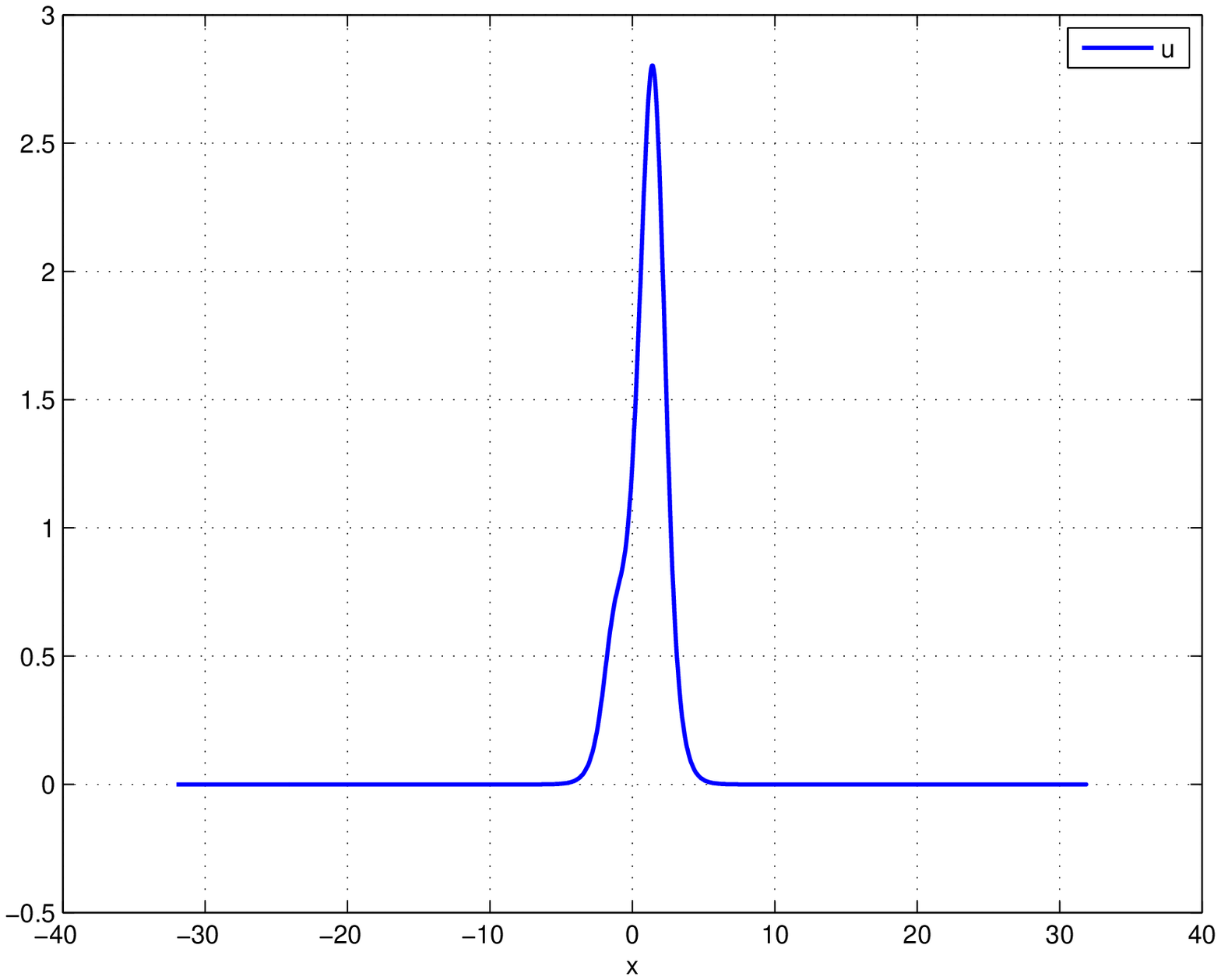}}
\subfigure[]{
\includegraphics[width=6.5cm]{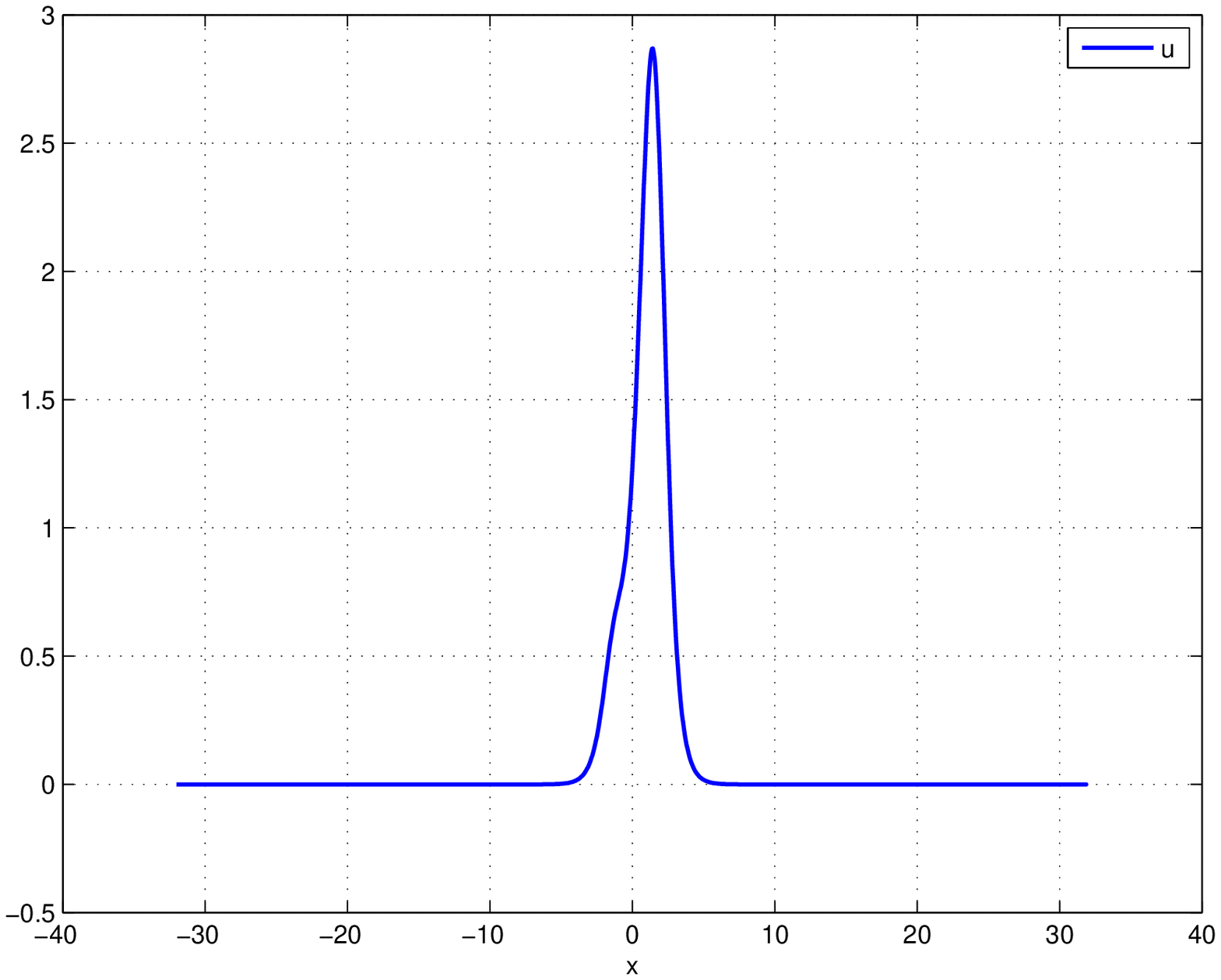}}
%\subfigure[]{
%\includegraphics[width=6.5cm]{nlsat2p1d.eps}}
%\subfigure[]{
%\includegraphics[width=6.5cm]{nlsat2p1c.eps}}
\caption{Computed wave profiles with  (\ref{lab22e}), (\ref{lab25e}) to approximate (\ref{gnls2b}). (a) $\mu=3.275$; (b) $\mu=3.289$.} \label{fexample_11}
\end{figure}

In equation (\ref{gnls2}), \cite{yang3}, transcritical bifurcations of solitary waves is found at $\kappa\approx0.01247946$ with a bifurcation point at $(\mu_{0},P_{0})\approx (3.28,14.35)$, where $P_{0}=P(\mu_{0})$ is given by (\ref{power}). The extended method (\ref{lab22e}), (\ref{lab25e}) has been checked close to this point; explicitly, solitary wave profiles have been generated for the values $\mu=3.275$ and $\mu=3.289$. They are in Figures \ref{fexample_11}(a) and (b), respectively. In both cases, the computed waves are anti-symmetric.
\begin{table}
\small
\begin{center}
\begin{tabular}{|c|c|}\hline\hline
\multicolumn{2}{|c|}
{$\mu=3.275$} \\
\hline
$S=L^{-1}N^{\prime}(U_{f})$&$F^{\prime}(U_{f})$\\\hline
1.5843078E+00&9.831590E-01\\
9.845328E-01&4.505360E-01\\
4.817932E-01&4.002412E-01\\
3.678772E-01&3.156180E-01+3.173535E-02i\\
1.979974E-01&3.156180E-01-3.173535E-02i\\
1.376471E-01&1.658000E-01\\\hline\hline
\multicolumn{2}{|c|}{$\mu=3.289$}\\
\hline
$S=L^{-1}N^{\prime}(U_{f})$&$F^{\prime}(U_{f})$\\\hline
1.857527E+00&9.429134E-01\\
9.430769E-01&4.720197E-01\\
3.759696E-01&1.766474E-01+2.501554E-01i\\
3.759696E-01&1.766474E-01-2.501554E-01i\\
1.483127E-01&1.836617E-01\\\hline
\end{tabular}
\end{center}
\caption{Six largest magnitude eigenvalues of the iteration matrices
(\ref{lab14b}) and (\ref{lab22e}), (\ref{lab25e}) for $\mu=3.275$ (first and second columns) and $\mu=3.289$ (third and fourth columns). $U_{f}$ stands for the last computed iterate. \label{tav_6}}
\end{table}

The convergent effect of the procedure (\ref{lab22e}), (\ref{lab25e}), compared to the classical fixed-point iteration, is shown in Table \ref{tav_6}. This displays the six largest magnitude eigenvalues of the iteration matrices (\ref{lab14b}) and the one of (\ref{lab22e}), (\ref{lab25e})
\begin{eqnarray}
F^{\prime}(U_{f})=\sum_{j=1}^{3} \frac{\gamma_{j}}{\langle N_{j}(U_{f}),U_{f}\rangle}N_{j}(U_{f})U_{f}^{T}(I-S),\label{itermat2}
\end{eqnarray}
at the last computed iterate $U_{f}$, for the two values of $\mu$ considered and where $N_{j}, \gamma_{j}, j=1,2,3$ are given by (\ref{lab41}) and (\ref{lab25e}), respectively. Again, (\ref{lab14b}) contains an only eigenvalue of modulus greater than one, while (\ref{itermat2}) translates the spectrum in such a way that the corresponding spectral radius is below one; this makes the iteration convergent.

\begin{figure}[htbp]
\centering
\subfigure[]{
\includegraphics[width=6.5cm]{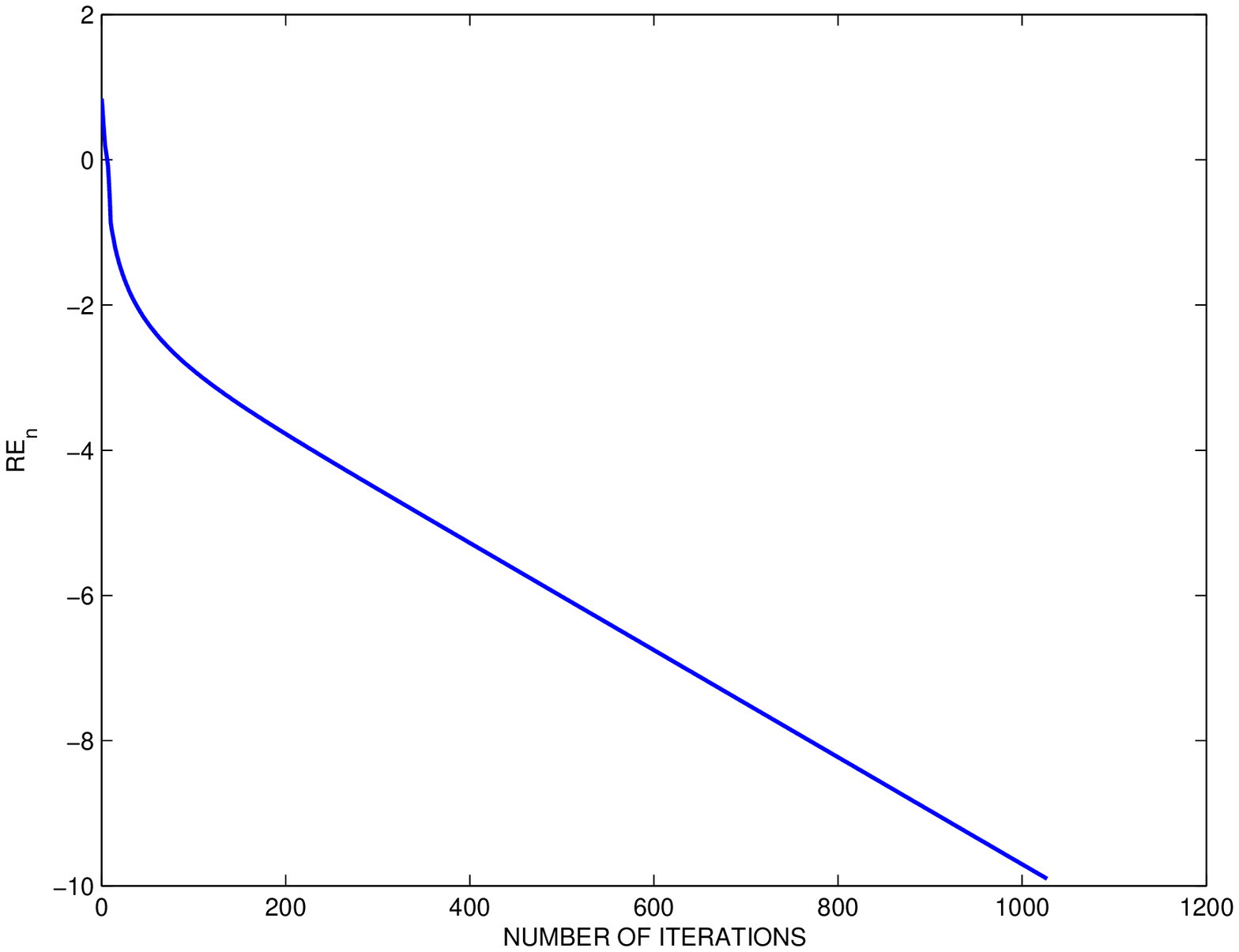}}
\subfigure[]{
\includegraphics[width=6.5cm]{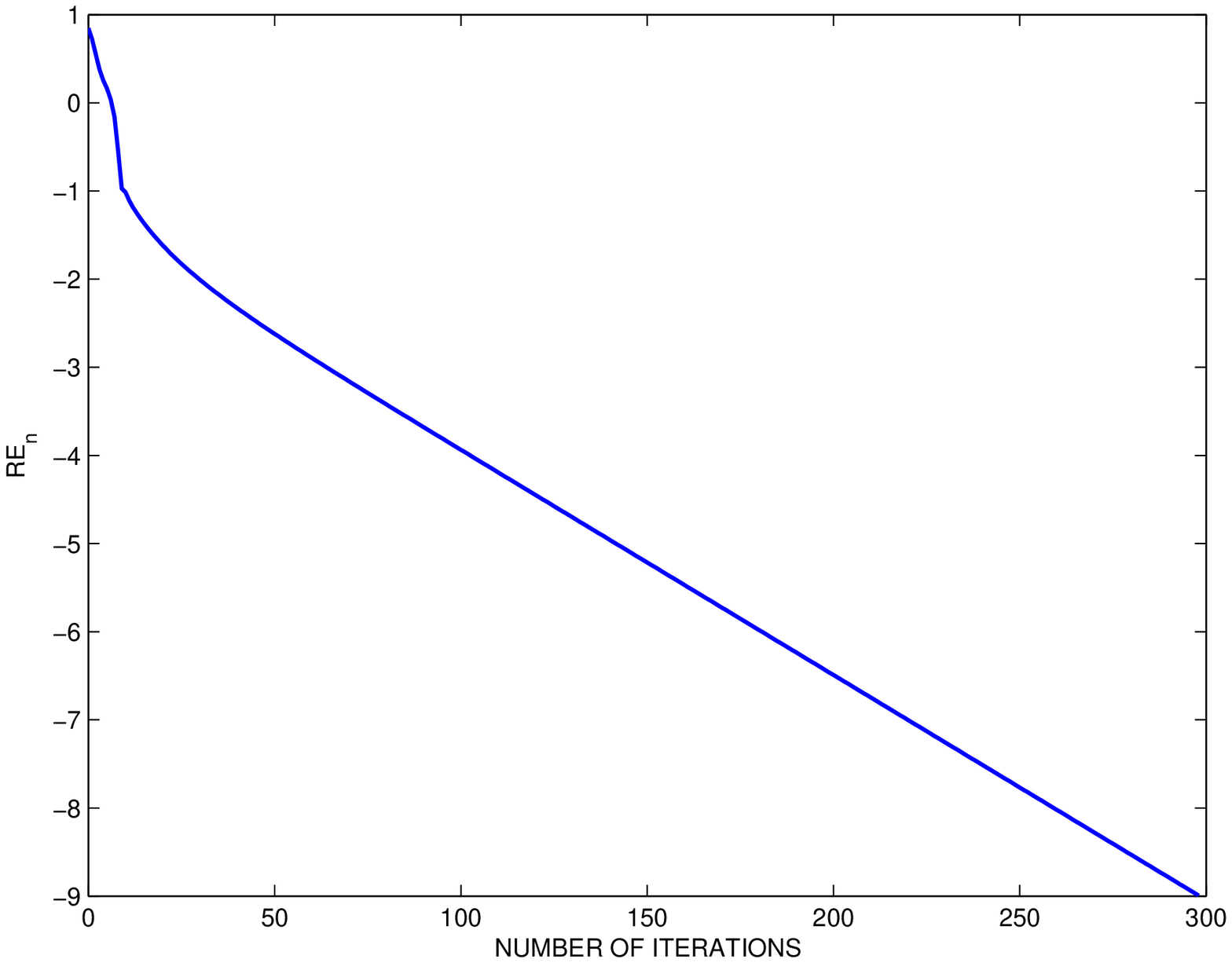}}
\caption{Residual error (\ref{res}) vs number of iterations for the generation of the profiles of Figure \ref{fexample_11}. (a) $\mu=3.275$; (b) $\mu=3.289$.} \label{fexample_11b}
\end{figure}

As far as the accuracy is concerned, the same quality controls as those of the previous examples are checked. Thus,
a second test to check the convergence  (if Table \ref{tav_6} is the first one) is shown in Figure \ref{fexample_11b}. This displays the behaviour of the residual error (\ref{res}), where $L$ and $N$ are now given by (\ref{lab41}), as function of the number of iterations, for (a) $\mu=3.275$ and (b) $\mu=3.289$. It is observed that the computation is harder than that of the previous example, cf. Figure \ref{fexample_6a}, since the number of iterations required to obtain a fixed level of error has increased, in some cases, in about two orders of magnitude.

\begin{figure}[htbp]
\centering \subfigure[]{
\includegraphics[width=6.5cm]{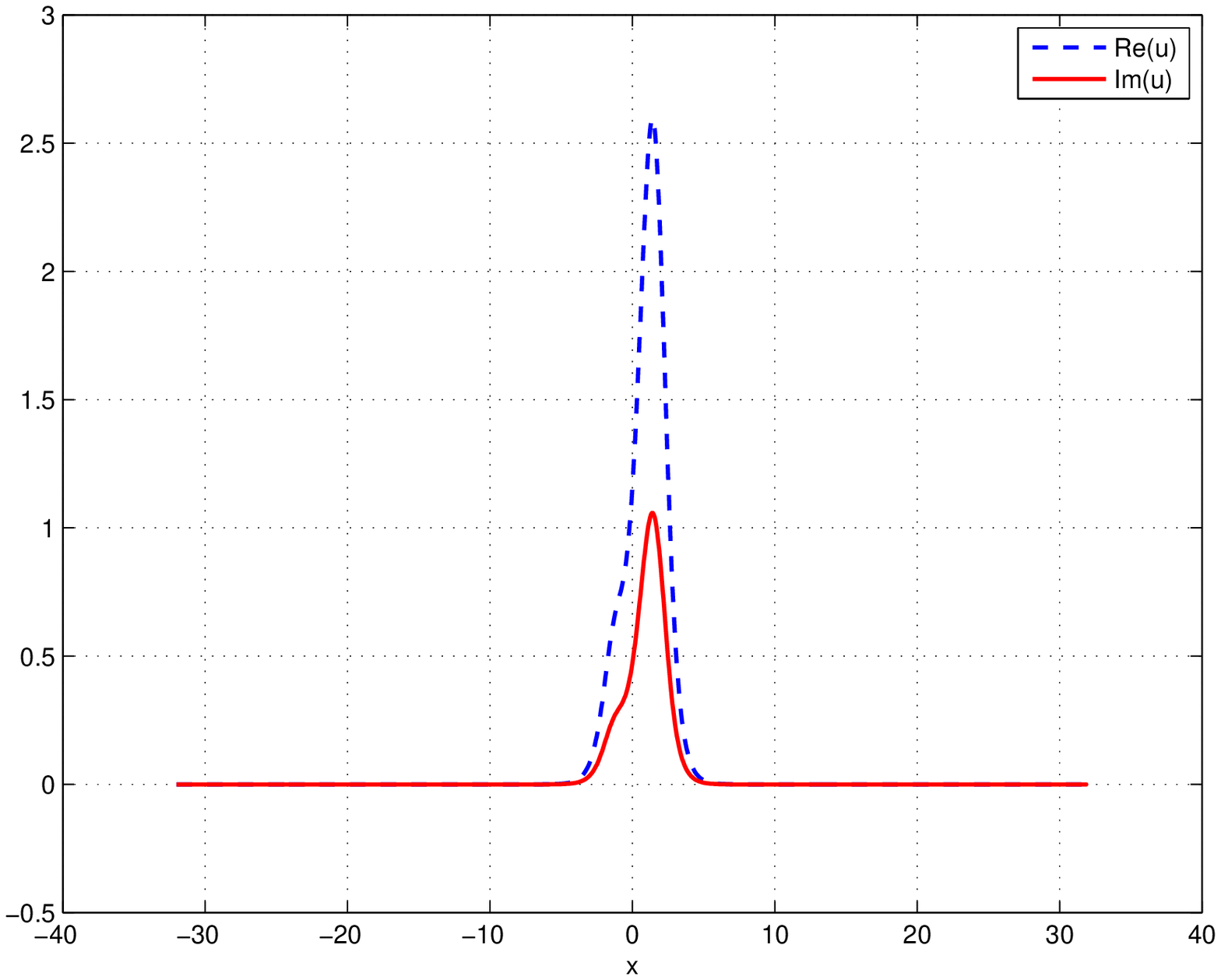}}
\subfigure[]{
\includegraphics[width=6.5cm]{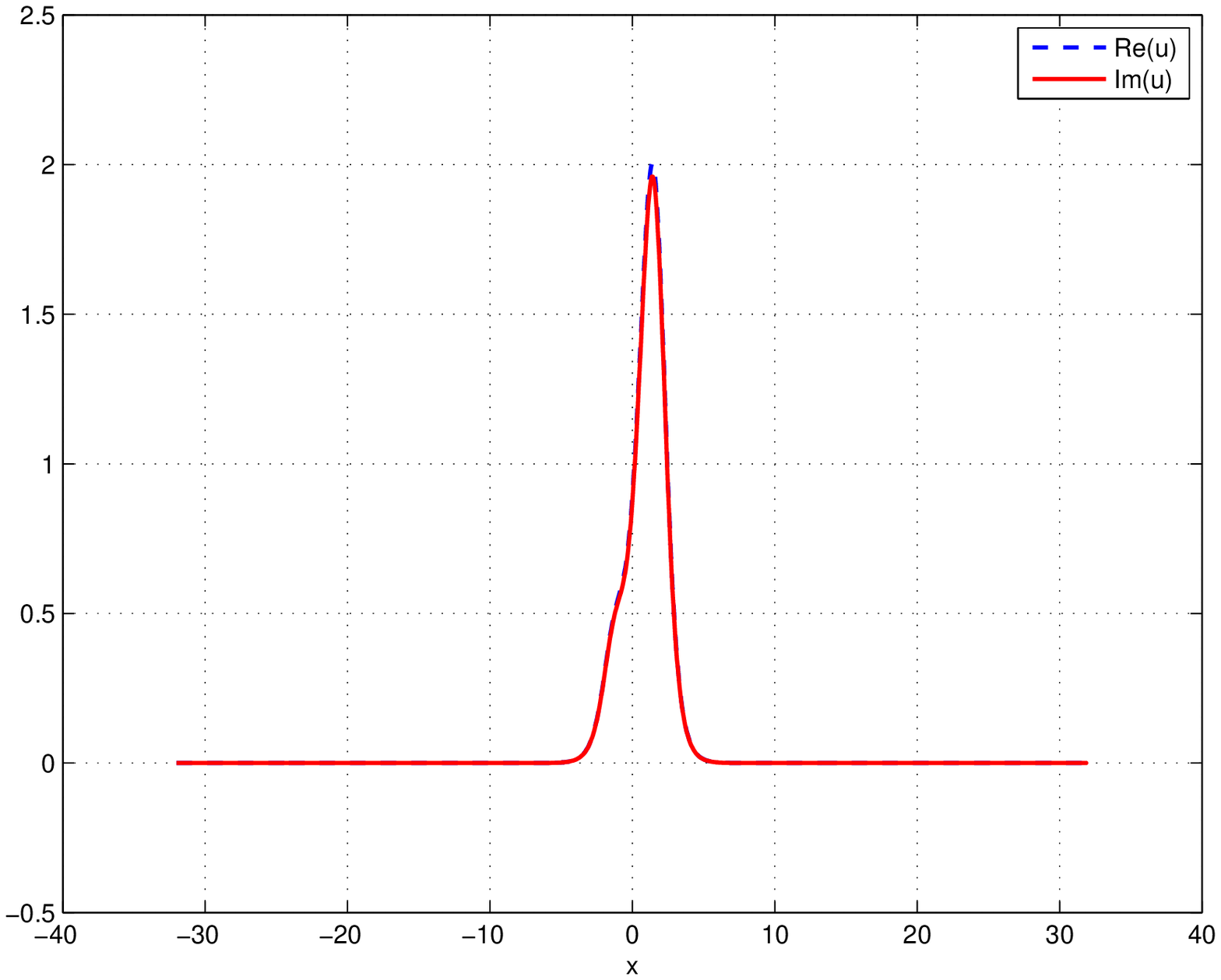}}
\subfigure[]{
\includegraphics[width=6.5cm]{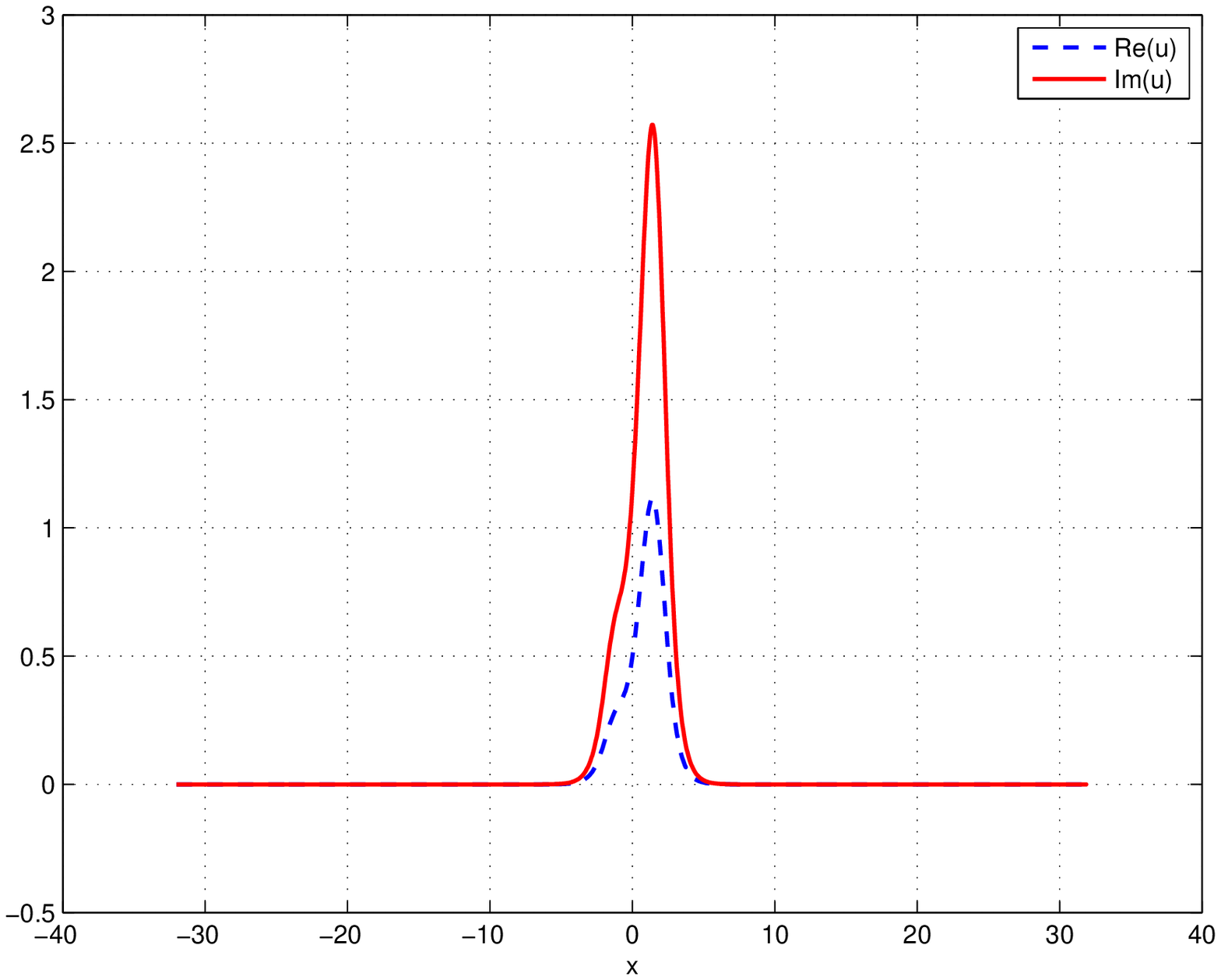}}
\subfigure[]{
\includegraphics[width=6.5cm]{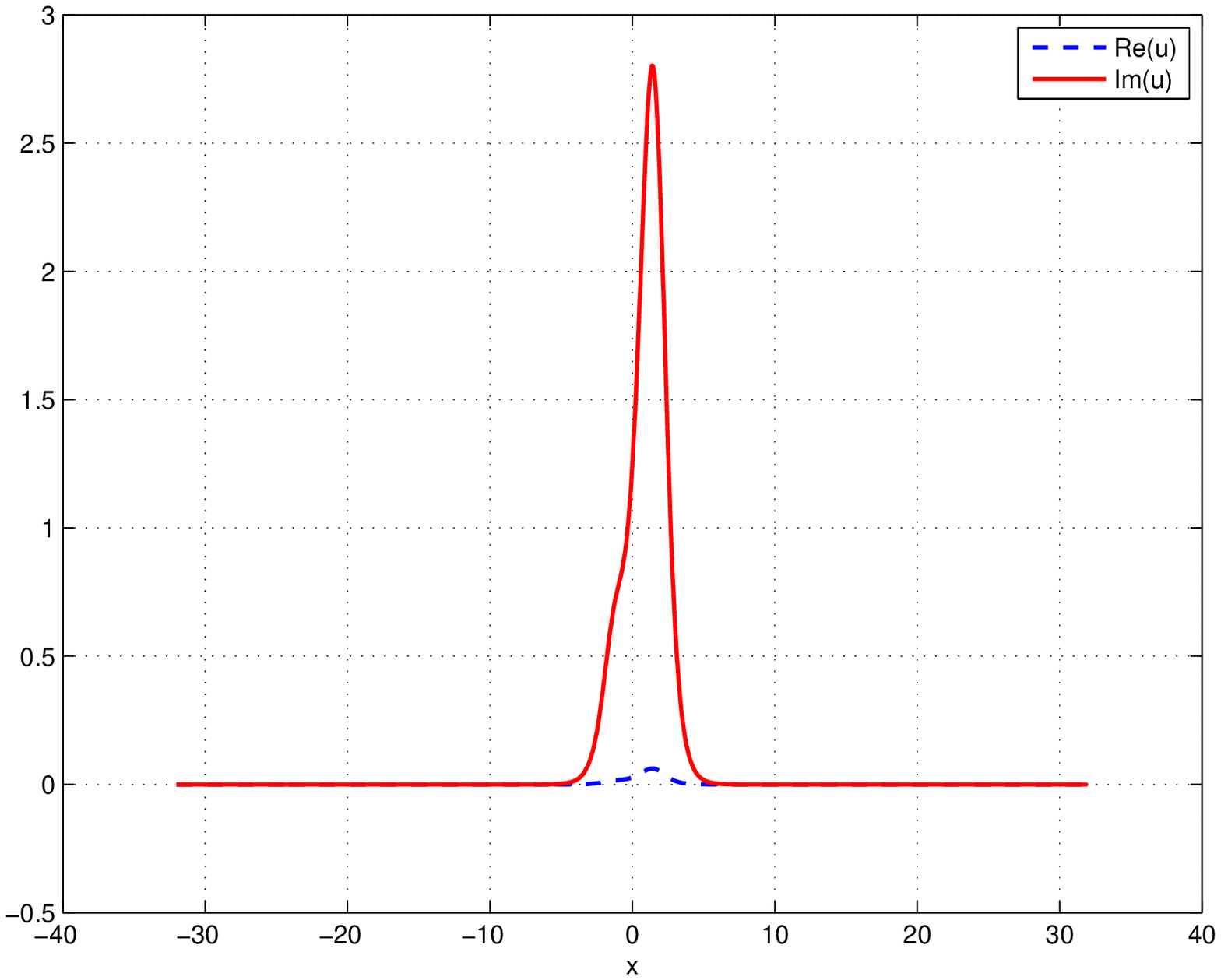}}
\caption{Evolution of the profile of Figure \ref{fexample_11}(a). Numerical solution (real part with solid line and imaginary part with dashed line) at times $t=50, 100, 150, 200$.} \label{fexample_12}
\end{figure}

\begin{figure}[htbp]
\centering \subfigure[]{
\includegraphics[width=6.5cm]{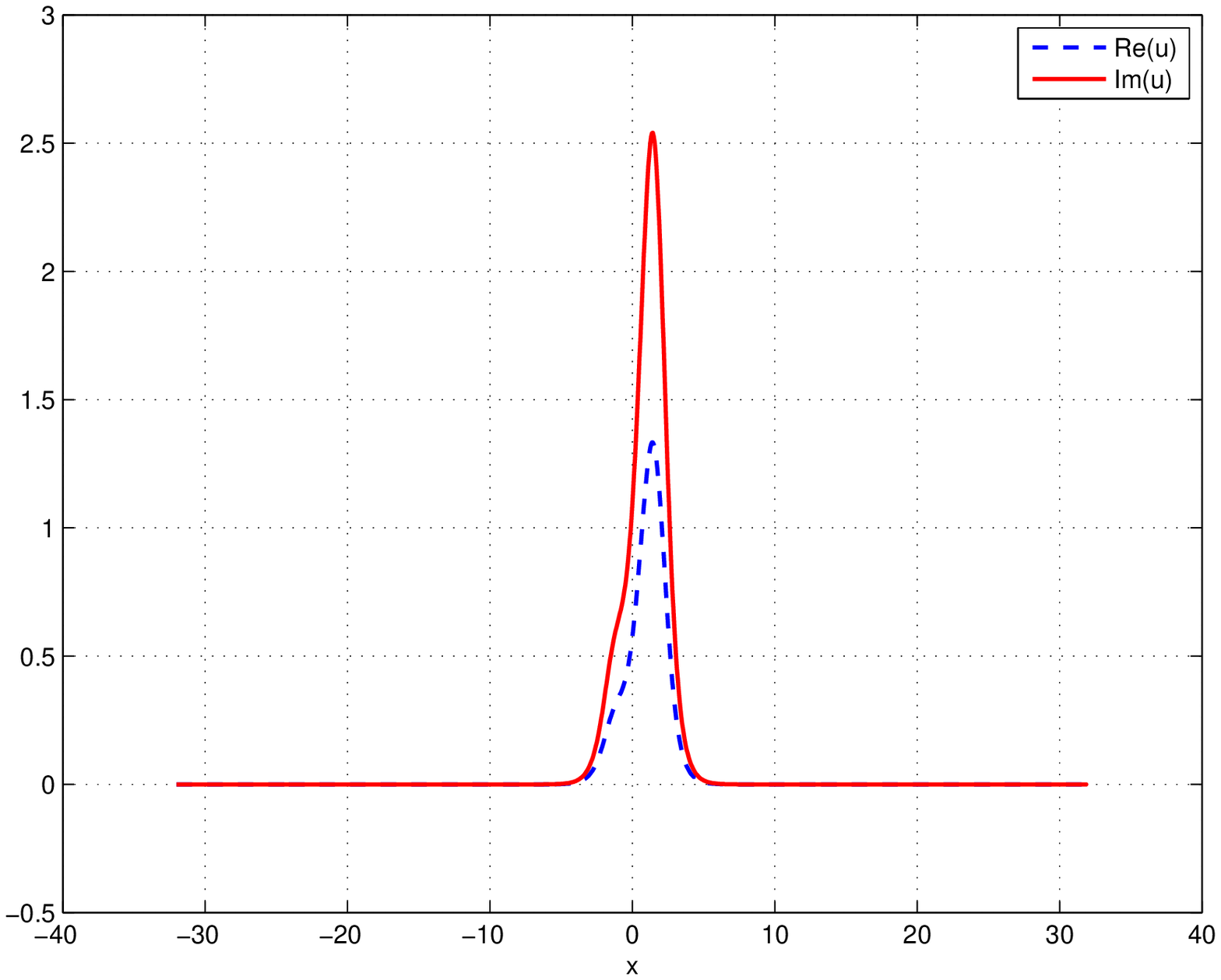}}
\subfigure[]{
\includegraphics[width=6.5cm]{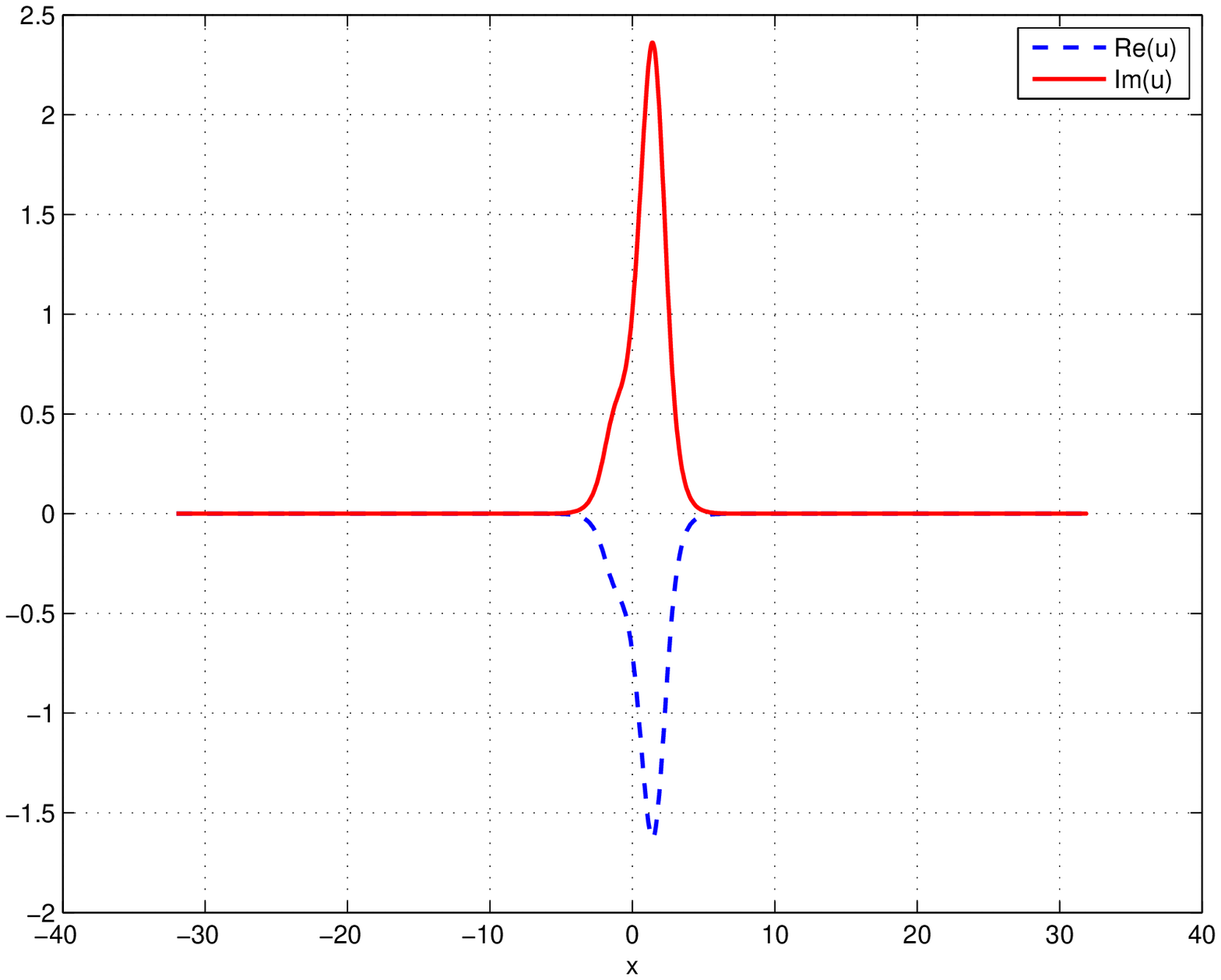}}
\subfigure[]{
\includegraphics[width=6.5cm]{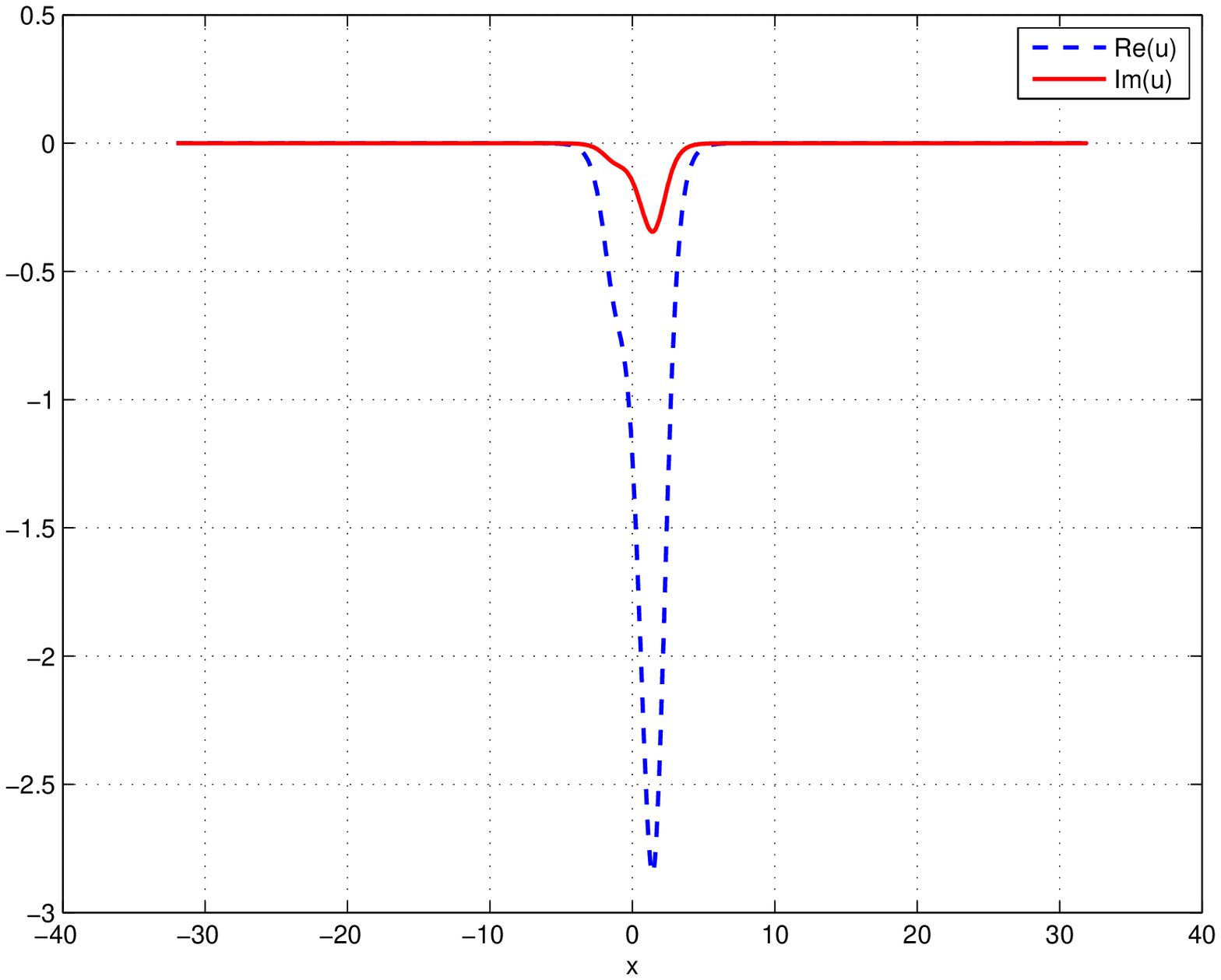}}
\subfigure[]{
\includegraphics[width=6.5cm]{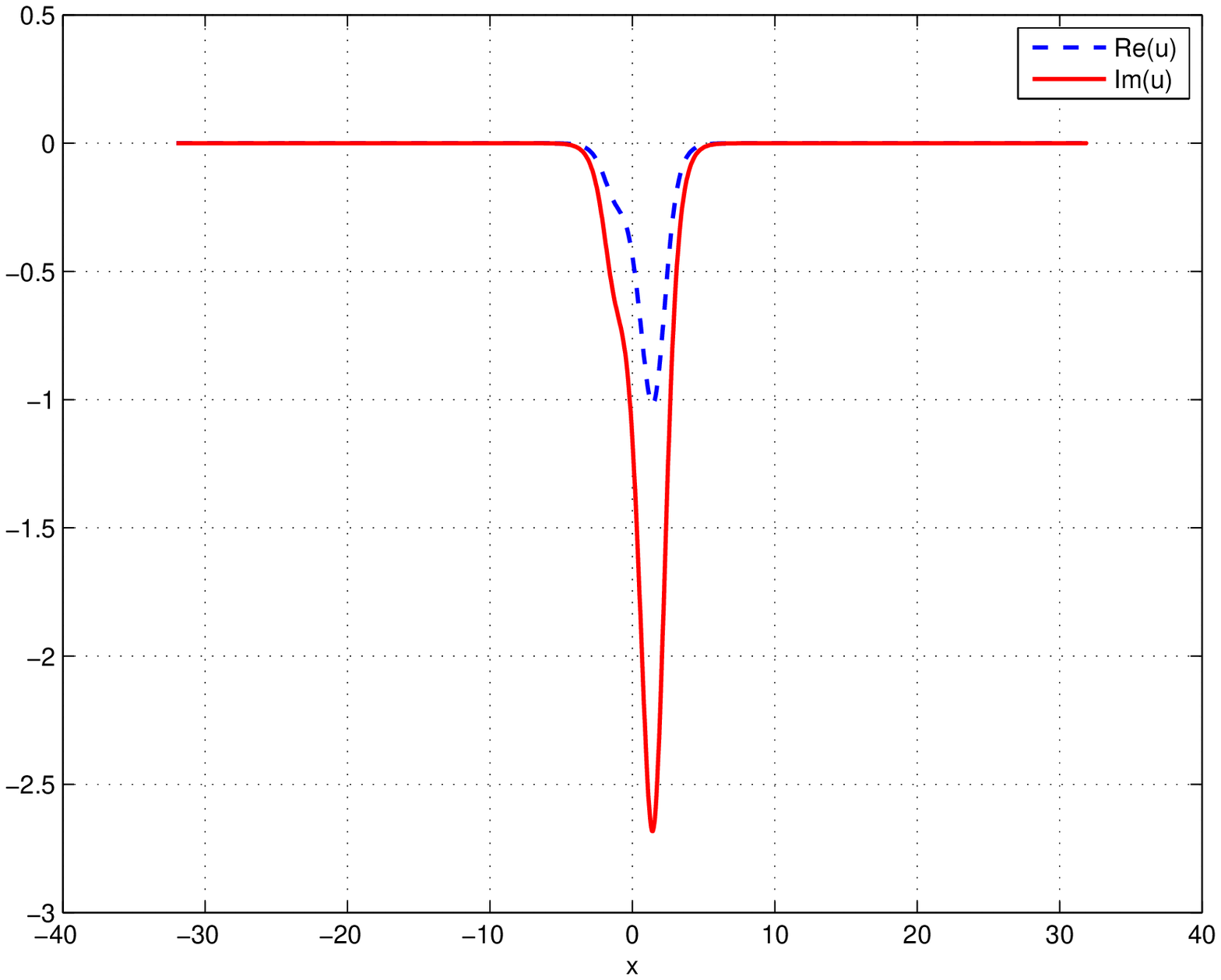}}
\caption{Evolution of the profile of Figure \ref{fexample_11}(b). Numerical solution (real part with solid line and imaginary part with dashed line) at times $t=50, 100, 150, 200$.} \label{fexample_12b}
\end{figure}

Finally, the accuracy of the computed profiles is checked in Figures \ref{fexample_12} (for $\mu=3.275$) and \ref{fexample_12b}  (for $\mu=3.289$). They correspond to considering the last iteration as initial condition of a time-stepping code for (\ref{gnls2}) and leaving the numerical solution to evolve. As in the previous example, the evolution of the real and imaginary parts is illustrated (with solid and dashed lines, respectively).

\begin{figure}[htbp]
\centering
\subfigure[]{
\includegraphics[width=6.5cm]{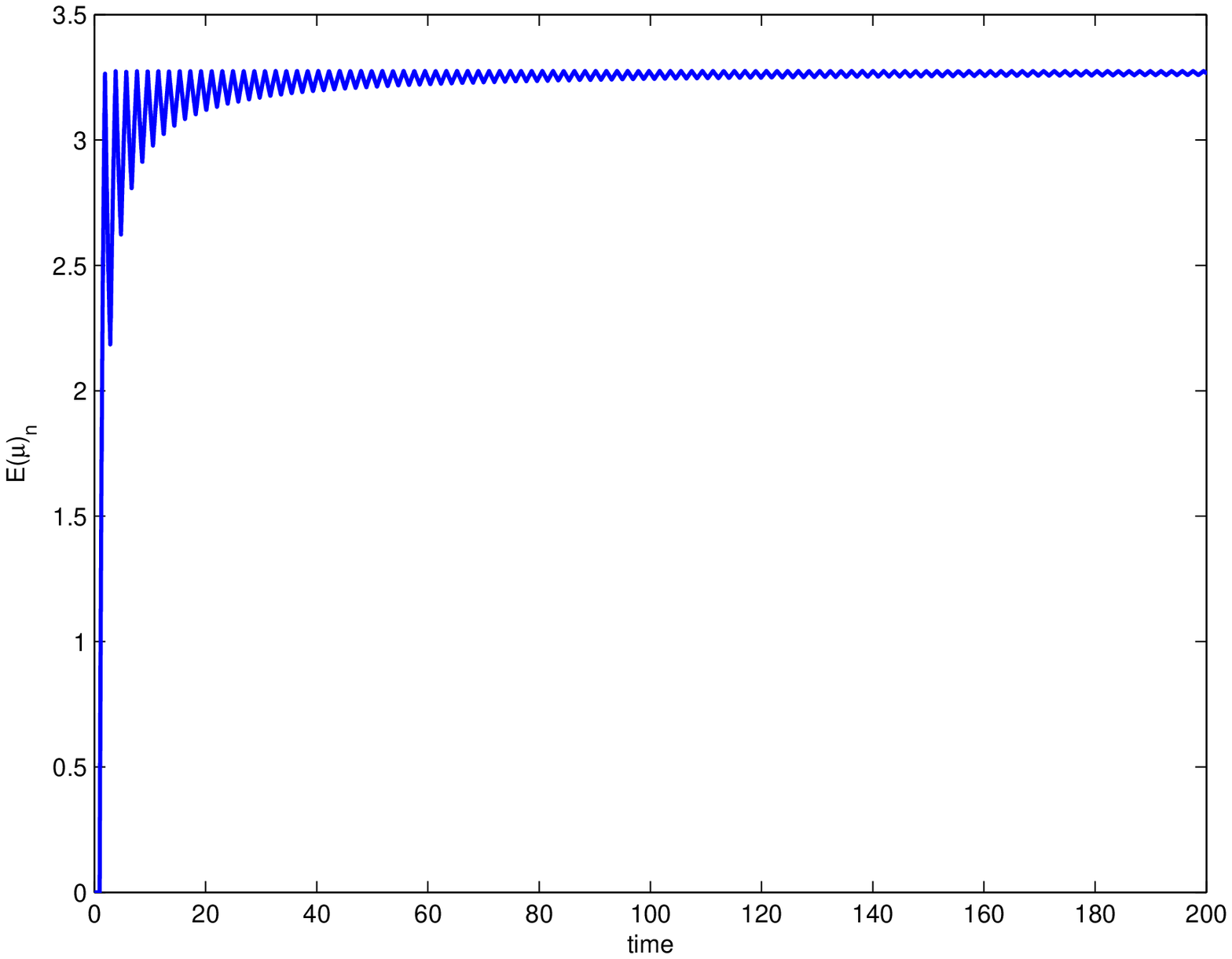}}
\subfigure[]{
\includegraphics[width=6.5cm]{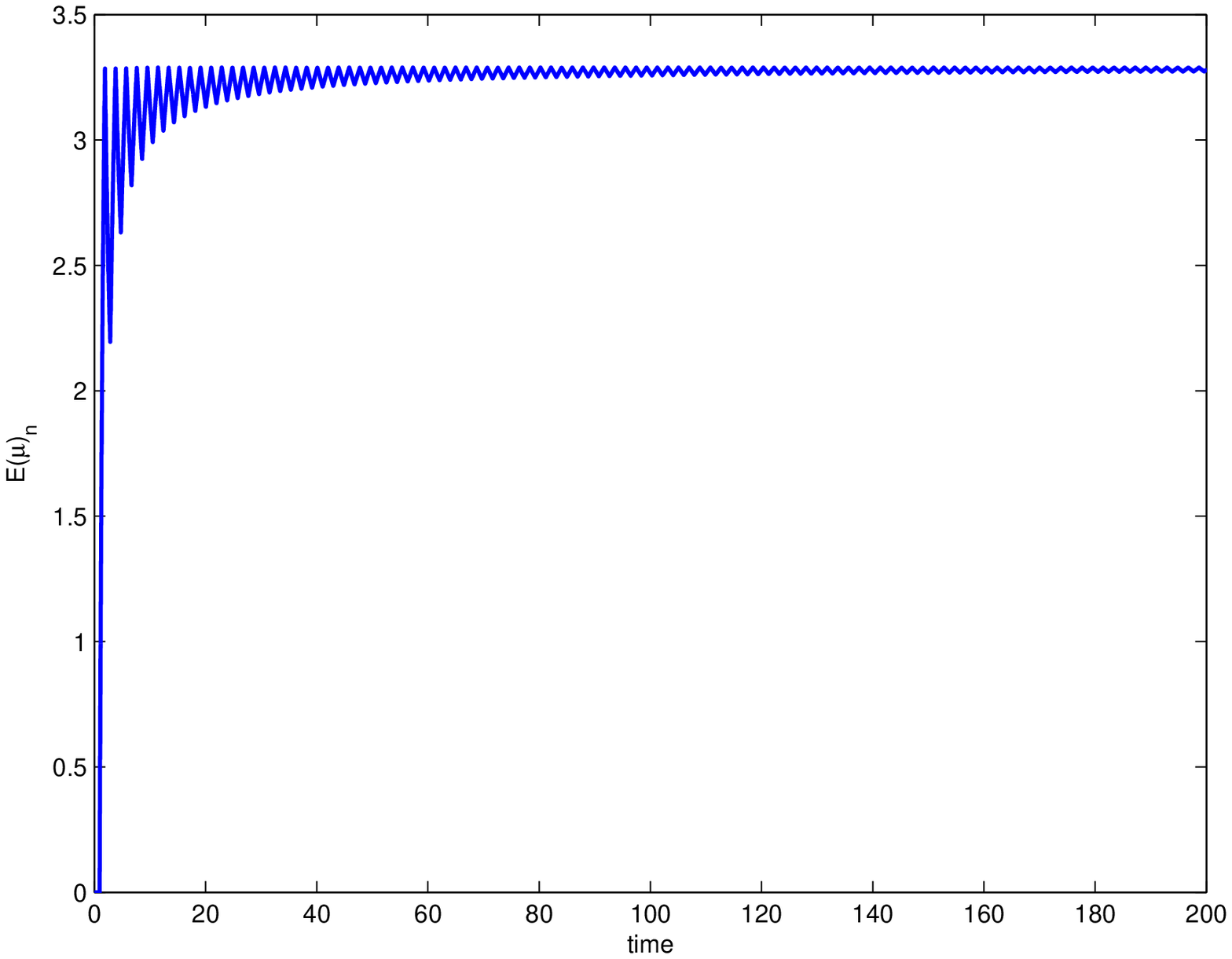}}
\caption{Evolution of the computed phase speed of the numerical solutions of Figure \ref{fexample_11}. (a) $\mu=3.275$; (b) $\mu=3.289$.} \label{fexample_13}
\end{figure}

The accuracy shown in Figures \ref{fexample_12} and \ref{fexample_12b}, is also confirmed, as in the previous example, by Figure \ref{fexample_13}. This displays the evolution of the computed angular velocity of the numerical approximation. This quantity is, in both cases ($\mu=3.275$ in Figure \ref{fexample_13}(a) and $\mu=3.289$ in Figure \ref{fexample_13}(b)) approaching the corresponding value of $\mu$.

\end{document}